\documentclass[10pt,a4paper]{article}

\usepackage{amsfonts,amsmath,mathrsfs,amssymb}
\usepackage{verbatim}
\usepackage{array}
\usepackage[pdftex]{graphics}
\usepackage[hang]{subfigure}

\usepackage{tikz}

\newtheorem{theorem}{Theorem}[section]
\newtheorem{lemma}[theorem]{Lemma}
\newtheorem{conj}[theorem]{Conjecture}
\newtheorem{defin}{Definition}[section]
\newtheorem{examp}{Example}[section]
\newtheorem{remark}{Remark}[section]
\newtheorem{cor}[theorem]{Corollary}

\definecolor{revision}{rgb}{1,0,0}

\begin{document}

\title {The orbit rigidity matrix of a symmetric framework}

\author{Bernd Schulze\footnote{Supported by the DFG Research Unit 565 `Polyhedral Surfaces'.}\\ Institute of Mathematics, MA 6-2\\ TU Berlin\\ Strasse des 17. Juni 136,\\ D-10623 Berlin, Germany\\
and\\ Walter Whiteley\footnote{Supported by a grant from NSERC (Canada).}\\
 Department of Mathematics and Statistics\\ York University\\ 4700 Keele Street\\ Toronto, ON M3J1P3, Canada\\ 
}

\maketitle

\begin{abstract}
A number of recent papers have studied when symmetry causes frameworks on a graph to become infinitesimally flexible, or stressed, and when it has no impact.  A number of other recent papers have studied special classes of frameworks on generically rigid graphs which are finite mechanisms.  Here we introduce a new tool, the orbit matrix, which connects these two areas and provides a matrix representation for   fully symmetric infinitesimal flexes, and fully symmetric stresses of symmetric frameworks.  The orbit matrix is a true analog of the standard rigidity matrix for general frameworks, and its analysis gives important insights into questions about the flexibility and rigidity of classes of symmetric frameworks, in all dimensions.

With this narrower focus on fully symmetric infinitesimal motions, comes the power to predict symmetry-preserving finite mechanisms - giving a simplified analysis which covers a wide range of the known mechanisms, and  generalizes the classes of known mechanisms.  This initial exploration of the properties of the orbit matrix also opens up a number of new questions and possible extensions of the previous results, including transfer of symmetry based results from Euclidean space to spherical, hyperbolic, and some other metrics with shared symmetry groups and underlying projective geometry.

\end{abstract}

\section{Introduction}
Over the last decade, a substantial theory on the interactions of symmetry and rigidity has been
{developed}
\cite{cfgsw,gsw,FGsymmax,KG1,KG2,owen,BS2,BS1,BS4}.   This includes descriptions of when symmetry changes generically rigid graphs into infinitesimally flexible frameworks, and when symmetry does not modify the behavior.  These analyses have used tools of representation theory  to analyze the stresses and motions of the symmetric realizations of a graph.
Some extensions have gone further to describe situations when the symmetry switches a graph into configurations with symmetry-preserving finite flexes \cite{FG4,BS4}.  These predictions of finite symmetric flexes turn out to focus on frameworks with fully symmetric infinitesimal flexes, and with fully symmetric self-stresses \cite{BS6}.

There is a companion, extensive literature on flexible frameworks built on generically rigid frameworks, starting with Bricard's flexible octahedra \cite{bricard,Stachel}, running though linkages such as Bottema's mechanism \cite{bottema} and other finitely flexible frameworks  \cite{consusp}
 and Connelly's flexible sphere \cite{concounter} to recent work on flexible cross-polytopes in 4-space \cite{Stachel4D}.  Some of this work has looked at creating analog examples of finite mechanisms in other metrics such as the spherical and hyperbolic space \cite{alex}.   In general, it is a difficult task to decide when a specific infinitesimal flex of a framework on a generically rigid graph extends to a finite flex.   However, on careful examination, many of these known examples have symmetries and infinitesimal flexes which preserve this symmetry \cite{BS6,tarnaiely}.    It is natural to seek  tools and connections that can simplify the creation and generation of such examples of finite mechanisms (linkages).

 In \cite{BS6}, one block of the block decomposition induced by the representations of the symmetry group was used to study the spaces of fully symmetric motions and fully symmetric self-stresses.  This analysis gave some initial results predicting finite flexes which remain fully symmetric throughout their path.   However, actual generation of this block in the decomposition required substantial machinery from representation theory, and the entries in the matrix
{were}
not transparent.

In this paper we present the {\em{orbit  matrix}} for a symmetric framework as an original, simplifying  tool for detecting this whole package of fully symmetric infinitesimal flexes, fully symmetric self-stresses, and predicting finite flexes for configurations which are generic within the symmetry.  In our proofs, we will actually show that the orbit matrix is equivalent to the matrix studied in \cite{BS6}, but the construction is transparent, and the entries in the matrix are explicitly derived.   For a symmetric graph, with symmetry group $S$, this orbit matrix has a set of columns for each orbit of vertices under the group action, and row for each orbit of edges under the group action.  We will give a detailed construction for this matrix in \S\ref{sec:orbitmatdef}, and show that the kernel is precisely the fully symmetric infinitesimal motions (\S\ref{sec:symmotstress}, \S\ref{sec:kernel}) and the row dependencies are exactly the fully symmetric self-stresses (\S\ref{sec:symmotstress}, \S\ref{sec:stress}).   This orbit matrix provides a powerful tool for investigating many aspects of the behavior of fully symmetric frameworks on the graph.

From the counts of the columns $c$, the rows $r$ and the dimension of the fully symmetric trivial infinitesimal motions $m$, we can give some immediate
sufficient
conditions for the presence of fully symmetric infinitesimal flexes (see \S\ref{sec:finite}).  Moreover, at configurations in which the representative vertices for the orbits are chosen `generically', the presence of a fully symmetric infinitesimal flex is a guarantee of a symmetry-preserving finite flex.   With these tools
{for counting under symmetry,}
 we have direct
predictions of flexible frameworks which capture many of the classical examples, including two types of flexible octahedra, the Bottema mechanism, and the flexible cross-polytopes.   A striking example of a general class covered by this analysis is the following:

{\bf Theorem~\ref{3dhalfturn}} {\it Given a graph  which is generically isostatic in $3$-space, and a framework  on the graph realized in $3$-space as generically as possible with 2-fold symmetry with no vertices or edges fixed by the rotation, the framework has a finite flex preserving the symmetry.}

Because this orbit matrix is a powerful symmetry adapted analog of the standard rigidity matrix, many of the questions from standard rigidity have extensions for fully symmetric stresses and motions.  Some of these questions and possible extensions are presented in \S\ref{sec:further}, along with brief discussions of the potential for
{symmetry adapted extensions of}
the techniques and results.  As one example, it is natural to seek analogs of Laman's Theorem to characterize necessary and sufficient conditions for the orbit matrix of a graph and symmetry group to be independent with maximal rank.  As a second example,
because key portions of the point group symmetries and the corresponding counting for this orbit matrix can be transferred to other metrics (such as the spherical, hyperbolic and Minkowski spaces), the methods developed here provide a uniform construction of mechanisms such as the Bricard octahedron, the flexible cross-polytope, and the Bottema mechanism and its generalizations across multiple metrics.

As a final comment, key results on the global rigidity of generic frameworks depend on self-stresses and the equivalence of finite flexes and infinitesimal flexes for generic frameworks.  We now have fully-symmetric versions of these tools, and can extract some analogs of the global rigidity results for symmetry-generic frameworks, both as conditions under which they are globally rigid within the class of fully-symmetric frameworks, and when they are globally rigid within the class of all frameworks.  Still there are additional conjectures and new results to be explored in this area.

We hope that this paper serves as an invitation
{for the reader}
to join in the further explorations of these many levels of interactions of symmetry, rigidity, and flexibility.

\section{Rigidity theoretic definitions and preliminaries}

All graphs considered in this paper are finite graphs without loops or multiple edges. The \emph{vertex set} of a  \emph{graph} $G$ is denoted by $V(G)$ and the \emph{edge set} of $G$ is denoted by $E(G)$.\\
\indent A \emph{framework}  in $\mathbb{R}^{d}$ is a pair $(G,p)$, where $G$ is a graph and  $p:V(G)\to \mathbb{R}^d$ is a map such that $p(u) \neq p(v)$ for all $\{u,v\} \in E(G)$.  We also say that $(G,p)$ is a $d$-dimensional \emph{realization} of the \emph{underlying graph} $G$ \cite{gss, W1}.
For $v\in V(G)$, we say that $p(v)$ is the \emph{joint} of $(G,p)$ corresponding to $v$, and for $e\in E(G)$, we say that $p(e)$ is the \emph{bar} of $(G,p)$ corresponding to $e$.\\\indent
For a framework $(G,p)$ whose underlying graph $G$ has the vertex set $V(G)=\{1,\ldots ,n\}$, we will frequently denote the vector $p(i)$ by $p_{i}$ for each $i$. The $k^{th}$ component of a vector $x$ is denoted by $(x)_{k}$. It is often useful to identify $p$ with a vector in $\mathbb{R}^{dn}$ by using the order on $V(G)$. In this case we also refer to $p$ as a \emph{configuration} of $n$ points in $\mathbb{R}^{d}$. Throughout this paper, we do not differentiate between an abstract vector and its coordinate column vector relative to the canonical basis.\\\indent
A framework $(G,p)$ in $\mathbb{R}^d$ with $V(G)=\{1, \ldots, n\}$ is \emph{flexible} if there exists a continuous path, called a \emph{finite flex} or \emph{mechanism}, $p(t):[0,1]\to \mathbb{R}^{dn}$  such that \begin{itemize}
\item[(i)] $p(0)=p$;
\item[(ii)]  $\|p(t)_i-p(t)_j\|=\|p_i-p_j\|$ for  all $0\leq t\leq 1$ and all $\{i,j\}\in E(G)$;
\item[(iii)]  $\|p(t)_k-p(t)_l\|\neq\|p_k-p_l\|$ for  all $0< t\leq 1$ and some pair $\{k,l\}$ of vertices of $G$. \end{itemize}
 Otherwise $(G,p)$ is said to be \emph{rigid}. For some alternate equivalent definitions of a rigid and flexible framework see \cite{asiroth, RW1}, for example.\\\indent
An \emph{infinitesimal motion} of a framework $(G,p)$ in $\mathbb{R}^d$ with $V(G)=\{1, \ldots, n\}$ is a function $u: V(G)\to \mathbb{R}^{d}$ such that
\begin{equation}
\label{infinmotioneq}
(p_i-p_j)^T (u_i-u_j)=0 \quad\textrm{ for all } \{i,j\} \in E(G)\textrm{,}\end{equation}
where $u_i$ denotes the column vector $u(i)$ for each $i$.
\\\indent An infinitesimal motion $u$ of $(G,p)$ is an \emph{infinitesimal rigid motion} (or \emph{trivial infinitesimal motion}) if there exists a skew-symmetric matrix $S$ (a rotation) and a vector $t$ (a translation) such that $u(v)=Sp(v)+t$ for all $v\in V(G)$. Otherwise $u$ is an \emph{infinitesimal flex} (or \emph{non-trivial infinitesimal motion}) of $(G,p)$.\\\indent
$(G,p)$ is \emph{infinitesimally rigid} if every infinitesimal motion of $(G,p)$ is an infinitesimal rigid motion. Otherwise $(G,p)$ is said to be \emph{infinitesimally flexible} \cite{gss, W1}.\\\indent
The
\emph{rigidity matrix} of $(G,p)$ is the $|E(G)| \times dn$ matrix
\begin{displaymath} \mathbf{R}(G,p)=\bordermatrix{& & & & i & & & & j & & & \cr & & & &  & & \vdots & &  & & &
\cr \{i,j\} & 0 & \ldots &  0 & (p_{i}-p_{j})^T & 0 & \ldots & 0 & (p_{j}-p_{i})^T &  0 &  \ldots&  0 \cr & & & &  & & \vdots & &  & & &
}
\textrm{,}\end{displaymath}
that is, for each edge $\{i,j\}\in E(G)$, $\mathbf{R}(G,p)$ has the row with
$(p_{i}-p_{j})_{1},\ldots,(p_{i}-p_{j})_{d}$ in the columns $d(i-1)+1,\ldots,di$, $(p_{j}-p_{i})_{1},\ldots,(p_{j}-p_{i})_{d}$ in
the columns $d(j-1)+1,\ldots,dj$, and $0$ elsewhere \cite{gss, W1}.\\\indent
Note that if we identify an infinitesimal motion $u$ of $(G,p)$ with a column vector in $\mathbb{R}^{dn}$ (by using the order on $V(G)$), then the equations in (\ref{infinmotioneq}) can be written as $\mathbf{R}(G,p)u=0$. So, the kernel of the rigidity matrix $\mathbf{R}(G,p)$ is the space of all infinitesimal motions of $(G,p)$. It is well known that a framework $(G,p)$ in $\mathbb{R}^d$ is infinitesimally rigid if and only if either the rank of its associated rigidity matrix $\mathbf{R}(G,p)$ is precisely $dn-\binom{d+1}{2}$, or $G$ is a complete graph $K_n$ and the points $p_i$, $i=1,\ldots, n$, are affinely independent \cite{asiroth}.\\\indent
While an infinitesimally rigid framework is always rigid, the converse does not hold in general. Asimov and Roth, however, showed that for `generic' configurations, infinitesimal rigidity and rigidity are in fact equivalent \cite{asiroth}.\\\indent
A \emph{self-stress} of a framework $(G,p)$ with $V(G)=\{1, \ldots, n\}$ is a function  $\omega:E(G)\to \mathbb{R}$ such that at each joint $p_i$ of $(G,p)$ we have
\begin{displaymath}
\sum_{j :\{i,j\}\in E(G)}\omega_{ij}(p_{i}-p_{j})=0 \textrm{,}
\end{displaymath}
where $\omega_{ij}$ denotes $\omega(\{i,j\})$ for all $\{i,j\}\in E(G)$. Note that if we identify a self-stress $\omega$ with a column vector in $\mathbb{R}^{|E(G)|}$ (by using the order on $E(G)$), then we have $\omega^{T}\mathbf{R}(G,p)=0$.
In structural engineering, the self-stresses are also called \emph{equilibrium stresses} as they record  tensions and compressions in the bars balancing at each vertex.
\\\indent If $(G,p)$ has a non-zero self-stress, then $(G,p)$ is said to be \emph{dependent} (since in this case there exists a linear dependency among the row vectors of $\mathbf{R}(G,p)$). Otherwise, $(G,p)$ is said to be \emph{independent}.  A framework which is both independent and infinitesimally rigid is called \emph{isostatic} \cite{CW, W4, W1}.

\section{Symmetry in frameworks}
\label{sec:sym}

\begin{figure}[htp]
\begin{center}
\begin{tikzpicture}[very thick,scale=1]
\tikzstyle{every node}=[circle, draw=black, fill=white, inner sep=0pt, minimum width=5pt];
        \path (1.5,-1) node (p2) [label = below left: $p_{2}$] {} ;
       \path (0.6,0.5) node (p1) [label = left: $p_{1}$] {} ;
   \path (2.4,0.5) node (p3) [label = right: $p_{3}$] {} ;
   \path (1.5,1.1) node (p4) [label = above left: $p_{4}$] {} ;
   \draw (p1) -- (p4);
      \draw (p3) -- (p4);
     \draw (p2) -- (p3);
      \draw (p2) -- (p1);
        \draw [dashed, thin] (1.5,-1.6) -- (1.5,1.6);
      \node [draw=white, fill=white] (a) at (1.5,-2.1) {(a)};
    \end{tikzpicture}
    \hspace{2cm}
        \begin{tikzpicture}[very thick,scale=1]
\tikzstyle{every node}=[circle, draw=black, fill=white, inner sep=0pt, minimum width=5pt];
    \path (-0.7,0.8) node (p1) [label = above left: $p_1$] {} ;
    \path (0.7,0.8) node (p4) [label = above right: $p_4$] {} ;
    \path (-1.6,-0.8) node (p2) [label = below left: $p_2$] {} ;
     \path (1.6,-0.8) node (p3) [label = below right: $p_3$] {} ;
      \draw (p1) -- (p4);
    \draw (p1) -- (p2);
    \draw (p3) -- (p4);
    \draw (p2) -- (p3);
     \draw [dashed, thin] (0,-1.6) -- (0,1.6);
      \node [draw=white, fill=white] (b) at (0,-2.1) {(b)};
        \end{tikzpicture}
\end{center}
\vspace{-0.3cm}
\caption{$2$-dimensional realizations of $K_{2,2}$ in $\mathscr{R}_{(K_{2,2},\mathcal{C}_s)}$ of different types: the framework in (a) is of type
$\Phi_{a}$, where $\Phi_{a}: \mathcal{C}_{s} \to \textrm{Aut}(K_{2,2})$ is the homomorphism defined by $\Phi_{a}(s)=
(1 \, 3)(2)(4)$ and the framework in (b) is of type
$\Phi_{b}$, where $\Phi_{b}: \mathcal{C}_{s} \to \textrm{Aut}(K_{2,2})$ is the homomorphism defined by $\Phi_{b}(s)=
(1 \, 4)(2\, 3)$.}
\label{K33types}
\end{figure}
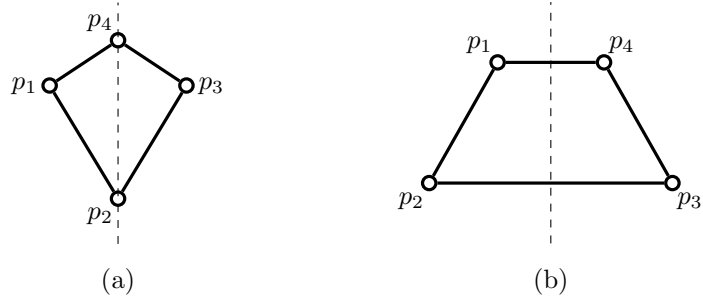

Let $G$ be a graph with $V(G)=\{1,\ldots,n\}$, and let $\textrm{Aut}(G)$ denote the automorphism group of $G$. A \emph{symmetry operation} of a framework $(G,p)$ in $\mathbb{R}^{d}$ is an isometry $x$ of $\mathbb{R}^{d}$ such that for some $\alpha\in \textrm{Aut}(G)$, we have
$x(p_i)=p_{\alpha(i)}$ for all $i\in V(G)$ \cite{Hall, BS2, BS4, BS1}.  \\\indent
The set of all symmetry operations of a framework $(G,p)$ forms a group under composition, called the \emph{point group} of $(G,p)$ \cite{bishop, Hall, BS4, BS1}. Since translating a framework does not change its rigidity properties, we may assume wlog that the point group of any framework in this paper is a \emph{symmetry group}, i.e., a subgroup of the orthogonal group $O(\mathbb{R}^{d})$ \cite{BS2, BS4, BS1}.\\\indent
We use the Schoenflies notation for the symmetry operations and symmetry groups considered in this paper, as this is one of the standard notations in the literature about symmetric structures (see \cite{bishop, cfgsw, FGsymmax, FG4, Hall, KG1, KG2, BS2, BS4, BS1}, for example). In this notation, the identity transformation is denoted by $Id$, a rotation  about a $(d-2)$-dimensional subspace of $\mathbb{R}^d$ by an angle of $\frac{2\pi}{m}$ is denoted by $C_m$, and a reflection in a $(d-1)$-dimensional subspace of $\mathbb{R}^d$ is denoted by $s$.\\\indent
While the general results of this paper apply to all symmetry groups, we will only analyze examples with four types of groups.
In the Schoenflies notation, they are denoted by $\mathcal{C}_{s}$, $\mathcal{C}_{m}$, $\mathcal{C}_{mv}$, and $\mathcal{C}_{mh}$. For any dimension $d$, $\mathcal{C}_{s}$ is a symmetry group consisting of the identity $Id$ and a single reflection $s$, and $\mathcal{C}_{m}$ is a cyclic group generated by a rotation $C_m$. The only other possible type of symmetry
  group in dimension 2 is the group $\mathcal{C}_{mv}$ which is a dihedral group generated by a pair $\{C_m,s\}$. In dimension $d>2$, $\mathcal{C}_{mv}$ denotes any symmetry group that is generated by a rotation $C_m$ and a reflection $s$ whose corresponding mirror contains the rotational axis of $C_m$, whereas a symmetry group $\mathcal{C}_{mh}$ is generated by a rotation $C_m$ and the reflection $s$ whose corresponding mirror is perpendicular to the $C_m$-axis. For further information about the Schoenflies notation we refer the reader to \cite{bishop, Hall, BS4}.
\\\indent
Given a symmetry group $S$ in dimension $d$ and a graph $G$, we let $\mathscr{R}_{(G,S)}$ denote the set of all $d$-dimensional realizations of $G$ whose point group is either equal to $S$ or contains $S$ as a subgroup \cite{BS2, BS4, BS1}. In other words, the set $\mathscr{R}_{(G,S)}$ consists of all realizations $(G,p)$ of $G$ for which there exists a map $\Phi:S\to \textrm{Aut}(G)$ so that
\begin{equation}\label{class} x\big(p_i\big)=p_{\Phi(x)(i)}\textrm{ for all } i\in V(G)\textrm{ and all } x\in S\textrm{.}\end{equation}
A framework $(G,p)\in \mathscr{R}_{(G,S)}$ satisfying the equations in (\ref{class}) for the map $\Phi:S\to \textrm{Aut}(G)$ is said to be \emph{of type $\Phi$}, and the set of all realizations in $\mathscr{R}_{(G,S)}$ which are of type $\Phi$ is denoted by $\mathscr{R}_{(G,S,\Phi)}$ (see again \cite{BS2, BS4, BS1} and Figure \ref{K33types}).
\\\indent It is shown in \cite{BS4, BS1} that if the map $p$ of a framework  $(G,p)\in \mathscr{R}_{(G,S)}$ is injective, then $(G,p)$ is of a unique type $\Phi$ and $\Phi$ is necessarily also a homomorphism. For simplicity, we therefore assume that the map $p$ of any framework $(G,p)$ considered in this paper is injective (i.e., $p_i\neq p_j$ if $i\neq j$). In particular, this allows us (with a slight abuse of notation) to use the terms $p_{x(i)}$ and  $p_{\Phi(x)(i)}$ interchangeably, where $i\in V(G)$ and $x\in S$. In general, if the type $\Phi$ is clear from the context, we often simply write $x(i)$ instead of $\Phi(x)(i)$.\\\indent
Let $(G,p)\in \mathscr{R}_{(G,S,\Phi)}$ and let  $x$ be a symmetry operation in $S$. Then the joint $p_i$ of $(G,p)$ is said to be \emph{fixed} by $x$ if $x(p_i)=p_i$ (or equivalently, $x(i)=i$),\\\indent Let the \emph{symmetry element} corresponding to  $x$ be the linear subspace $F_{x}$ of $\mathbb{R}^{d}$ which consists of all points $a\in\mathbb{R}^{d}$ with $x(a)=a$. Then the joint $p_i$ of any framework $(G,p)$ in $\mathscr{R}_{(G,S,\Phi)}$ must lie in the linear subspace \begin{displaymath}U(p_i)=\bigcap_{x\in S: x(p_i)=p_i}F_x\textrm{.}\end{displaymath} Note that $U(p_i)\neq \emptyset$, because the origin is fixed by every symmetry operation in the symmetry group $S$.

\begin{examp} The joint $p_1$ of the framework $(K_{2,2},p)\in \mathscr{R}_{(K_{2,2},\mathcal{C}_s,\Phi_a)}$ depicted in Figure \ref{K33types} (a) is fixed by the identity $Id\in \mathcal{C}_s$, but not by the reflection $s\in \mathcal{C}_s$, so that $U(p_1)=F_{Id}=\mathbb{R}^d$. The joint $p_2$ of $(K_{2,2},p)$, however, is fixed by both the identity $Id$ and the reflection $s$ in  $\mathcal{C}_s$, so that $U(p_2)=F_{Id}\cap F_s=F_s$. In other words, $U(p_2)$ is the mirror line corresponding to $s$.
\end{examp}

Note that if we choose a set of representatives $\mathscr{O}_{V(G)}=\{1,\ldots, k\}$ for the orbits  $S(i)=\{\Phi(x)(i)|\,x\in S\}$ of vertices of $G$, then the positions of \emph{all} joints of  $(G,p)\in \mathscr{R}_{(G,S,\Phi)}$ are uniquely determined by the positions of the joints $p_1,\ldots,p_k$ and the symmetry constraints imposed by $S$ and $\Phi$. Thus, any framework in $\mathscr{R}_{(G,S,\Phi)}$ may be constructed by first choosing positions $p_i\in U(p_{i})$ for each $i=1,\ldots, k$, and then letting $S$ and $\Phi$ determine the positions of the remaining joints. In particular, by placing the vertices of $\mathscr{O}_{V(G)}$ into `generic' positions within their associated subspaces $U(p_i)$ we obtain an \emph{$(S,\Phi)$-generic} realization of $G$ (i.e., a realization of $G$ that is as `generic' as possible within the set $\mathscr{R}_{(G,S,\Phi)}$)  in this way. For a precise definition of $(S,\Phi)$-generic, and further information about  $(S,\Phi)$-generic frameworks, we refer the reader to \cite{BS4, BS1}.

\section{Fully symmetric motions and self-stresses}
\label{sec:symmotstress}

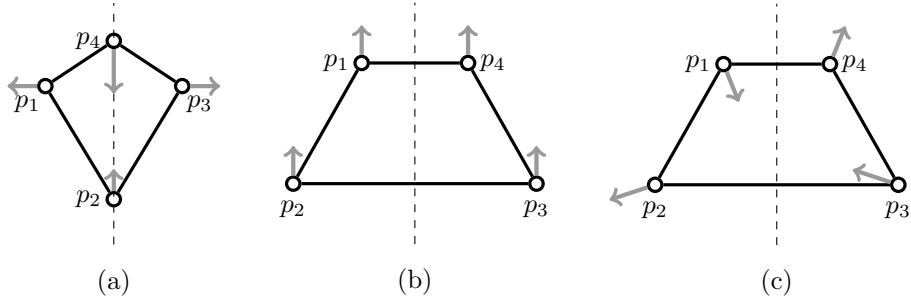
\begin{figure}[htp]
\begin{center}
\begin{tikzpicture}[very thick,scale=1]
\tikzstyle{every node}=[circle, draw=black, fill=white, inner sep=0pt, minimum width=5pt];
        \path (1.5,-1) node (p2)[label = left: $p_{2}$] {} ;
       \path (0.6,0.5) node (p1)  [label = below left: $p_{1}$]{} ;
   \path (2.4,0.5) node (p3) [label = below right: $p_{3}$] {} ;
   \path (1.5,1.1) node (p4) [label = left: $p_{4}$] {} ;
   \draw (p1) -- (p4);
      \draw (p3) -- (p4);
     \draw (p2) -- (p3);
      \draw (p2) -- (p1);
             \draw [dashed, thin] (1.5,-1.6) -- (1.5,1.6);
              \draw [ultra thick, ->, black!40!white] (p1) -- (0.1,0.5);
      \draw [ultra thick, ->, black!40!white] (p3) -- (2.9,0.5);
      \draw [ultra thick, ->, black!40!white] (p4) -- (1.5,0.4);
      \draw [ultra thick, ->, black!40!white] (p2) -- (1.5,-0.6);
      \node [draw=white, fill=white] (a) at (1.5,-2.1) {(a)};
    \end{tikzpicture}
    \hspace{0.5cm}
            \begin{tikzpicture}[very thick,scale=1]
\tikzstyle{every node}=[circle, draw=black, fill=white, inner sep=0pt, minimum width=5pt];
    \path (-0.7,0.8) node (p1) [label = left: $p_{1}$] {} ;
    \path (0.7,0.8) node (p4) [label = right: $p_{4}$]{} ;
    \path (-1.6,-0.8) node (p2) [label = below: $p_{2}$] {} ;
     \path (1.6,-0.8) node (p3) [label = below: $p_{3}$] {} ;
      \draw (p1) -- (p4);
    \draw (p1) -- (p2);
    \draw (p3) -- (p4);
    \draw (p2) -- (p3);
     \draw [dashed, thin] (0,-1.6) -- (0,1.6);
     \draw [ultra thick, ->, black!40!white] (p1) -- (-0.7,1.3);
      \draw [ultra thick, ->, black!40!white] (p4) -- (0.7,1.3);
      \draw [ultra thick, ->, black!40!white] (p2) -- (-1.6,-0.3);
      \draw [ultra thick, ->, black!40!white] (p3) -- (1.6,-0.3);
      \node [draw=white, fill=white] (b) at (0,-2.1) {(b)};
        \end{tikzpicture}
        \hspace{0.5cm}
        \begin{tikzpicture}[very thick,scale=1]
\tikzstyle{every node}=[circle, draw=black, fill=white, inner sep=0pt, minimum width=5pt];
    \path (-0.7,0.8) node (p1) [label = left: $p_{1}$]  {} ;
    \path (0.7,0.8) node (p4)[label = right: $p_{4}$] {} ;
    \path (-1.6,-0.8) node [label = below: $p_{2}$](p2)  {} ;
     \path (1.6,-0.8) node [label = below: $p_{3}$](p3)  {} ;
      \draw (p1) -- (p4);
    \draw (p1) -- (p2);
    \draw (p3) -- (p4);
    \draw (p2) -- (p3);
     \draw [dashed, thin] (0,-1.6) -- (0,1.6);
     \draw [ultra thick, ->, black!40!white] (p1) -- (-0.5,0.3);
      \draw [ultra thick, ->, black!40!white] (p4) -- (0.9,1.3);
      \draw [ultra thick, ->, black!40!white] (p2) -- (-2.2,-1);
      \draw [ultra thick, ->, black!40!white] (p3) -- (1,-0.6);
      \node [draw=white, fill=white] (b) at (0,-2.1) {(c)};
        \end{tikzpicture}
\end{center}
\vspace{-0.3cm}
\caption{Infinitesimal motions of frameworks in the plane: (a) a fully $(\mathcal{C}_s,\Phi_a)$-symmetric infinitesimal flex of $(K_{2,2},p) \in \mathscr{R}_{(K_{2,2},\mathcal{C}_s, \Phi_a)}$; (b) a fully $(\mathcal{C}_s,\Phi_b)$-symmetric infinitesimal rigid motion of $(K_{2,2},p) \in \mathscr{R}_{(K_{2,2},\mathcal{C}_s, \Phi_b)}$; (c) an infinitesimal flex of $(K_{2,2},p) \in \mathscr{R}_{(K_{2,2},\mathcal{C}_s, \Phi_b)}$ which is not fully $(\mathcal{C}_s,\Phi_b)$-symmetric.}
\label{fulsym}
\end{figure}

An infinitesimal motion $u$ of a framework $(G,p)\in \mathscr{R}_{(G,S,\Phi)}$  is \emph{fully $(S,\Phi)$-symmetric} if \begin{equation}\label{fulsymmot} x\big(u_i\big)=u_{\Phi(x)(i)}\textrm{ for all } i\in V(G)\textrm{ and all } x\in S\textrm{,}\end{equation} i.e., if $u$ is unchanged under all symmetry operations in $S$ (see also Figure \ref{fulsym}(a) and (b)).

Note that it follows immediately from (\ref{fulsymmot}) that if $u$ is a fully $(S,\Phi)$-symmetric infinitesimal motion of $(G,p)$, then $u_i$ is an element of $U(p_i)$ for each $i$. Moreover, $u$ is uniquely determined by the velocity vectors $u_1,\ldots, u_k$ whenever  $\mathscr{O}_{V(G)}=\{1,\ldots, k\}$ is a set of representatives for the vertex orbits  $S(i)=\{\Phi(x)(i)|\,x\in S\}$ of $G$.

\begin{examp}\label{quadcs1} Consider the framework shown in Figure \ref{fulsym}(a). With $p_1^T=(a,b), p_2^T=(0,c), p_3^T= (-a,b),$ and $p_4^T=(0,d)$ the rigidity matrix of $(K_{2,2},p)$ has the form
\begin{displaymath}\bordermatrix{
                &1&2&3=s(1) &4\cr
                \{1,2\}&(a, b-c) &  (-a,c-b)  &0\ 0 & 0 \ 0\cr
                 \{1,4\}& (a,b-d) & 0 \ 0  & 0 \ 0  & (-a,d-b)\cr
                s( \{1,2\})& 0 \ 0  & (a,c-b)   & (-a,b-c) & 0 \ 0  \cr
                s( \{1,4\})&  0 \ 0 &   0 \ 0       & (-a,b-d) & (a,d-b)}\end{displaymath}
This matrix has rank 4, and hence leaves a space of $8-4 = 4$ infinitesimal motions. Thus, there exists a $1$-dimensional space of  infinitesimal flexes of $(K_{2,2},p)$ spanned by $u=\left(\begin{array} {cccccccc}-1 & 0 & 0 & \frac{a}{c-b} & 1 & 0 &  0 & \frac{a}{d-b}\end{array}\right)^T$. This infinitesimal flex is clearly fully $(\mathcal{C}_s,\Phi_a)$-symmetric.
\end{examp}

\begin{examp}\label{quadcs2}
The rigidity matrix of the framework  $(K_{2,2},p)$ shown in Figure \ref{fulsym}(b,c) with $p_1^T=(a,b), p_2^T=(c,d), p_3^T= (-c,d),$ and $p_4^T=(-a,b)$ has the form
\begin{displaymath}\bordermatrix{
                &1&2&3=s(2) &4=s(1)\cr
                \{1,2\}&(a-c, b-d) &  (c-a,d-b)  &0\ 0 & 0 \ 0\cr
                 \{1,4\}& (2a,0) & 0 \ 0  & 0 \ 0  & (-2a,0)\cr
                 \{2,3\}& 0 \ 0  & (2c,0)   & (-2c,0) & 0 \ 0  \cr
                s( \{1,2\})&  0 \ 0 &   0 \ 0       & (a-c,d-b) & (c-a,b-d)}\end{displaymath}
This matrix has again rank 4, and leaves a space of $8-4 = 4$ infinitesimal motions. The $1$-dimensional space of  infinitesimal flexes of $(K_{2,2},p)$ is spanned by $u=\left(\begin{array} {cccccccc}1 & -1 & -1 & \frac{2(c-a)+b-d}{d-b} & -1 & -\frac{2(c-a)+b-d}{d-b} &  1 & 1\end{array}\right)^T$. This infinitesimal flex is clearly \emph{not} fully $(\mathcal{C}_s,\Phi_b)$-symmetric.
\end{examp}

\begin{figure}[htp]
\begin{center}
 \begin{tikzpicture}[very thick,scale=1]
\tikzstyle{every node}=[circle, draw=black, fill=white, inner sep=0pt, minimum width=5pt];
    \path (-0.7,0.8) node (p1) {} ;
    \path (0.7,0.8) node (p4)  {} ;
    \path (-1.6,-0.8) node (p2)  {} ;
     \path (1.6,-0.8) node (p3)  {} ;
      \draw (p1) -- (p4);
    \draw (p1) -- (p2);
    \draw (p3) -- (p4);
    \draw (p2) -- (p3);
     \draw (p2) -- (p4);
    \draw (p1) -- (p3);
     \draw [dashed, thin] (0,-1.6) -- (0,1.6);
     \node [draw=white, fill=white] (a) at (-1.4,0) {$\alpha$};
     \node [draw=white, fill=white] (a) at (1.4,0) {$\alpha$};
     \node [draw=white, fill=white] (a) at (0,1.1) {$\beta$};
     \node [draw=white, fill=white] (a) at (0,-1.1) {$\delta$};
     \node [draw=white, fill=white] (a) at (0.5,-0.3) {$\gamma$};
     \node [draw=white, fill=white] (a) at (-0.5,-0.3) {$\gamma$};
      \node [draw=white, fill=white] (b) at (0,-2.1) {(a)};
        \end{tikzpicture}
        \hspace{2cm}
        \begin{tikzpicture}[very thick,scale=1]
\tikzstyle{every node}=[circle, draw=black, fill=white, inner sep=0pt, minimum width=5pt];
    \path (0,1) node (p1)  {} ;
    \path (-1,0) node (p2)  {} ;
    \path (1,0) node (p3)  {} ;
   \path (-1,-1) node (p4) {} ;
   \path (1,-1) node (p5)  {} ;
      \draw (p1) -- (p4);
     \draw (p1) -- (p5);
     \draw (p1) -- (p2);
     \draw (p1) -- (p3);
     \draw (p2) -- (p4);
     \draw (p3) -- (p5);
     \draw (p2) -- (p5);
     \draw (p3) -- (p4);
       \draw [dashed, thin] (0,-1.6) -- (0,1.6);
        \node [rectangle, draw=white, fill=white] (a) at (-0.7,0.6) {$\alpha$};
     \node [rectangle, draw=white, fill=white] (a) at (0.9,0.6) {$-\alpha$};
     \node [rectangle, draw=white, fill=white] (a) at (-0.25,0.05) {$\beta$};
     \node [rectangle, draw=white, fill=white] (a) at (0.2,0.05) {$-\beta$};
     \node [rectangle, draw=white, fill=white] (a) at (0.5,-1.1) {$-\gamma$};
     \node [rectangle,draw=white, fill=white] (a) at (-0.5,-1.1) {$\gamma$};
       \node [rectangle, draw=white, fill=white] (a) at (-1.3,-0.5) {$\delta$};
     \node [rectangle, draw=white, fill=white] (a) at (1.3,-0.5) {$-\delta$};
       \node [draw=white, fill=white] (b) at (0,-2.1) {(b)};
        \end{tikzpicture}
        \end{center}
\vspace{-0.3cm}
\caption{Self-stressed frameworks in the plane: (a) a fully $(\mathcal{C}_s,\Phi)$-symmetric self-stress of $(K_{4},p) \in \mathscr{R}_{(K_{4},\mathcal{C}_s, \Phi)}$; (b) a  self-stress of $(G,p) \in \mathscr{R}_{(G,\mathcal{C}_s, \Psi)}$ which is not fully $(\mathcal{C}_s,\Psi)$-symmetric. The types $\Phi$ and $\Psi$ are uniquely determined by the injective realizations \cite{BS4,BS1}.}
\label{fulsymstr}
\end{figure}
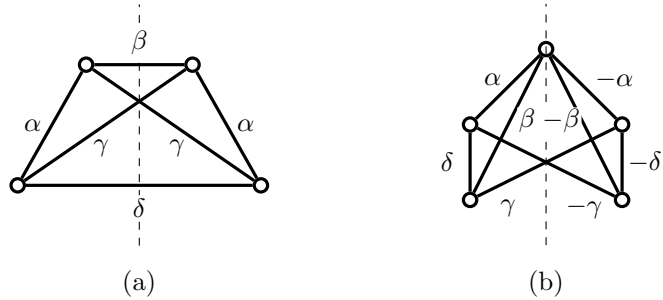

A self-stress $\omega$ of a  framework $(G,p)\in \mathscr{R}_{(G,S,\Phi)}$  is \emph{fully $(S,\Phi)$-symmetric} if $(\omega)_e=(\omega)_f$ whenever $e$ and $f$ belong to the same orbit $S(e)=\{\Phi(x)(e)|\,x\in S\}$ of edges of $G$ (see also Figure \ref{fulsymstr}(a)). \\\indent
Note that a fully $(S,\Phi)$-symmetric self-stress is clearly uniquely determined by the components $(\omega)_1,\ldots,(\omega)_r$, whenever $\mathscr{O}_{E(G)}=\{e_1,\ldots, e_r\}$ is a set of representatives for the edge orbits  $S(e)=\{\Phi(x)(e)|\,x\in S\}$ of $G$.

It is shown in \cite{KG2, BS2} that the rigidity matrix of a framework $(G,p)\in \mathscr{R}_{(G,S,\Phi)}$ can be transformed into a block-diagonalized form using techniques from group representation theory. In this block-diagonalization of $\mathbf{R}(G,p)$, the submatrix block $\tilde{\mathbf{R}}_1(G,p)$ that corresponds to the trivial irreducible representation of $S$ describes the relationship between external displacement vectors on the joints and resulting internal distortion vectors in the bars of $(G,p)$ that are fully $(S,\Phi)$-symmetric. So, the submatrix block $\tilde{\mathbf{R}}_1(G,p)$ comprises all the information regarding the fully $(S,\Phi)$-symmetric infinitesimal rigidity properties of $(G,p)$. The orbit rigidity matrix  of  $(G,p)$ which we will introduce in the next section will have the same properties as the submatrix block $\tilde{\mathbf{R}}_1(G,p)$; however, we will see that the orbit rigidity matrix  allows a significantly simplified fully $(S,\Phi)$-symmetric infinitesimal rigidity analysis of $(G,p)$, since it can be set up directly without finding the block-diagonalization of $\mathbf{R}(G,p)$
or using other group representation techniques.

\section{The orbit rigidity matrix}
\label{sec:orbitmatdef}

To make the general definition of the orbit rigidity matrix more transparent, we first consider a few simple examples.

\begin{examp}\label{c2quadexamp}
Consider the $2$-dimensional framework $(K_{2,2},p)\in\mathscr{R}_{(K_{2,2},\mathcal{C}_2,\Phi)}$ depicted in Figure \ref{quadc2pic}, where $\Phi:\mathcal{C}_2\to \textrm{Aut}(K_{2,2})$ is the homomorphism defined by $\Phi(C_2)=(1\, 3)(2 \, 4)$.\\\indent
If we denote $p_1^T=(a,b)$, $p_2^T=(c,d)$, $p_3^T=(-a,-b)$, and $p_4^T=(-c,-d)$, then the rigidity matrix of $(K_{2,2},p)$ is
\begin{displaymath}\bordermatrix{
                &1&2& {3}=C_2(1) &{4}=C_2(2) \cr
                \{1,2\}&(a-c,b-d) &  (c-a,d-b)  &0\ 0 & 0 \ 0\cr
                 \{1,4\}& (a+c,b+d) & 0 \ 0  & 0 \ 0 & (-a-c,-b-d)\cr
             C_2\{1,2\}& 0 \ 0  & 0 \ 0 & (c-a,d-b) &  (a-c,b-d) \cr
             C_2\{1,4\}& 0 \ 0  &  (a+c,b+d)       & (-a-c,-b-d) & 0 \ 0}
\end{displaymath}
\begin{figure}[htp]
\begin{center}
\begin{tikzpicture}[very thick,scale=1]
\tikzstyle{every node}=[circle, draw=black, fill=white, inner sep=0pt, minimum width=5pt];
        \path (0,0) node (p1) [label = left: $p_{1}$] {} ;
       \path (0,-1.1) node (p2) [label = left: $p_{2}$] {} ;
   \path (2.5,-1.8) node (p3) [label = right: $p_{3}$] {} ;
   \path (2.5,-0.7) node (p4) [label = right: $p_{4}$] {} ;
   \draw (p1) -- (p4);
      \draw (p3) -- (p4);
     \draw (p2) -- (p3);
      \draw (p2) -- (p1);
\filldraw[fill=black, draw=black]
    (1.25,-0.9) circle (0.01cm);
\node [rectangle, draw=white, fill=white] (a) at (01.25,-0.7) {\small center};
         \end{tikzpicture}
\end{center}
\vspace{-0.3cm}
\caption{A framework $(K_{2,2},p)\in\mathscr{R}_{(K_{2,2},\mathcal{C}_2,\Phi)}$.}
\label{quadc2pic}
\end{figure}
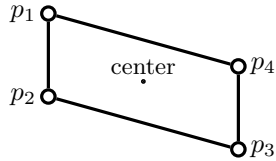
This matrix has rank 4, and hence leaves a space of $8-4 = 4$ infinitesimal motions. Thus, there exists a $1$-dimensional space of  infinitesimal flexes of $(K_{2,2},p)$ spanned by $u=\left(\begin{array} {cccccccc}-1 & 0 & x & y & 1 & 0 & -x &  -y \end{array}\right)^T$, where $ x=\frac{cd-ab}{ad-bc}$ and $y=-\frac{c^2-a^2}{ad-bc}$. This infinitesimal flex is clearly fully $(\mathcal{C}_2,\Phi)$-symmetric.\\\indent
Note that if we are only interested in infinitesimal motions and self-stresses of $(K_{2,2},p)$ that are fully $(\mathcal{C}_2,\Phi)$-symmetric, then it suffices to focus on the first two rows of $\mathbf{R}(K_{2,2},p)$ (i.e., the rows corresponding to the representatives $\{1,2\}$ and $\{1,4\}$ for the edge orbits $S(e)=\{\Phi(x)(e)|\,x\in \mathcal{C}_2\}$ of $K_{2,2}$). The other two rows are redundant in this fully symmetric context. So, the orbit rigidity matrix for $(K_{2,2},p)$ will have two rows, one for each representative of the edge orbits under the action of $\mathcal{C}_2$. Further, the orbit rigidity matrix will have only four columns, because each of the joints $p_1$ and $p_2$ has two degrees of freedom, and the displacement vectors at the joints $p_3=C_2(p_1)$ and $p_4=C_2(p_2)$ are uniquely determined by the displacement vectors at the joints $p_1$ and $p_2$ and the symmetry constraints given by $\mathcal{C}_2$ and $\Phi$. We write the orbit rigidity matrix of $(K_{2,2},p)$ as follows:
\begin{displaymath}\bordermatrix{
                &1&2 \cr
                \{1,2\}&(p_1-p_2)^T &  (p_2-p_1)^T\cr
                 \{1,4\}& \big(p_1-C_2(p_2)\big)^T & \big(p_2-C_2^{-1}(p_1)\big)^T \cr
             }=\bordermatrix{
                &1&2 \cr
                &(a-c,b-d) &  (c-a,d-b)\cr
                 & (a+c,b+d) & (c+a,d+b) \cr
             }
\end{displaymath}
\end{examp}

\begin{examp}\label{1csquadexamp1}
The orbit rigidity matrix for the framework $(K_{2,2},p)\in \mathscr{R}_{(K_{2,2},\mathcal{C}_s,\Phi_a)}$ in Example \ref{quadcs1}
(Figure~2(a))
has again two rows, since $K_{2,2}$ has two edge orbits (each of size 2) under the  action of $\mathcal{C}_s$. The vertex orbits are represented by the vertices $1$, $2$ and $4$, for example. Clearly, the joint $p_1$ has two degrees of freedom, which gives rise to two columns in the orbit matrix. The joints $p_2$ and $p_4$, however, are fixed by the reflection $s$ in $\mathcal{C}_s$, so that any fully $(\mathcal{C}_s,\Phi_a)$-symmetric displacement vectors at $p_2$ and $p_4$ must lie on the mirror corresponding to $s$ (i.e., on the $y$-axis). Thus, the orbit rigidity matrix of $(K_{2,2},p)$ has only one column for each of the joints $p_2$ and $p_4$:
\begin{eqnarray}& & \bordermatrix{
                &1&2 & 4\cr
                \{1,2\}&(p_1-p_2)^T &  (p_2-p_1)^T\binom{0}{1} & 0\cr
                 \{1,4\}& (p_1-p_4)^T & 0 &  (p_4-p_1)^T \binom{0}{1}\cr
             }\nonumber\\
             & = & \bordermatrix{
                &1&2 & 4\cr
                &(a,b-c) &  (c-b) & 0\cr
                 & (a,b-d) & 0 & (d-b) \cr
             }\nonumber
\end{eqnarray}
\end{examp}

\begin{examp}\label{2csquadexamp2}
The orbit rigidity matrix for the framework $(K_{2,2},p)\in \mathscr{R}_{(K_{2,2},\mathcal{C}_s,\Phi_b)}$ in Example \ref{quadcs2}
(Figure~2(b))
 is a $3\times4$ matrix, since there are three edge orbits - represented by the edges $\{1,2\}$, $\{1,4\}$, and  $\{2,3\}$, for example - and two vertex orbits - represented by the vertices $1$ and $2$, for example, and none of the joints of $(K_{2,2},p)$ are fixed by the reflection in $\mathcal{C}_s$. Note, however, that the end-vertices of the edge $\{1,4\}$
lie in the same vertex orbit for $1$ under the action of $\mathcal{C}_s$
 and that the end-vertices of the edge $\{2,3\}$ lie in the same vertex orbit for $2$ under the action of $\mathcal{C}_s$. Thus, for the orbit rigidity matrix of $(K_{2,2},p)$ we write
\begin{displaymath}\bordermatrix{
                &1&2 \cr
                \{1,2\}&(p_1-p_2)^T &  (p_2-p_1)^T\cr
                 \{1,4\}& 2\big(p_1-s(p_1)\big)^T & 0 \ 0 \cr
                 \{2,3\}& 0 \ 0 & 2\big(p_2-s(p_2)\big)^T \cr
             }=\bordermatrix{
                &1&2 \cr
                &(a-c,b-d) &  (c-a,d-b)\cr
                 & (4a,0) & 0 \ 0 \cr
                 & 0 \ 0 & (4c,0) \cr
             }
\end{displaymath}
\end{examp}

We now give the general definition of the orbit rigidity matrix of a symmetric framework.

\begin{defin} Let $G$ be a graph, $S$ be a symmetry group in dimension $d$, $\Phi:S\to \textrm{Aut}(G)$ be a homomorphism, and $(G,p)$ be a framework in $\mathscr{R}_{(G,S,\Phi)}$. Further, let $\mathscr{O}_{V(G)}=\{1,\ldots, k\}$ be a set of representatives for the orbits  $S(i)=\{\Phi(x)(i)|\,x\in S\}$ of vertices of $G$. We construct the \emph{orbit rigidity matrix} (or in short, \emph{orbit matrix}) $\mathbf{O}(G,p,S)$ of $(G,p)$ so that it has exactly one row for each orbit $S(e)=\{\Phi(x)(e)|\,x\in S\}$ of edges of $G$ and exactly $c_i=\textrm{dim }\big(U(p_i)\big)$ columns for each vertex $i\in \mathscr{O}_{V(G)}$.
\\\indent Given an edge orbit $S(e)$ of $G$, there are two possibilities for the corresponding row in $\mathbf{O}(G,p,S)$:

\begin{description}
\item[Case 1:] The two end-vertices of the edge $e$ lie in distinct vertex orbits. Then there exists an edge in $S(e)$ that is of the form $\{a,x(b)\}$ for some $x\in S$, where $a,b\in\mathscr{O}_{V(G)}$. Let a basis $\mathscr{B}_{a}$ for $U(p_a)$ and a basis  $\mathscr{B}_{b}$ for $U(p_b)$ be given and let $\mathbf{M}_{a}$ and  $\mathbf{M}_{b}$ be the matrices whose columns are the coordinate vectors of  $\mathscr{B}_{a}$ and $\mathscr{B}_{b}$ relative to the canonical basis of $\mathbb{R}^d$, respectively.
The row we write in $\mathbf{O}(G,p,S)$ is:
        \begin{displaymath}\renewcommand{\arraystretch}{0.8}
    \bordermatrix{  & &  a &  &  b &  \cr & 0 \ldots 0 & \big(p_a-x(p_b)\big)^T\mathbf{M}_{a} & 0  \ldots  0 & \big(p_b-x^{-1}(p_a)\big)^T\mathbf{M}_{b} & 0  \ldots  0}\textrm{.}
    \end{displaymath}
\item[Case 2:] The two end-vertices of the edge $e$ lie in the same vertex orbit. Then there exists an edge in $S(e)$ that is of the form $\{a,x(a)\}$ for some $x\in S$, where $a\in\mathscr{O}_{V(G)}$. Let a basis $\mathscr{B}_{a}$ for $U(p_a)$ be given and let $\mathbf{M}_{a}$  be the matrix whose columns are the coordinate vectors of  $\mathscr{B}_{a}$ relative to the canonical basis of $\mathbb{R}^d$. The row we write in $\mathbf{O}(G,p,S)$ is:
\begin{displaymath}
\bordermatrix{ & &  a & \cr & 0  \ldots  0 & \big(2p_a-x(p_a)-x^{-1}(p_a)\big)^T\mathbf{M}_{a} & 0  \ldots  0}\textrm{.}
    \end{displaymath}
In particular, if $x(p_a)= x^{-1}(p_a)$, this row becomes
\begin{displaymath}
\bordermatrix{ & &  a & \cr & 0  \ldots  0 & 2\big(p_a-x(p_a)\big)^T\mathbf{M}_{a} & 0  \ldots  0}\textrm{.}
   \end{displaymath}
\end{description}
\end{defin}

\begin{remark}\label{rankrem}
\emph{Note that the rank of the orbit rigidity matrix $\mathbf{O}(G,p,S)$ is clearly independent of the choice of bases  for the spaces $U(p_a)$ (and their corresponding matrices $\mathbf{M}_a$), $a=1,\ldots, k$.}
\end{remark}

\begin{remark}\label{nothingfixedrem}
\emph{If none of the joints of $(G,p)$ are fixed by any non-trivial symmetry operation in $S$, then the orbit rigidity matrix $\mathbf{O}(G,p,S)$ of $(G,p)$ has  $dk=d|\mathscr{O}_{V(G)}|$ columns, and each of the matrices $\mathbf{M}_{a}$ and $\mathbf{M}_{b}$ may be chosen to be the $d\times d$ identity matrix. In this case, the matrix $\mathbf{O}(G,p,S)$ becomes particularly easy to construct (see Examples  \ref{c2quadexamp} and \ref{2csquadexamp2}).}
\end{remark}


\section{The kernel of $\mathbf{O}(G,p,S)$}
\label{sec:kernel}

In this section, we show that the kernel of the orbit rigidity matrix $\mathbf{O}(G,p,S)$ of a symmetric framework $(G,p)\in \mathscr{R}_{(G,S,\Phi)}$ is the space of all fully $(S,\Phi)$-symmetric infinitesimal motions of $(G,p)$, restricted to the set  $\mathscr{O}_{V(G)}$ of representatives for the vertex orbits  $S(i)$ of $G$ (Theorem \ref{thm:kernel}). It follows from this result that we can detect whether $(G,p)$ has a fully $(S,\Phi)$-symmetric infinitesimal flex by simply computing the rank of $\mathbf{O}(G,p,S)$.

\begin{theorem}
\label{thm:kernel}
Let $G$ be a graph, $S$ be a symmetry group in dimension $d$, $\Phi:S\to \textrm{Aut}(G)$ be a homomorphism, $\mathscr{O}_{V(G)}=\{1,\ldots, k\}$ be a set of representatives for the orbits  $S(i)=\{\Phi(x)(i)|\,x\in S\}$ of vertices of $G$, and $(G,p)\in \mathscr{R}_{(G,S,\Phi)}$. Further, for each $i=1,\ldots k$, let a basis $\mathscr{B}_{i}$ for $U(p_i)$ be given and let $\mathbf{M}_{i}$  be the  $d\times c_i$  matrix whose columns are the coordinate vectors of  $\mathscr{B}_{i}$ relative to the canonical basis of $\mathbb{R}^d$. Then \begin{displaymath}\tilde{u}=\left(
\begin{array} {c}
\tilde{u}_1\\ \vdots \\ \tilde{u}_k
\end{array}
\right)\in \mathbb{R}^{c_1}\times \ldots \times \mathbb{R}^{c_k}\end{displaymath} lies in the kernel of $\mathbf{O}(G,p,S)$ if and only if \begin{displaymath}\overline{u}=\left(
\begin{array} {c}
\mathbf{M}_{1}\tilde{u}_1\\ \vdots \\ \mathbf{M}_{k}\tilde{u}_k
\end{array}
\right)\in \mathbb{R}^{dk}\end{displaymath} is the restriction $u|_{\mathscr{O}_{V(G)}}$ of a fully $(S,\Phi)$-symmetric infinitesimal motion $u$ of $(G,p)$ to $\mathscr{O}_{V(G)}$.
\end{theorem}
\textbf{Proof.} 
 Suppose there exists an edge $e=\{a,x(b)\}$ in $G$ whose two end-vertices  lie in distinct vertex orbits (see Case 1 in the definition of the orbit rigidity matrix). The row equation of the matrix $\mathbf{O}(G,p,S)$ for the edge orbit represented by $e$ is then of the form
 \begin{displaymath}
   \big(p_a-\mathbf{X}p_b\big)^T\big(\mathbf{M}_{a}\tilde{u}_a\big) + \big(p_b-\mathbf{X}^{-1}p_a\big)^T\big(\mathbf{M}_{b}\tilde{u}_b\big) = 0\textrm{,}
    \end{displaymath}
where $\mathbf{X}$ is the matrix that represents $x$ with respect to the canonical basis of $\mathbb{R}^d$. Since the inner product in the second summand is invariant under the orthogonal transformation $x\in S$, we have
 \begin{displaymath}
   \big(p_a-\mathbf{X}p_b\big)^T\big(\mathbf{M}_{a}\tilde{u}_a\big) + \big(\mathbf{X}p_b-p_a\big)^T \big(\mathbf{X}\mathbf{M}_{b}\tilde{u}_b\big) = 0\textrm{,}
    \end{displaymath}
which is the row equation of the standard rigidity matrix $\mathbf{R}(G,p)$ for $e=\{a,x(b)\}$.
\\\indent Similarly, for any other edge $y(\{a,x(b)\})$, $y\in S$, that lies in the edge orbit $S(e)$, we have
\begin{displaymath}
   \Big(\mathbf{Y}p_a-\mathbf{Y}\mathbf{X}p_b\Big)^T\big(\mathbf{Y}\mathbf{M}_{a}\tilde{u}_a\big) + \Big(\mathbf{Y}\mathbf{X}p_b-\mathbf{Y}p_a\Big)^T \big(\mathbf{Y}\mathbf{X} \mathbf{M}_{b}\tilde{u}_b\big) = 0\textrm{,}
    \end{displaymath}
where $\mathbf{Y}$ is the matrix that represents $y$ with respect to the canonical basis of $\mathbb{R}^d$. This is the standard row equation of $\mathbf{R}(G,p)$ for the edge $y(\{a,x(b)\})$.\\\indent
Suppose next that there exists a bar $e=\{a,x(a)\}$ in $G$ whose two end-vertices  lie in the same vertex orbit 
 (see Case 2 in the definition of the orbit rigidity matrix). The row equation of the matrix $\mathbf{O}(G,p,S)$ for the edge orbit represented by $e$ is then of the form
 \begin{displaymath}
  \big(p_a-\mathbf{X}p_a\big)^T\big(\mathbf{M}_{a}\tilde{u}_a\big)  + \big(p_a-\mathbf{X}^{-1}p_a\big)^T\big(\mathbf{M}_{a}\tilde{u}_a\big) = 0\textrm{.}
    \end{displaymath}
Since the inner product in the second summand is invariant under the orthogonal transformation $x\in S$, we have
 \begin{displaymath}
  \big(p_a-\mathbf{X}p_a\big)^T\big(\mathbf{M}_{a}\tilde{u}_a\big)  + \big(\mathbf{X}p_a-p_a\big)^T\big(\mathbf{X}\mathbf{M}_{a}\tilde{u}_a\big) = 0\textrm{,}
    \end{displaymath}
which is the standard row equation of $\mathbf{R}(G,p)$ for $e=\{a,x(a)\}$.
\\\indent Similarly, for any other edge $y(\{a,x(a)\})$, $y\in S$, that lies in the edge orbit $S(e)$, we have
 \begin{displaymath}
  \big(\mathbf{Y}p_a-\mathbf{Y}\mathbf{X}p_a\big)^T\big(\mathbf{Y}\mathbf{M}_{a}\tilde{u}_a\big)  + \big(\mathbf{Y}\mathbf{X}p_a-\mathbf{Y}p_a\big)^T\big(\mathbf{Y}\mathbf{X}\mathbf{M}_{a}\tilde{u}_a\big) = 0\textrm{,}
    \end{displaymath}
which is the standard row equation of $\mathbf{R}(G,p)$ for the edge $y(\{a,x(a)\})$.
\\\indent It follows that $\tilde{u}$ lies in the kernel of $\mathbf{O}(G,p,S)$ if and only if  $\overline{u}$ is the restriction $u|_{\mathscr{O}_{V(G)}}$ of a fully $(S,\Phi)$-symmetric infinitesimal motion $u$ of $(G,p)$ to $\mathscr{O}_{V(G)}$. $\square$

\begin{examp} Consider the framework $(K_{2,2},p)\in \mathscr{R}_{(K_{2,2},\mathcal{C}_s,\Phi_a)}$ from Examples \ref{quadcs1} and \ref{1csquadexamp1}. The vector $\tilde{u}=\left(\begin{array} {cccccccc}-1 & 0 & \frac{a}{c-b} & \frac{a}{d-b}\end{array}\right)^T$ clearly lies in the kernel of $\mathbf{O}(K_{2,2},p,\mathcal{C}_s)$, and the vector
\begin{displaymath}\overline{u}=
\left(\begin{array} {c}\mathbf{M}_1\binom{-1}{0}\\ \mathbf{M}_2\frac{a}{c-b}\\ \mathbf{M}_4\frac{a}{d-b}= \end{array}\right)
= \left(\begin{array} {c}\mathbf{Id}\binom{-1}{0}\\ \binom{0}{1}\frac{a}{c-b}\\ \binom{0}{1}\frac{a}{d-b}= \end{array}\right)
= \left(\begin{array} {c}-1 \\0\\ 0 \\\frac{a}{c-b}\\ 0 \\ \frac{a}{d-b} \end{array}\right)
\end{displaymath}
is the restriction of the fully $(\mathcal{C}_s,\Phi_a)$-symmetric  infinitesimal flex $u=\left(\begin{array} {cccccccc}-1 & 0 & 0 & \frac{a}{c-b} & 1 & 0 &  0 & \frac{a}{d-b}\end{array}\right)^T$ to the set $\{1,2,4\}$ of representatives for the vertex orbits $\mathcal{C}_s(i)=\{\Phi_a(x)(i)|\,x\in \mathcal{C}_s\}$.
\end{examp}

\section{Symmetry-preserving finite flexes}
\label{sec:finite} %
The following extension of the theorem of Asimov and Roth (see \cite{asiroth}) to frameworks that possess non-trivial symmetries was
derived
in \cite{BS6} (see also \cite{BS4}):

\begin{theorem}
\label{thm:flexes}
Let $G$ be a graph, $S$ be a symmetry group in dimension $d$, $\Phi:S\to \textrm{Aut}(G)$ be a homomorphism, and $(G,p)$ be a framework in $\mathscr{R}_{(G,S,\Phi)}$ whose joints span all of $\mathbb{R}^{d}$. If $(G,p)$ is $(S,\Phi)$-generic and $(G,p)$ has a fully $(S,\Phi)$-symmetric infinitesimal flex, then there also exists a finite flex of $(G,p)$ which preserves the symmetry of $(G,p)$ throughout the path.
\end{theorem}

\begin{remark} \emph{It is also shown in \cite{BS6} that  the condition that $(G,p)$ is $(S,\Phi)$-generic in Theorem \ref{thm:flexes} may be replaced  by the weaker condition that the submatrix block $\tilde{\mathbf{R}}_1(G,p)$ of the block-diagonalized rigidity matrix $\tilde{\mathbf{R}}(G,p)$ which corresponds to the trivial irreducible representation of $S$ (or, equivalently, the orbit rigidity matrix  $\mathbf{O}(G,p,S)$) has maximal rank in some neighborhood of the configuration $p$ within the space of configurations that satisfy the symmetry constraints given by $S$ and $\Phi$. In particular, this says that if the rows of the orbit rigidity matrix $\mathbf{O}(G,p,S)$ are linearly independent and $(G,p)$ has a fully $(S,\Phi)$-symmetric infinitesimal flex, then $(G,p)$ also has a symmetry-preserving finite flex.}
\end{remark}

In combination with Theorem \ref{thm:flexes}, Theorem \ref{thm:kernel} gives rise to a simple new method for detecting finite flexes in symmetric frameworks. In the following
subsections, we elaborate on this new method and apply it to a number of interesting examples.


\subsection{Detection of finite flexes from the size of the orbit rigidity matrix}

First, we consider situations where knowledge of only the \emph{size} of the orbit rigidity matrix already allows us to detect finite flexes in symmetric frameworks.

The following result is an immediate consequence of Theorem \ref{thm:kernel}

\begin{theorem}\label{symMaxthm}
Let $G$ be a graph, $S$ be a symmetry group in dimension $d$, $\Phi:S\to \textrm{Aut}(G)$ be a homomorphism, and $(G,p)$ be a framework in $\mathscr{R}_{(G,S,\Phi)}$. Further, let $r$ and $c$ denote the number of rows and columns of the orbit rigidity matrix $\mathbf{O}(G,p,S)$, respectively, and  let $m$ denote the dimension of the space of fully $(S,\Phi)$-symmetric infinitesimal rigid motions of $(G,p)$. If \begin{equation}\label{eqflex} r<c - m \textrm{,}\end{equation} then $(G,p)$ has a fully $(S,\Phi)$-symmetric infinitesimal flex.
\end{theorem}

Recall from Section \ref{sec:orbitmatdef} that  the number of rows, $r$, and the number of columns, $c$, of the orbit rigidity matrix $\mathbf{O}(G,p,S)$ of $(G,p)$ do not depend on the configuration $p$, but only on the graph $G$ and the prescribed symmetry constraints given by $S$ and $\Phi$. As shown in \cite{BS4}, the dimension $m$ of the space of fully $(S,\Phi)$-symmetric infinitesimal rigid motions of $(G,p)$ is also independent of $p$, provided that the joints of $(G,p)$ span all of $\mathbb{R}^{d}$. So, suppose the set $\mathscr{R}_{(G,S,\Phi)}$ contains a framework whose joints span all of $\mathbb{R}^{d}$. Then, as shown in \cite{BS1}, the joints of all $(S,\Phi)$-generic realizations of $G$ also span all of $\mathbb{R}^{d}$. Thus, if (\ref{eqflex}) holds, then all $(S,\Phi)$-generic realizations of $G$ have a fully $(S,\Phi)$-symmetric infinitesimal flex, and hence, by Theorem \ref{thm:flexes}, also a finite symmetry-preserving flex.

The dimension $m$ of the space of fully $(S,\Phi)$-symmetric infinitesimal rigid motions of $(G,p)$ can easily be computed using the techniques described in \cite{BS2, BS4}. In particular, in dimension 2 and 3,  $m$ can be deduced immediately from the character tables given in \cite{cfgsw}. Thus, in order to check condition (\ref{eqflex}) it is only left to determine the size of the orbit rigidity matrix $\mathbf{O}(G,p,S)$ which basically requires only a simple count of the vertex orbits and edge orbits of the graph $G$.


Alternatively, the values of $r$ and $c$ in (\ref{eqflex}) can also be found by expressing the characters of the `internal' and `external' matrix representation  for the group $S$ (see \cite{FGsymmax,KG2, BS2, BS1}, for example)  as linear combinations of the characters of the irreducible representations of $S$: the numbers $r$ and $c$ are the respective coefficients corresponding to the trivial irreducible representation in these linear combinations (see \cite{BS4, BS6} for details). However, our new `orbit approach' is much simpler than this method of computing characters, since it allows us to determine $r$ and $c$ directly without using any techniques from group representation theory.

\subsubsection{Examples in dimension 2}

Let's apply our new method to the symmetric quadrilaterals we discussed in Section \ref{sec:orbitmatdef}.

\begin{figure}
    \begin{center}
  \subfigure[] {\includegraphics [width=.30\textwidth]{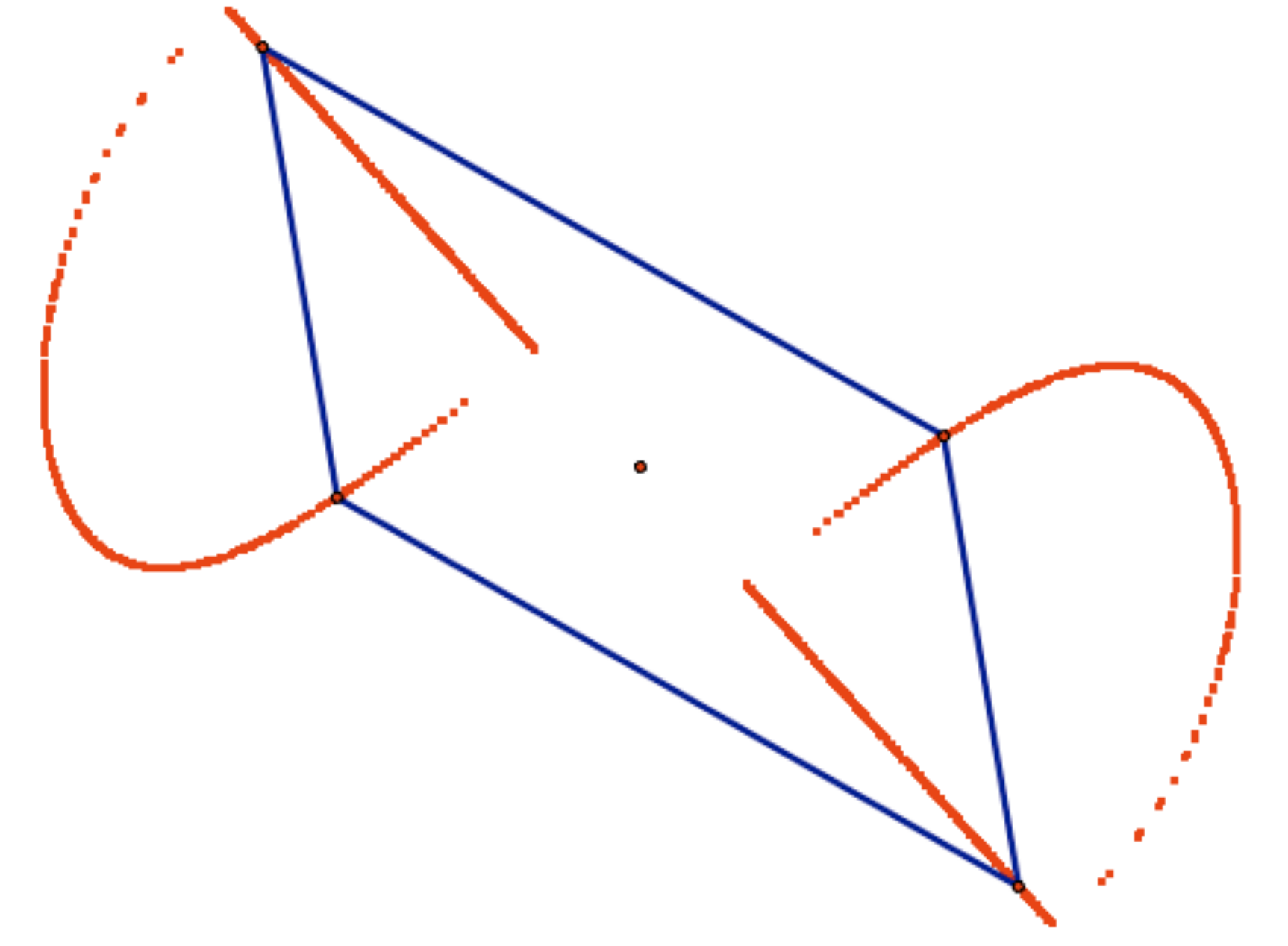}}\quad\quad\quad
  \subfigure[] {\includegraphics [width=.20\textwidth]{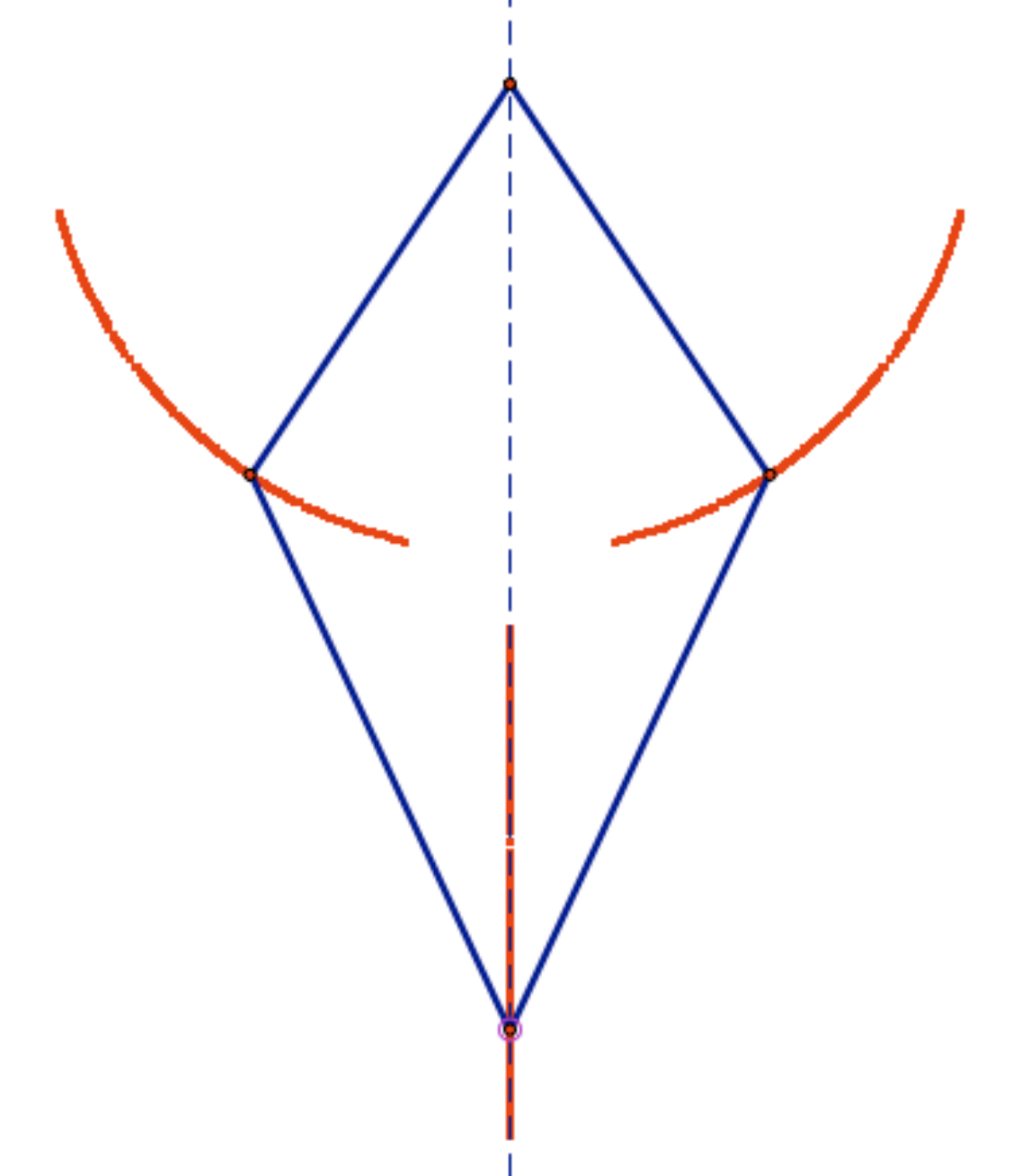}}
       \end{center}
    \caption{Symmetry-preserving mechanisms of the quadrilaterals from Examples \ref{c2quadexamp} (a) and \ref{quadcs1} (b).}
    \label{fig:mechanisms}
    \end{figure}

We first consider the quadrilateral with point group $\mathcal{C}_2$ from Example \ref{c2quadexamp} (see Figure \ref{quadc2pic}). There are two vertex orbits, represented by the vertices $1$ and $2$, for example, and we have $c_i=\textrm{dim}\big(U(p_i)\big)=2$ for $i=1,2$. Further, there are two edge orbits, and we have $m=1$, since the only infinitesimal rigid motions that are fully symmetric with respect to this half-turn symmetry are the ones that correspond to rotations about the origin (see \cite{BS4} for details). Thus, we have
\begin{displaymath}
r=2<3=c - m\textrm{.}
\end{displaymath}
So, we may conclude that $(\mathcal{C}_2,\Phi)$-generic realizations of $K_{2,2}$ have a symmetry-preserving mechanism (see also Figure \ref{fig:mechanisms}(a)).

This result can easily be generalized to predict the flexibility of a whole class of $2$-dimensional frameworks with half-turn symmetry.\\\indent
Recall from Section \ref{sec:sym} that a joint $p_i$ of $(G,p)\in \mathscr{R}_{(G,S,\Phi)}$ is said to be \emph{fixed} by $x\in S$ if $x(p_i)=p_i$. Similarly, we say that a bar $\{p_i,p_j\}$ of $(G,p)$ is \emph{fixed} by $x$ if either $x(p_i)=p_i$ and $x(p_j)=p_j$ or $x(p_i)=p_j$ and $x(p_j)=p_i$. The number of joints and bars of $(G,p)$ that are fixed by $x$ are denoted by $j_x$ and $b_x$, respectively.

\begin{theorem}\label{2dhalfturn}
Let $G$ be a graph with $|E(G)|=2|V(G)|-4$, $\mathcal{C}_2=\{Id, C_2\}$ be the half-turn symmetry group in dimension $2$, and $\Phi:\mathcal{C}_2\to \textrm{Aut}(G)$ be a homomorphism. If $j_{C_2}=b_{C_2}=0$, then $(\mathcal{C}_2,\Phi)$-generic realizations of $G$ have a symmetry-preserving mechanism.
\end{theorem}
\textbf{Proof.} Since $j_{C_2}=0$ we have $c=2\frac{|V(G)|}{2}=|V(G)|$, and since $b_{C_2}=0$ we have $r=\frac{|E(G)|}{2}=|V(G)|-2$. As mentioned above, we have  $m=1$, so that
\begin{displaymath} r=|V(G)|-2< |V(G)|-1=c-m \textrm{.}\nonumber
\end{displaymath}
Thus, by Theorem \ref{thm:kernel} and \ref{thm:flexes}, $(\mathcal{C}_2,\Phi)$-generic realizations of $G$ have a symmetry-preserving mechanism. $\square$

Next, we consider the quadrilateral with point group $\mathcal{C}_s$ from Example \ref{quadcs1} (see Figure \ref{fulsym}(a)). There are three vertex orbits, represented by the vertices $1$, $2$, and $4$, for example, and we have $c_i=\textrm{dim}\big(U(p_i)\big)=1$ for $i=2,4$ and $c_1=\textrm{dim}\big(U(p_1)\big)=2$. Further, there are two edge orbits, and we have $m=1$, since the only infinitesimal rigid motions that are fully symmetric with respect to this mirror symmetry are the ones that correspond to translations along the mirror line \cite{BS4}. Thus, we have
\begin{displaymath}
r=2<3=c - m\textrm{.}
\end{displaymath}
So, we may conclude that $(\mathcal{C}_s,\Phi_a)$-generic realizations of $K_{2,2}$ have a symmetry-preserving mechanism  (see also Figure \ref{fig:mechanisms}(b)).

The following theorem provides a generalized version of this result.

\begin{theorem}\label{2dmirror}
Let $G$ be a graph with $|E(G)|=2|V(G)|-4$, $\mathcal{C}_s=\{Id, s\}$ be a `reflectional' symmetry group in dimension $2$, and $\Phi:\mathcal{C}_s\to \textrm{Aut}(G)$ be a homomorphism. If  $b_{s}=0$, then $(\mathcal{C}_s,\Phi)$-generic realizations of $G$ have a symmetry-preserving mechanism.
\end{theorem}
\textbf{Proof.} Let $|V'|$ be the number of vertex orbits of size 2 and $|E'|$ be the number of edge orbits of size 2. Then we have
\begin{eqnarray}
|E(G)| & = & 2|V(G)|-4\nonumber\\
2|E'|+b_{s} & = & 2(2|V'|+j_{s})-4\nonumber\\
|E'|+\frac{1}{2}b_{s} & = & 2|V'|+j_{s}-2\nonumber\\
|E'|+b_{s} & = & 2|V'|+j_{s}-1+\frac{1}{2}b_s-1\nonumber\\
r & = & c-m+\frac{1}{2}b_s-1 \textrm{.} \nonumber
\end{eqnarray}
So if $b_{s}<2$ (or, equivalently, $b_{s}=0$, since $b_{s}=1$ contradicts the count $|E(G)|=2|V(G)|-4$)  then $r<c-m$.
The result now follows from Theorems \ref{thm:kernel} and \ref{thm:flexes}. $\square$

Similar results to Theorems \ref{2dhalfturn} and \ref{2dmirror} can of course also be established for  other symmetry groups in dimension $2$.

Note that the orbit count for the quadrilateral with mirror symmetry from Example \ref{quadcs2} (see Figure \ref{fulsym}(b,c)) is
\begin{displaymath}
r=3=c - m\textrm{,}
\end{displaymath}
and by computing the rank of the corresponding orbit rigidity matrix explicitly, it is easy to verify that this quadrilateral does in fact not have any fully $(\mathcal{C}_s,\Phi_b)$-symmetric infinitesimal flex, let alone a symmetry-preserving mechanism.  However, it does have a mechanism that breaks the mirror symmetry.

\subsubsection{Examples in dimension 3}

The Bricard octahedra \cite{bricard} are famous examples of flexible frameworks in 3-space.
While it follows from Cauchy's Theorem (\cite{Cauchy}) that convex realizations of the octahedral graph are isostatic, the French engineer R. Bricard found three types of octahedra with self-intersecting faces whose realizations as frameworks are flexible. Two of these three types of Bricard octahedra possess non-trivial symmetries: Bricard octahedra of the first type have a half-turn symmetry and Bricard octahedra of the second type have a mirror symmetry.
In the following, we use our new `orbit approach' to not only show that both of these types of frameworks are flexible, but also that they possess a `symmetry-preserving' finite flex. Various other treatments of the Bricard octahedra can be found in \cite{Baker, BS6, Stachel}, for example. R. Connelly's celebrated counterexample to Euler's rigidity conjecture from 1776  is also based on a flexible Bricard octahedron (of the first type) \cite{concounter}.
However, since the flexible Connelly sphere - as well as Steffen's modified Connelly sphere - break the half-turn symmetry as they flex, our methods do not apply to these particular examples.

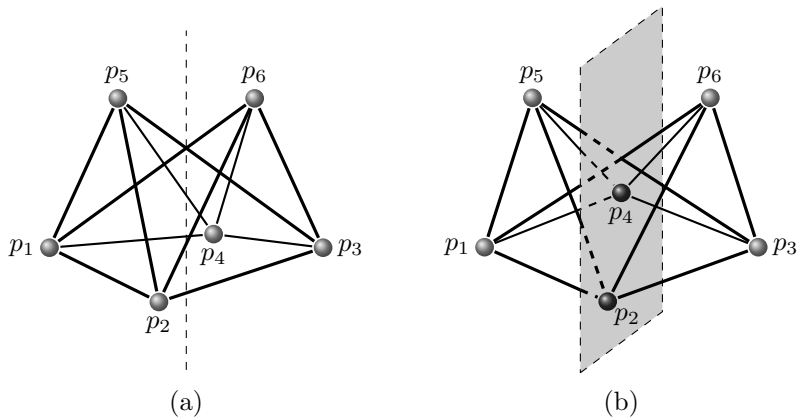
\begin{figure}[htp]
\begin{center}
\begin{tikzpicture}[very thick,scale=0.9]
\tikzstyle{every node}=[circle, fill=white, inner sep=0pt, minimum width=5pt];
\draw [dashed, thin] (0,-1.5) -- (0,3.5);
\linespread{1.0}
\node [circle, shade, ball color=black!40!white, inner sep=0pt, minimum width=7pt](p1) at (-2,0.3) {};
\node [circle, shade, ball color=black!40!white, inner sep=0pt, minimum width=7pt](p3) at (2,0.3) {};
\node [circle, shade, ball color=black!40!white, inner sep=0pt, minimum width=7pt](p2) at (-0.4,-0.5) {};
\node [circle, shade, ball color=black!40!white, inner sep=0pt, minimum width=7pt](p4) at (0.4,0.5) {};
\node [circle, shade, ball color=black!40!white, inner sep=0pt, minimum width=7pt](p5) at (-1,2.5) {};
\node [circle, shade, ball color=black!40!white, inner sep=0pt, minimum width=7pt](p6) at (1,2.5) {};
\draw (p1) -- (p2)node[rectangle, draw=white, anchor=south, below=5pt] {$p_2$};
\draw (p2) -- (p3) node[rectangle, draw=white,anchor=south, right=5pt] {$p_3$};
\draw [thick](p3) -- (p4) node[rectangle, draw=white,anchor=south, below=5pt] {$p_4$};
\draw [thick](p4) -- (p1) node[rectangle, draw=white,anchor=south, left=5pt] {$p_1$};
\draw (p5)node[rectangle, draw=white,anchor=south, above=5pt] {$p_5$} -- (p1);
\draw (p5) -- (p2);
\draw (p5) -- (p3);
\draw [thick](p5) -- (p4);
\draw (p6)node[rectangle, draw=white,anchor=south, above=5pt] {$p_6$} -- (p1);
\draw (p6) -- (p2);
\draw (p6) -- (p3);
\draw [thick](p6) -- (p4);
\node [draw=white, fill=white] (a) at (0,-2) {(a)};
\end{tikzpicture}
\hspace{0.8cm}
\begin{tikzpicture}[very thick,scale=0.9]
\tikzstyle{every node}=[circle, fill=white, inner sep=0pt, minimum width=5pt];
\filldraw[fill=black!20!white, draw=black, thin, dashed]
    (-0.6,-1.55) -- (0.6,-0.65) -- (0.6,3.85) -- (-0.6,2.95) -- (-0.6,-1.55);
\node [circle, shade, ball color=black!40!white, inner sep=0pt, minimum width=7pt](p1) at (-2,0.3) {};
\node [circle, shade, ball color=black!40!white, inner sep=0pt, minimum width=7pt](p3) at (2,0.3) {};
\node [circle, shade, ball color=black!80!white, inner sep=0pt, minimum width=7pt](p2) at (-0.2,-0.5) {};
\node [circle, shade, ball color=black!80!white, inner sep=0pt, minimum width=7pt](p4) at (0,1.1) {};
\node [circle, shade, ball color=black!40!white, inner sep=0pt, minimum width=7pt](p5) at (-1.3,2.5) {};
\node [circle, shade, ball color=black!40!white, inner sep=0pt, minimum width=7pt](p6) at (1.3,2.5) {};
\draw (p1) -- (-0.6,-0.33);
\draw [dashed] (-0.6,-0.33) -- (p2)node[rectangle, fill=black!20!white,anchor=south, below right=2pt] {$p_2$};
\draw [thick](p3) -- (p4) node[rectangle, fill=black!20!white, anchor=south, below=5pt] {$p_4$};
\draw (p2) -- (p3) node[rectangle, draw=white,anchor=south, right=5pt] {$p_3$};
\draw [dashed, thick](p4) -- (-0.6,0.86);
\draw [thick](-0.6,0.86) -- (p1) node[rectangle, draw=white,anchor=south, left=5pt] {$p_1$};
\draw (p5)node[rectangle, draw=white,anchor=south, above=5pt] {$p_5$} -- (p1);
\draw (p5) -- (-0.6,0.69);
\draw [dashed](-0.6,0.69) -- (p2);
\draw (p5) -- (-0.6,2.02);
\draw [dashed](-0.6,2.02)--(0,1.6053);
\draw (0,1.6053) -- (p3);
\draw [thick](p5) -- (-0.6,1.73);
\draw [dashed, thick] (-0.6,1.73) -- (p4);
\draw (p6) node[rectangle, draw=white,anchor=south, above=5pt] {$p_6$}-- (0,1.605);
\draw [dashed](0,1.605)--(-0.6,1.205);
\draw (-0.605,1.205) -- (p1);
\draw (p6) -- (p2);
\draw (p6) -- (p3);
\draw [thick](p6) -- (p4);
\node [draw=white, fill=white] (a) at (0,-2) {(b)};
\end{tikzpicture}
\end{center}
\vspace{-0.3cm}
\caption{Flexible Bricard octahedra with point group $\mathcal{C}_2$ (a) and $\mathcal{C}_s$ (b).}
\label{figurebricardcs}
\end{figure}

We let $G$ be the graph of the octahedron, $\mathcal{C}_{2}$ be a `half-turn' symmetry group in dimension 3, and $\Phi_a:\mathcal{C}_{2}\to \textrm{Aut}(G)$ be the homomorphism defined by
\begin{eqnarray}\Phi_a(Id)& = &id\nonumber\\ \Phi_a(C_2) &=& (1\,3)(2\,4)(5\,6)\textrm{.}\nonumber
\end{eqnarray}
Then there are three vertex orbits - represented by the vertices $1$, $2$, and $5$, for example  (see also Figure \ref{figurebricardcs}(a)). Since none of the joints $p_1$, $p_2$, and $p_5$ are fixed by the half-turn $C_2$, the number of columns of the orbit rigidity matrix $\mathbf{O}(G,p,\mathcal{C}_2)$ is
\begin{displaymath}
c=3 \cdot 3 =9\textrm{.}
\end{displaymath}
Since there are clearly six edge orbits, each of size $2$, we have
\begin{displaymath}
r=6\textrm{.}
\end{displaymath}
Finally, as shown in \cite{BS4}, we have
\begin{displaymath}
m=2\textrm{,}
\end{displaymath}
since the fully $(\mathcal{C}_2,\Phi_a)$-symmetric infinitesimal rigid motions are those that arise from translations along the $C_2$-axis and rotations about the $C_2$-axis. It follows that
\begin{displaymath}
r=6<7=c - m \textrm{,}
\end{displaymath}
so that we may conclude that $(\mathcal{C}_2,\Phi_a)$-generic realizations of the octahedral graph possess a symmetry-preserving finite flex.

This result can be generalized as follows.

\begin{theorem}\label{3dhalfturn}
Let $G$ be a graph with $|E(G)|=3|V(G)|-6$, $\mathcal{C}_2=\{Id, C_2\}$ be a half-turn symmetry group in dimension $3$, and $\Phi:\mathcal{C}_2\to \textrm{Aut}(G)$ be a homomorphism. If $j_{C_2}=b_{C_2}=0$, then $(\mathcal{C}_2,\Phi)$-generic realizations of $G$ have a symmetry-preserving mechanism.
\end{theorem}
\textbf{Proof.} Since $j_{C_2}=0$ we have $c=3\frac{|V(G)|}{2}$, and since $b_{C_2}=0$ we have $r=\frac{|E(G)|}{2}$. As mentioned above, we have  $m=2$, and hence
\begin{displaymath}r  =   3\frac{|V(G)|}{2}-3<3\frac{|V(G)|}{2}-2 =c-m \textrm{.}\nonumber
\end{displaymath}
So, by Theorems \ref{thm:kernel} and \ref{thm:flexes}, $(\mathcal{C}_2,\Phi)$-generic realizations of $G$ have a symmetry-preserving mechanism. $\square$


Next, let $G$ again be the graph of the octahedron, $\mathcal{C}_{s}$ be a `reflectional' symmetry group in dimension 3, and $\Phi_b:\mathcal{C}_{s}\to \textrm{Aut}(G)$ be the homomorphism defined by
\begin{eqnarray}\Phi_b(Id)& = &id\nonumber\\ \Phi_b(s) &=& (1\,3)(2)(4)(5\,6)\textrm{.}\nonumber
\end{eqnarray}
Then there are four vertex orbits - represented by the vertices $1$, $2$, $4$, and $5$, for example (see also Figure \ref{figurebricardcs}(b)). Since the joints $p_2$ and $p_4$ are fixed by the reflection $s$, and the joints $p_1$ and $p_5$ are not, the number of columns of the orbit rigidity matrix $\mathbf{O}(G,p,\mathcal{C}_s)$ is
\begin{displaymath}
c=2 \cdot 3 +2\cdot 2 =10\textrm{.}
\end{displaymath}
Since there are clearly six edge orbits, each of size $2$, we have
\begin{displaymath}
r=6\textrm{.}
\end{displaymath}
Finally, as shown in \cite{BS4}, we have
\begin{displaymath}
m=3\textrm{,}
\end{displaymath}
since the fully $(\mathcal{C}_s,\Phi_b)$-symmetric infinitesimal rigid motions are those that arise from translations along the mirror plane and rotations about the axis through the origin which is perpendicular to the mirror \cite{BS4}. It follows that
\begin{displaymath}
r=6<7=c - m \textrm{,}
\end{displaymath}
so that we may conclude that $(\mathcal{C}_s,\Phi_b)$-generic realizations of the octahedral graph also possess a symmetry-preserving finite flex.

More generally, we have the following result.

\begin{theorem}\label{3dmirror}
Let $G$ be a graph with $|E(G)|=3|V(G)|-6$, $\mathcal{C}_s=\{Id, s\}$ be a reflectional symmetry group in dimension $3$, and $\Phi:\mathcal{C}_s\to \textrm{Aut}(G)$ be a homomorphism. If $j_{s}>b_{s}$, then $(\mathcal{C}_s,\Phi)$-generic realizations of $G$ have a symmetry-preserving mechanism.
\end{theorem}
\textbf{Proof.} Let $|V'|$ be the number of vertex orbits of size 2 and $|E'|$ be the number of edge orbits of size 2. Then we have
\begin{eqnarray}
|E(G)| & = & 3|V(G)|-6\nonumber\\
2|E'|+b_{s} & = & 3(2|V'|+j_{s})-6\nonumber\\
|E'|+\frac{1}{2}b_{s} & = & 3|V'|+\frac{3}{2}j_{s}-3\nonumber\\
|E'|+b_{s} & = & 3|V'|+2j_{s}+\frac{1}{2}(b_{s}-j_{s})-3\nonumber\\
r & = & c+\frac{1}{2}(b_{s}-j_{s})-m \textrm{.} \nonumber
\end{eqnarray}
So if $j_{s}>b_{s}$, then $r<c-m$. Thus, by Theorems \ref{thm:kernel} and \ref{thm:flexes}, $(\mathcal{C}_s,\Phi)$-generic realizations of $G$ have a symmetry-preserving mechanism. $\square$

The flexibility of the Bricard octahedra shown in Figures \ref{figurebricardcs}(a) and (b) follow immediately from Theorems \ref{3dhalfturn} and \ref{3dmirror}.

Note that Theorems~\ref{3dhalfturn} and \ref{3dmirror} also prove the existence of a symmetry-preserving finite flex in a number of other famous and interesting frameworks in 3-space.
For example, Theorem~\ref{3dhalfturn} applies to symmetric `double-suspensions' which are frameworks that consist of an arbitrary $2n$-gon and two `cone-vertices' that are linked to each of the joints of the $2n$-gon \cite{consusp} and to some symmetric ring structures and reticulated cylinder structures like the ones analyzed in \cite{FG4}, \cite{tarnaiely}, and \cite{BS4}.  Similarly, Theorem~\ref{3dmirror}  applies to the famous `double-banana' (see \cite{gss}, for example) with mirror symmetry (with the two connecting vertices of the two `bananas' lying on the mirror), to various bipartite frameworks (such as $3$-dimensional realizations of $K_{4,6}$ with mirror symmetry, where all the joints of either one of the partite sets lie in the mirror \cite{BS4}), and to some other symmetric ring structures and reticulated cylinder structures.
%

Similar results to Theorems \ref{3dhalfturn} and \ref{3dmirror} can of course also be established for
some
other symmetry groups in dimension $3$.
In particular, it can be shown that realizations of the octahedral graph which are generic with respect to the dihedral symmetry arising from the half-turn symmetry and the mirror symmetry defined in the examples above also possess a finite flex which preserves the dihedral symmetry throughout the path (see also \cite{BS4}).
Notice that the configurations which are symmetry generic for the dihedral symmetry are not symmetry generic for the half-turn, or the mirror, alone, so this separate analysis is needed.

%

\subsubsection{Examples in dimension $d>3$}

Since the orbit counts remain remarkably simple in \emph{all} dimensions, our new method becomes particularly useful in analyzing symmetric structures in higher dimensional space. We demonstrate this by giving a very simple proof for the flexibility of the 4-dimensional cross-polytope described in \cite{Stachel4D}. As we will discuss in Section \ref{sec:metrictransf}, this example can also immediately be transferred to various other metrics.

For the graph $G$ of the $4$-dimensional cross-polytope, we have $|V(G)|=8$ and $|E(G)|=24$, and hence  $|E(G)|-(4|V(G)|-10)=24-(4\cdot 8-10)=2$. Therefore, there will always be at least two linearly independent self-stresses in any $4$-dimensional realization of $G$. However, it turns out that certain \emph{symmetric} $4$-dimensional cross-polytopes still
become flexible: Consider $4$-dimensional realizations of $G$ with dihedral symmetry $\mathcal{C}_{2v}$
which are constructed by placing two joints in each of the two perpendicular mirrors that correspond to the two reflections in $\mathcal{C}_{2v}$,
and adding their mirror reflections,
 and
 then connecting each of these
eight vertices to all other vertices of $G$ except its own mirror image. This gives rise to four vertex orbits - each of size 2. Since each mirror is a $3$-dimensional hyperplane, we have \begin{displaymath}
c=4\cdot 3 =12\textrm{.}
\end{displaymath}
Further, it is easy to check that there are $r=8$ edge orbits (4 orbits of size 4 and 4 orbits of size 2) and that $m=3$ (the fully symmetric infinitesimal rigid motions are the ones that arise from translations along the symmetry element (the `plane of rotation') of $C_2$  and rotations about the plane perpendicular to the symmetry element of $C_2$). It follows that
\begin{displaymath}
r=8<9=c - m\textrm{,}
\end{displaymath}
which, by Theorems \ref{3dhalfturn} and \ref{3dmirror}, implies that $4$-dimensional cross-polytopes which are generic with respect to this type of dihedral symmetry possess a symmetry-preserving finite flex.

Next, we provide some general counting results for frameworks with point groups $\mathcal{C}_2$ and $\mathcal{C}_s$ in dimension $d>3$ whose underlying graphs satisfy the necessary count $|E(G)|=d|V(G)|-\binom{d+1}{2}$ to be generically isostatic in dimension $d$.


Using the techniques described in \cite{BS4} it is easy to show that for a $d$-dimensional framework $(G,p)$ with point group $\mathcal{C}_2$ ($\mathcal{C}_s$), the space of fully symmetric infinitesimal translations has dimension $(d-2)$ ($d-1$, respectively) and the space of fully symmetric infinitesimal rotations has dimension $\binom{d-2}{2}+1$ ($\binom{d-1}{2}$, respectively), so that the dimension $m$ of fully symmetric infinitesimal rigid motions is equal to $d-2+\binom{d-2}{2}+1=1+\binom{d-1}{2}$ ($d-1+\binom{d-1}{2}=\binom{d}{2}$, respectively), provided that the joints of $(G,p)$ span all of $\mathbb{R}^d$. Alternatively, this can be shown by computing the dimension of the kernel of the corresponding orbit rigidity matrix $\mathbf{O}(K_n,p,S)$, where $K_n$ is the complete graph on the vertex set of $G$.

\begin{theorem}\label{4dhalf-turn}
Let $G$ be a graph with $|E(G)|=d|V(G)|-\binom{d+1}{2}$, $\mathcal{C}_2=\{Id, C_2\}$ be a half-turn symmetry group in dimension $d$, and $\Phi:\mathcal{C}_2\to \textrm{Aut}(G)$ be a homomorphism.
\begin{itemize}
\item[(i)] If $d=4$ and $b_{C_2}=0$, then $(\mathcal{C}_2,\Phi)$-generic realizations of $G$ have a symmetry-preserving mechanism;
\item[(ii)] if $d>4$ and   $j_{C_2}>\frac{b_{C_2}}{d-4}+\frac{d(d-7)+8}{2(d-4)}$, then $(\mathcal{C}_2,\Phi)$-generic realizations of $G$ have a symmetry-preserving mechanism.
\end{itemize}
\end{theorem}
\textbf{Proof.} Let $|V'|$ be the number of vertex orbits of size 2 and $|E'|$ be the number of edge orbits of size 2.

 (i)   If $d=4$, then with $b_{C_2}=0$ we have
\begin{displaymath}
r =  2|V(G)|-5< 2|V(G)|-4=2(2|V'|+j_s)-4=c-m \textrm{.}
\end{displaymath}

(ii) If $d>4$, then
\begin{eqnarray}
|E(G)| & = & d|V(G)|-\binom{d+1}{2}\nonumber\\
2|E'|+b_{C_2} & = & d(2|V'|+j_{C_2})-\binom{d+1}{2}\nonumber\\
|E'|+b_{C_2} & = & d|V'|+\frac{d}{2}j_{C_2}+\frac{1}{2}b_{C_2}-\frac{d(d+1)}{4}\nonumber\\
r & = & \big(d|V'|+(d-2)j_{C_2}\big) - \frac{d-4}{2}j_{C_2}+\frac{1}{2}b_{C_2}-\frac{d(d+1)}{4}\nonumber\\
r & = & c - \frac{d-4}{2}j_{C_2}+\frac{1}{2}b_{C_2}-\frac{d(d+1)}{4}\nonumber\\
r & = & c+\Big(\frac{1}{2}b_{C_2}- \frac{d-4}{2}j_{C_2}+\frac{d(d-7)}{4}+2\Big)-\Big(1+\binom{d-1}{2}\Big)\nonumber\\
r & = & c+\Big(\frac{1}{2}b_{C_2}- \frac{d-4}{2}j_{C_2}+\frac{d(d-7)}{4}+2\Big)-m\textrm{.}\nonumber
\end{eqnarray}
It is now easy to verify that if $j_{C_2}>\frac{b_{C_2}}{d-4}+\frac{d(d-7)+8}{2(d-4)}$, then $r<c-m$. The result now follows from Theorems \ref{thm:kernel} and \ref{thm:flexes}. $\square$
\medskip

In the formula $\Big(\frac{1}{2}b_{C_2}- \frac{d-4}{2}j_{C_2}+\frac{d(d-7)}{4}+2\Big)$, the case $d=4$,  $b_{C_2}=0$ matches part (i).
For $d=5$, the formula becomes $\frac{1}{2}(b_{C_2}- j_{C_2}-1)$, so that the count  $j_{C_2}>b_{C_2}$ guarantees the existence of a symmetry-preserving finite flex in a $(\mathcal{C}_2,\Phi)$-generic realization of $G$.


\begin{theorem}\label{4dmirror}
Let $G$ be a graph with $|E(G)|=d|V(G)|-\binom{d+1}{2}$, $\mathcal{C}_s=\{Id, s\}$ be a reflectional symmetry group in dimension $d>3$, and $\Phi:\mathcal{C}_s\to \textrm{Aut}(G)$ be a homomorphism. If $j_s>\frac{b_s}{d-2}+\frac{d(d-3)}{2(d-2)}$, then $(\mathcal{C}_s,\Phi)$-generic realizations of $G$ have a symmetry-preserving mechanism.
\end{theorem}
\textbf{Proof.} Let $|V'|$ be the number of vertex orbits of size 2 and $|E'|$ be the number of edge orbits of size 2. Then, analogously to the proof of Theorem \ref{4dhalf-turn}, we have
\begin{eqnarray}
|E(G)| & = & d|V(G)|-\binom{d+1}{2}\nonumber\\
|E'|+b_{s} & = & d|V'|+\frac{d}{2}j_{s}+\frac{1}{2}b_{s}-\frac{d(d+1)}{4}\nonumber\\
r & = & \big(d|V'|+(d-1)j_s\big) - \frac{d-2}{2}j_{s}+\frac{1}{2}b_{s}-\frac{d(d+1)}{4}\nonumber\\
r & = & c+\Big(\frac{1}{2}b_{s}- \frac{d-2}{2}j_{s}+\frac{d(d-3)}{4}\Big)-\Big(d-1+\binom{d-1}{2}\Big)\nonumber\\
r & = & c+\Big(\frac{1}{2}b_{s}- \frac{d-2}{2}j_{s}+\frac{d(d-3)}{4}\Big)-m\textrm{.}\nonumber
\end{eqnarray}
It is now easy to verify that if $j_s>\frac{b_s}{d-2}+\frac{d(d-3)}{2(d-2)}$, then $r<c-m$. The result now follows from Theorems \ref{thm:kernel} and \ref{thm:flexes}. $\square$

\subsubsection{Coning the counts for $d\geq 3$}
\label{sec:coning}

Coning of a framework embeds the framework in one higher dimension, then adds a new vertex attached to all previous vertices.  This is a technique that takes the counts for a finite flex in a generic framework in dimension $d$, to the counts for a finite flex of the cone in dimension $d+1$ \cite{WWcones}.   It is natural to consider how this impacts the counts of the orbit matrix.   First - the symmetry-coning will transfer the symmetry group by adding the new vertex on the normal to the origin into the new dimension, extending the axes of any rotations, and the mirrors of any reflections in the symmetry group into the larger space.   With this placement, the cone vertex is fixed by the entire symmetry group, the prior edges and vertices have the same orbits, and the new edges from the cone vertex to the prior vertices have orbits corresponding to their end points - that is one for each of the prior vertex orbits.  This process transfers the counts which guaranteed symmetry-preserving finite flexes of the original graph in dimension $d$ to counts on the cone which guarantee symmetry-preserving finite flexes of the cone graph in dimension $d+1$.  Combined, the orbit matrices and coning provide a powerful tool to construct flexible polytopes in all dimensions.
%

\indent Consider, for example, the graph $G$ of the octahedron which is generically isostatic in dimension $3$. If we cone $G$ (i.e., we add a new vertex to $G$ and connect it to each of the original vertices of $G$), then we obtain a new graph $G^*$ which is generically isostatic in dimension $4$. If we now realize $G^*$ `generically' with half-turn symmetry in $4$-space so that no bar is fixed by the half-turn, and the cone vertex is the only vertex that lies on the ($2$-dimensional) half-turn axis, then the resulting framework possesses a symmetry-preserving mechanism (the orbit counts are $r=6+3=9<c-m=3\cdot 4+2-4=10$).
\\\indent In general, if $d>3$, and we repeatedly cone the octahedral graph  $(d-3)$ times, then the resulting graph $G^*$ is generically isostatic in dimension $d$. However, if we realize $G^*$
   `generically' with half-turn symmetry in dimension $d$ so that the $(d-3)$ cone vertices all lie on the $(d-2)$-dimensional half-turn axis, then we have $r=6+3\cdot(d-3)+\binom{d-3}{2}$ (for the $6$ edge orbits of $G$, the $(d-3)$ connections from the cone vertices to each of the vertices of $G$, and the $\binom{d-3}{2}$ edges in the half-turn axis for the complete graph on the cone vertices), $c=d\cdot 3+(d-2)(d-3)$, and $m=1+\binom{d-1}{2}$, and hence $r=(c-m)-1$. Thus, we obtain flexible polytopes with half-turn symmetry in all dimensions in this way.
 \\\indent  Analogously, based on the realization of the octahedron with point group $\mathcal{C}_s$ in Figure \ref{figurebricardcs}(b), we may construct flexible polytopes with mirror symmetry in all dimensions. Of course we may also symmetrically cone other flexible polytopes (e.g., the cross-polytope) to produce new flexible polytopes in the next higher dimension (see also Section \ref{sec:metrictransf}).

\subsection{Detection of finite flexes from the rank of the orbit rigidity matrix}

We have seen that the count $r\geq c - m$ is a necessary condition for a symmetric framework $(G,p)\in \mathscr{R}_{(G,S,\Phi)}$ to have no fully $(S,\Phi)$-symmetric infinitesimal flex (Theorem \ref{symMaxthm}). However, it is not a sufficient condition, so that if $(G,p)$ satisfies the count $r\geq c - m$, we need to compute the actual rank of the orbit rigidity matrix $\mathbf{O}(G,p,S)$ to determine whether $(G,p)$ has a fully $(S,\Phi)$-symmetric infinitesimal flex.

Alternatively, one could use group representation theory to block-diagonalize the rigidity matrix $\mathbf{R}(G,p)$ as described in \cite{KG2,BS2}, and then compute the rank of the submatrix block which corresponds to the trivial irreducible representation of $S$. This approach, however, requires significantly more work than simply finding the rank of the orbit rigidity matrix.

In the following, we demonstrate the simplicity of our new method for detecting finite flexes via the rank of the orbit rigidity matrix with the help of some examples.

\subsubsection{Examples in dimension $2$}

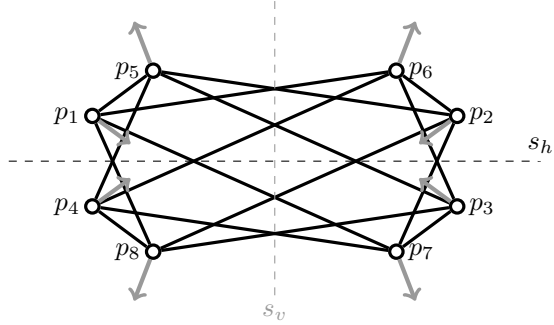
\begin{figure}[htp]
\begin{center}
\begin{tikzpicture}[very thick,scale=1]
\tikzstyle{every node}=[circle, draw=black, fill=white, inner sep=0pt, minimum width=5pt];
        \path (-2.4,-0.6) node (p1) [label =  left: $p_4$] {} ;
        \path (2.4,-0.6) node (p2) [label =  right: $p_3$] {} ;
        \path (2.4,0.6) node (p3) [label =  right: $p_2$] {} ;
        \path (-2.4,0.6) node (p4) [label =  left: $p_1$] {} ;
        \path (-1.6,-1.2) node (p5) [label = left: $p_8$] {} ;
        \path (1.6,-1.2) node (p6) [label =  right: $p_7$] {} ;
        \path (1.6,1.2) node (p7) [label =  right: $p_6$] {} ;
        \path (-1.6,1.2) node (p8) [label =  left: $p_5$] {} ;
        \draw (p1) -- (p5);
        \draw (p1) -- (p6);
        \draw (p1) -- (p7);
        \draw (p1) -- (p8);
        \draw (p2) -- (p5);
        \draw (p2) -- (p6);
        \draw (p2) -- (p7);
        \draw (p2) -- (p8);
        \draw (p3) -- (p5);
        \draw (p3) -- (p6);
        \draw (p3) -- (p7);
        \draw (p3) -- (p8);
        \draw (p4) -- (p5);
        \draw (p4) -- (p6);
        \draw (p4) -- (p7);
        \draw (p4) -- (p8);
        \draw [->, ultra thick, black!40!white] (p4) -- (-1.9,0.23);
        \draw [->, ultra thick, black!40!white] (p3) -- (1.9,0.23);
        \draw [->, ultra thick, black!40!white] (p2) -- (1.9,-0.23);
        \draw [->, ultra thick, black!40!white] (p1) -- (-1.9,-0.23);
        \draw [->, ultra thick, black!40!white] (p8) -- (-1.85,1.85);
        \draw [->, ultra thick, black!40!white] (p7) -- (1.85,1.85);
        \draw [->, ultra thick, black!40!white] (p6) -- (1.85,-1.85);
        \draw [->, ultra thick, black!40!white] (p5) -- (-1.85,-1.85);
           \draw [dashed, thin, black!40!white] (0,-2.2)node[draw=white,anchor=south] {$s_v$} -- (0,2.2);
\draw [dashed, thin] (-3.5,0) -- (3.5,0)node[draw=white,anchor=south] {$s_h$};
       \end{tikzpicture}
\end{center}
\vspace{-0.3cm}
\caption{A fully $(\mathcal{C}_{2v},\Phi)$-symmetric infinitesimal flex of a framework in $\mathscr{R}_{(K_{4,4},\mathcal{C}_{2v},\Phi)}$.}
\label{k44fl}
\end{figure}

As a first example, we consider the complete bipartite graph $K_{4,4}$ with partite sets $\{1,2,3,4\}$ and $\{5,6,7,8\}$. Note that the graph $K_{4,4}$ is generically rigid in dimension $2$. Moreover, any $2$-dimensional realization of $K_{4,4}$ has three linearly independent self-stresses since $|E(K_{4,4})|-(2|V(K_{4,4})|-3)= 16-(2\cdot 8 -3)=3$. However, as we will see, under certain symmetry conditions, $2$-dimensional realizations of $K_{4,4}$ become flexible.\\\indent Let $\mathcal{C}_{2v}=\{Id, C_2, s_h, s_v\}$ be the dihedral symmetry group in dimension 2 generated by the reflections $s_h$ and $s_v$ whose corresponding mirror lines are the $x$-axis and $y$-axis, respectively, and let $\Phi:\mathcal{C}_{2v}\to \textrm{Aut}(K_{4,4})$  be the homomorphism defined by
\begin{eqnarray}\Phi(Id)& = &id\nonumber\\ \Phi(C_2) &=& (1\,3)(2\,4)(5\,7)(6\,8)\nonumber\\ \Phi(s_h) &=& (1\,4)(2\,3)(5\,8)(6\,7)\nonumber\\ \Phi(s_v) &=& (1\,2)(3\,4)(5\,6)(7\,8)\textrm{.}\nonumber
\end{eqnarray}
A framework $(K_{4,4},p)$ in the set $\mathscr{R}_{(K_{4,4},\mathcal{C}_{2v},\Phi)}$ is depicted in Figure \ref{k44fl}.
Let's first compute the \emph{size} of the orbit rigidity matrix $\mathbf{O}(K_{4,4},p,\mathcal{C}_{2v})$. There are two vertex orbits - represented by the vertices $1$ and $5$, for example - and also four edge orbits - represented by the edges $\{1,5\}$, $\{1,6\}$, $\{1,7\}$, and $\{1,8\}$, for example. Since $m$ is clearly equal to zero, and $c_1=\textrm{dim}\big(U(p_1)\big)=2$ and $c_5=\textrm{dim}\big(U(p_5)\big)=2$, we have \begin{displaymath}r=4=c-m\textrm{.}\end{displaymath}
So, to determine whether $(K_{4,4},p)$ possesses a fully $(\mathcal{C}_{2v},\Phi)$-symmetric infinitesimal flex, we need to set up the matrix $\mathbf{O}(K_{4,4},p,\mathcal{C}_{2v})$ explicitly. If we denote $(p_1)^T=(a,b)$ and $(p_5)^T=(c,d)$, then we have
\begin{eqnarray}
\mathbf{O}(K_{4,4},p,\mathcal{C}_{2v}) &=&\bordermatrix{ & 1 & 5 \cr \{1,5\} & \big(p_1-p_5\big)^T & \big(p_5-p_1\big)^T \cr \{1,6\} & \big(p_1-s_v(p_5)\big)^T & \big(p_5-s_v^{-1}(p_1)\big)^T \cr  \{1,7\} & \big(p_1-C_2(p_5)\big)^T & \big(p_5-C_2^{-1}(p_1)\big)^T \cr \{1,8\} &  \big(p_1-s_h(p_5)\big)^T & \big(p_5-s_h^{-1}(p_1)\big)^T}\nonumber\\
& = & \bordermatrix{& &  \cr  & (a-c, b-d) & (c-a, d-b) \cr  & (a+c, b-d) & (c+a, d-b) \cr  & (a+c, b+d) & (c+a, d+b)\cr  & (a-c, b+d) & (c-a, d+b)} \textrm{.}\nonumber
\end{eqnarray}


It is now easy to see that for any choice of $a,b,c$, and $d$, the rows of $\mathbf{O}(K_{4,4},p,\mathcal{C}_{2v})$ are linearly dependent (the sum of the first and third row vector minus the sum of the second and fourth row vector is equal to zero). Thus, the kernel of $\mathbf{O}(K_{4,4},p,\mathcal{C}_{2v})$ is non-trivial, and since $m=0$, it follows that any realization of $K_{4,4}$ in $\mathscr{R}_{(K_{4,4},\mathcal{C}_{2v},\Phi)}$ possesses a fully $(\mathcal{C}_{2v},\Phi)$-symmetric infinitesimal flex. (By computing an element in the kernel of $\mathbf{O}(K_{4,4},p,\mathcal{C}_{2v})$ explicitly, it can be seen that all the velocity vectors of this infinitesimal flex are orthogonal to the conic determined by the joints of $(K_{4,4},p)$ (see also Figure \ref{k44fl})). By Theorems \ref{3dhalfturn} and \ref{3dmirror}, it follows that $(\mathcal{C}_{2v},\Phi)$-generic realizations of $K_{4,4}$ possess a symmetry-preserving finite flex. This flex is also known as `Bottema's!
  mechanism \cite{bottema}.


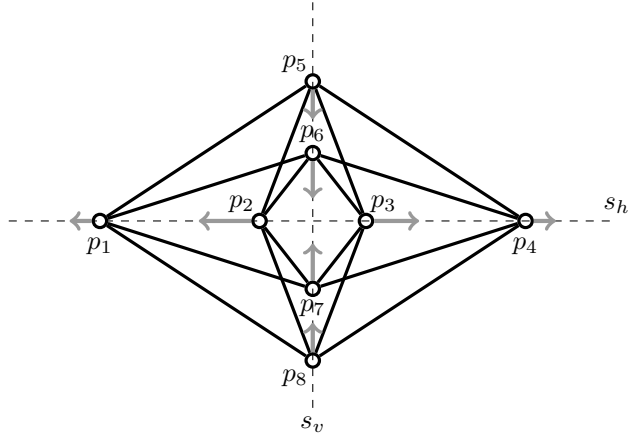
\begin{figure}[htp]
\begin{center}
\begin{tikzpicture}[very thick,scale=1]
\tikzstyle{every node}=[circle, draw=black, fill=white, inner sep=0pt, minimum width=5pt];
\linespread{1.0};
        \path (-2.8,0) node (p1) [label =  below: $p_1$] {} ;
        \path (-0.7,0) node (p2) [label =  above left: $p_2$] {} ;
        \path (0.7,0) node (p3) [label =  above right: $p_3$] {} ;
        \path (2.8,0) node (p4) [label =  below: $p_4$] {} ;
        \path (0,-1.85) node (p5) [label = below left: $p_8$] {} ;
        \path (0,-0.9) node (p6) {} ;
        \path (0,0.9) node (p7)  {} ;
        \path (0,1.85) node (p8) [label =  above left: $p_5$] {} ;
        \node[rectangle, draw=white](labelp6) at (0,-1.15) {$p_7$};
        \node[rectangle, draw=white](labelp7) at (0,1.15) {$p_6$};
        \draw (p1) -- (p5);
        \draw (p1) -- (p6);
        \draw (p1) -- (p7);
        \draw (p1) -- (p8);
        \draw (p2) -- (p5);
        \draw (p2) -- (p6);
        \draw (p2) -- (p7);
        \draw (p2) -- (p8);
        \draw (p3) -- (p5);
        \draw (p3) -- (p6);
        \draw (p3) -- (p7);
        \draw (p3) -- (p8);
        \draw (p4) -- (p5);
        \draw (p4) -- (p6);
        \draw (p4) -- (p7);
        \draw (p4) -- (p8);
              \draw [dashed, thin] (0,-2.9) node[draw=white,anchor=south] {$s_v$} -- (0,2.9);
\draw [dashed, thin] (-4,0) -- (4,0)node[draw=white,anchor=south] {$s_h$} ;
        \draw [->, ultra thick, black!40!white] (p4) -- (3.2,0);
        \draw [->, ultra thick, black!40!white] (p3) -- (1.4,0);
        \draw [->, ultra thick, black!40!white] (p2) -- (-1.5,-0);
        \draw [->, ultra thick, black!40!white] (p1) -- (-3.2,-0);
        \draw [->, ultra thick, black!40!white] (p8) -- (0,1.35);
        \draw [->, ultra thick, black!40!white] (p7) -- (0,0.3);
        \draw [->, ultra thick, black!40!white] (p6) -- (0,-0.3);
        \draw [->, ultra thick, black!40!white] (p5) -- (0,-1.35);
       \end{tikzpicture}
\end{center}
\vspace{-0.3cm}
\caption{A fully $(\mathcal{C}_{2v},\Psi)$-symmetric infinitesimal flex of a framework in $\mathscr{R}_{(K_{4,4},\mathcal{C}_{2v},\Psi)}$.}
\label{fig:k44vertfix}
\end{figure}

Figure \ref{fig:k44vertfix} shows another type of realization of $K_{4,4}$ in the plane with dihedral symmetry. This framework is an element of the set $\mathscr{R}_{(K_{4,4},\mathcal{C}_{2v},\Psi)}$, where $\Psi:\mathcal{C}_{2v}\to \textrm{Aut}(K_{4,4})$  is the homomorphism defined by
\begin{eqnarray}\Psi(Id)& = &id\nonumber\\ \Psi(C_2) &=& (1\,4)(2\,3)(5\,8)(6\,7)\nonumber\\ \Psi(s_h) &=& (1)(2)(3)(4)(5\,8)(6\,7)\nonumber\\ \Psi(s_v) &=& (1\,4)(2\,3)(5)(6)(7)(8)\textrm{.}\nonumber
\end{eqnarray}
The vertex orbits are represented by the set $\{1,2,5,6\}$, for example, and we have $c_i=\textrm{dim}\big(U(p_i)\big)=1$ for $i=1,2,5,6$. Further, the edge orbits are represented by the set $\big\{\{1,5\},\{1,6\},\{2,5\},\{2,6\}\big\}$, for example. Thus, the orbit count is again \begin{displaymath}r=4=c-m\textrm{.}\end{displaymath}
To determine whether $(K_{4,4},p)$ possesses a fully $(\mathcal{C}_{2v},\Psi)$-symmetric infinitesimal flex, we set up the orbit matrix $\mathbf{O}(K_{4,4},p,\mathcal{C}_{2v})$. With $(p_1)^T=(a,0)$, $(p_2)^T=(b,0)$, $(p_5)^T=(0,c)$, and $(p_6)^T=(0,d)$, we have
\begin{eqnarray}
& &\mathbf{O}(K_{4,4},p,\mathcal{C}_{2v})\nonumber\\
& =& \bordermatrix{ & 1 & 2 & 5 & 6\cr \{1,5\} & \big(p_1-p_5\big)^T\binom{1}{0} & 0 & \big(p_5-p_1\big)^T\binom{0}{1} & 0\cr \{1,6\} & \big(p_1-p_6\big)^T\binom{1}{0} & 0 &  0  & \big(p_6-p_1\big)^T\binom{0}{1}  \cr  \{2,5\} &  0  & \big(p_2-p_5\big)^T\binom{1}{0} & \big(p_5-p_2\big)^T\binom{0}{1} & 0\cr \{2,6\} &   0  & \big(p_2-p_6\big)^T\binom{1}{0} &  0  & \big(p_6-p_2\big)^T\binom{0}{1} }\nonumber\\
& = & \bordermatrix{& & & & \cr  & a & 0 & c & 0 \cr  & a & 0 & 0 & d \cr  & 0 & b & c & 0\cr  & 0 & b & 0 & d} \textrm{.}\nonumber
\end{eqnarray}
Clearly, the rows of $\mathbf{O}(K_{4,4},p,\mathcal{C}_{2v})$ are linearly dependent. Thus, analogously as above, we may conclude that $(\mathcal{C}_{2v},\Psi)$-generic realizations of $K_{4,4}$ also possess a symmetry-preserving finite flex.

\subsubsection{Examples in dimension $d>2$}

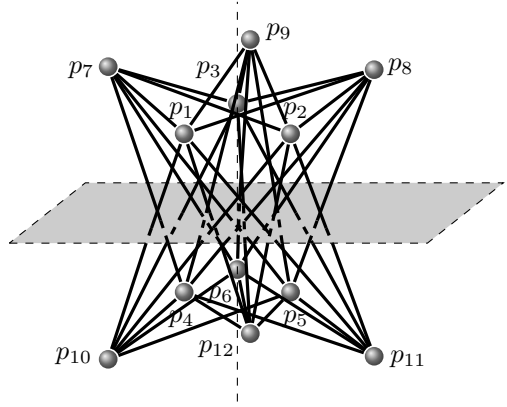
\begin{figure}[htp]
\begin{center}
\begin{tikzpicture}[very thick,scale=1]
\tikzstyle{every node}=[circle, fill=white, inner sep=0pt, minimum width=5pt];
\filldraw[fill=black!20!white, draw=black, thin, dashed]
    (-3,-0.4) -- (2.5,-0.4) -- (3.5,0.4) -- (-2,0.4) -- (-3,-0.4);
\node [circle, shade, ball color=black!40!white, inner sep=0pt, minimum width=7pt](p1) at (-0.7,1.05) {};
\node [circle, shade, ball color=black!40!white, inner sep=0pt, minimum width=7pt](p2) at (0.7,1.05) {};
\node [circle, shade, ball color=black!40!white, inner sep=0pt, minimum width=7pt](p3) at (0,1.45) {};
\node [circle, shade, ball color=black!40!white, inner sep=0pt, minimum width=7pt](p4) at (-0.7,-1.05) {};
\node [circle, shade, ball color=black!40!white, inner sep=0pt, minimum width=7pt](p5) at (0.7,-1.05) {};
\node [circle, shade, ball color=black!40!white, inner sep=0pt, minimum width=7pt](p6) at (0,-0.75) {};
\node [circle, shade, ball color=black!40!white, inner sep=0pt, minimum width=7pt](p7) at (-1.7,1.94) {};
\node [circle, shade, ball color=black!40!white, inner sep=0pt, minimum width=7pt](p8) at (1.8,1.9) {};
\node [circle, shade, ball color=black!40!white, inner sep=0pt, minimum width=7pt](p9) at (0.17,2.3) {};
\node [circle, shade, ball color=black!40!white, inner sep=0pt, minimum width=7pt](p10) at (-1.7,-1.94) {};
\node [circle, shade, ball color=black!40!white, inner sep=0pt, minimum width=7pt](p11) at (1.8,-1.9) {};
\node [circle, shade, ball color=black!40!white, inner sep=0pt, minimum width=7pt](p12) at (0.17,-1.6) {};
\node [rectangle,draw=white, fill=white] (a) at (-0.73,1.36) {$p_1$};
\node [rectangle,draw=white, fill=white] (a) at (0.77,1.36) {$p_2$};
\node [rectangle,draw=white, fill=white] (a) at (-0.37,1.9) {$p_3$};
\node [rectangle,draw=white, fill=white] (a) at (-0.73,-1.4) {$p_4$};
\node [rectangle,draw=white, fill=white] (a) at (0.77,-1.4) {$p_5$};
\node [rectangle,draw=white, fill=white] (a) at (-0.22,-1.08) {$p_6$};
\node [rectangle,draw=white, fill=white] (a) at (-2.05,1.9) {$p_7$};
\node [rectangle,draw=white, fill=white] (a) at (2.15,1.94) {$p_8$};
\node [rectangle,draw=white, fill=white] (a) at (0.55,2.4) {$p_9$};
\node [rectangle,draw=white, fill=white] (a) at (-2.15,-1.9) {$p_{10}$};
\node [rectangle,draw=white, fill=white] (a) at (2.25,-1.94) {$p_{11}$};
\node [rectangle,draw=white, fill=white] (a) at (-0.27,-1.8) {$p_{12}$};
\draw[dashed, thin] (0,-2.5)--(0,2.8);
\draw (p1) -- (p7);
\draw (p1) -- (p8);
\draw (p1) -- (p9);
\draw  (p1) --(-1.13,-0.2) ;
\draw [dashed](-1.13,-0.2)--(-1.2,-0.4);
\draw (-1.2,-0.4)--(p10);
\draw (p1) --(0.36,-0.2) ;
 \draw [dashed] (0.36,-0.2) -- (0.53,-0.4);
\draw (0.53,-0.4) -- (p11);
\draw (p1) -- (-0.36,0);
\draw [dashed]  (-0.36,0) --(-0.22,-0.4);
\draw (-0.22,-0.4) -- (p12);
\draw (p2) -- (p7);
\draw (p2) -- (p8);
\draw (p2) -- (p9);
\draw (p2) -- (-0.32,-0.2);
\draw [dashed] (-0.32,-0.2)--(-0.48,-0.4);
\draw (-0.48,-0.4) -- (p10);
\draw  (p2) -- (1.168,-0.2);
\draw [dashed](1.168,-0.2)--(1.25,-0.4);
\draw (1.25,-0.4) -- (p11);
\draw (p2) -- (0.49,0);
\draw [dashed] (0.49,0) --(0.4,-0.4);
\draw (0.4,-0.4) -- (p12);
\draw (p3) -- (p7);
\draw (p3) -- (p8);
\draw (p3) -- (p9);
\draw  (p3) -- (-0.722,0);
\draw [dashed] (-0.722,0) -- (-0.93,-0.4);
\draw (-0.93,-0.4) -- (p10);
\draw  (p3) -- (0.77,0);
\draw [dashed] (0.77,0) -- (0.99,-0.4);
\draw (0.99,-0.4) -- (p11);
\draw  (p3) -- (0.08,0.1);
\draw [dashed] (0.08,0.1) -- (0.1,-0.4);
\draw (0.1,-0.4) -- (p12);
\draw (p4) -- (-0.91,-0.4) ;
\draw [dashed] (-0.91,-0.4) -- (-1.01,-0.1);
\draw  (-1.01,-0.1) -- (p7);
\draw (p4) -- (-0.165,-0.4) ;
\draw [dashed] (-0.165,-0.4) -- (0.1,-0.1);
\draw  (0.1,-0.1) -- (p8);
\draw (p4) -- (-0.515,-0.4) ;
\draw [dashed] (-0.515,-0.4) -- (-0.41,0);
\draw  (-0.41,0) -- (p9);
\draw (p4) -- (p10);
\draw (p4) -- (p11);
\draw (p4) -- (p12);
\draw (p5) -- (0.19,-0.4) ;
\draw [dashed] (0.19,-0.4) -- (-0.07,-0.1);
\draw  (-0.07,-0.1) -- (p7);
\draw (p5) -- (0.94,-0.4) ;
\draw [dashed] (0.94,-0.4) -- (1.05,-0.1);
\draw  (1.05,-0.1) -- (p8);
\draw (p5) -- (0.59,-0.4) ;
\draw [dashed] (0.59,-0.4) -- (0.533,0);
\draw  (0.533,0) -- (p9);
\draw (p5) -- (p10);
\draw (p5) -- (p11);
\draw (p5) -- (p12);
\draw (p6) -- (-0.22,-0.4) ;
\draw [dashed] (-0.22,-0.4) -- (-0.48,0);
\draw  (-0.48,0) -- (p7);
\draw (p6) -- (0.22,-0.4) ;
\draw [dashed] (0.22,-0.4) -- (0.5,0);
\draw  (0.5,0) -- (p8);
\draw (p6) -- (0.01,-0.4) ;
\draw [dashed] (0.01,-0.4) -- (0.03,0.2);
\draw  (0.03,0.2) -- (p9);
\draw (p6) -- (p10);
\draw (p6) -- (p11);
\draw (p6) -- (p12);
\end{tikzpicture}
\end{center}
\vspace{-0.3cm}
\caption{A framework $(K_{6,6},p)$ in $\mathscr{R}_{(K_{6,6},\mathcal{C}_{3h},\Phi)}$.}
\label{k66flex}
\end{figure}

We first describe an example of a flexible framework in $3$-space which can be thought of as the
$3$-dimensional analog of Bottema's mechanism in the plane. Consider the complete bipartite graph $K_{6,6}$ with partite sets $\{1,\ldots,6\}$ and $\{7,\ldots,12\}$. This graph is generically rigid in dimension $3$. Moreover, every $3$-dimensional realization of $K_{6,6}$ possesses at least 6 linearly independent self-stresses, because $|E(K_{6,6})|-(3|V(K_{6,6})|-6)= 36-(3\cdot 12 -6)=6$. However, with the help of the orbit rigidity matrix it is easy to see that certain \emph{symmetric} $3$-dimensional realizations of $K_{6,6}$ become flexible.\\\indent Let $\mathcal{C}_{3h}$ be the symmetry group in dimension 3 which is generated by the reflection $s$ whose corresponding mirror plane is the $xy$-plane and the $3$-fold rotation $C_3$ whose corresponding rotational axis is the $z$-axis. Further, we let $\Phi:\mathcal{C}_{3h}\to \textrm{Aut}(K_{6,6})$  be the homomorphism defined by
\begin{eqnarray}\Phi(C_3) &=& (1\,2\,3)(4\,5\,6)(7\,8\,9)(10\,11\,12)\nonumber\\ \Phi(s) &=& (1\,4)(2\,5)(3\,6)(7\,10)(8\,11)(9\,12)\textrm{,}\nonumber
\end{eqnarray}
and let $(K_{6,6},p)$ be a $(\mathcal{C}_{3h},\Phi)$-generic realization of $K_{6,6}$ (see also Figure \ref{k66flex}).

We first compute the \emph{size} of the orbit rigidity matrix $\mathbf{O}(K_{6,6},p,\mathcal{C}_{3h})$. There are two vertex orbits - represented by the vertices $1$ and $7$, for example - and six edge orbits - represented by the edges $\{1,i\}$, $i=7,\ldots,12$, for example.  Since rotations about the $C_3$-axis are clearly the only infinitesimal rigid motions that are fully $(\mathcal{C}_{3h},\Phi)$-symmetric, we have $m=1$. Since we also have $c_1=c_7=3$, it follows that \begin{displaymath}r=6>5=c-m\textrm{.} \end{displaymath}
So we detect a fully $(\mathcal{C}_{3h},\Phi)$-symmetric self-stress, but no fully $(\mathcal{C}_{3h},\Phi)$-symmetric infinitesimal flex of $(K_{6,6},p)$ with this count. If we want to show that $(K_{6,6},p)$ possesses a fully $(\mathcal{C}_{3h},\Phi)$-symmetric infinitesimal flex, we need to prove that the rank of $\mathbf{O}(K_{6,6},p,\mathcal{C}_{3h})$ is at most 4. We assume wlog that $(p_1)^T=(\sqrt{3},0,1)$ and $(p_7)^T=(a,b,d)$. Then we have
\begin{eqnarray}
& & \mathbf{O}(K_{6,6},p,\mathcal{C}_{3h})\nonumber\\
&=&\bordermatrix{ & 1 & 7 \cr \{1,7\} & \big(p_1-p_7\big)^T & \big(p_7-p_1\big)^T \cr \{1,C_3(7)\} & \big(p_1-C_3(p_7)\big)^T & \big(p_7-C_3^{2}(p_1)\big)^T \cr  \{1,C_3^2(7)\} & \big(p_1-C_3^2(p_7)\big)^T & \big(p_7-C_3(p_1)\big)^T \cr \{1,s(7)\} &  \big(p_1-s(p_7)\big)^T & \big(p_7-s(p_1)\big)^T \cr  \{1,sC_3(7)\} & \big(p_1-sC_3(p_7)\big)^T & \big(p_7-sC_3^2(p_1)\big)^T \cr \{1,sC_3^2(7)\} &  \big(p_1-sC_3^2(p_7)\big)^T & \big(p_7-sC_3(p_1)\big)^T}\nonumber\\
& = & \bordermatrix{& & \cr  & (\sqrt{3}-a, -b , 1-c) & (a-\sqrt{3}, b, c-1) \cr  & (\sqrt{3}+\frac{a+\sqrt{3}b}{2}, \frac{-\sqrt{3}a-b}{2}, 1-c) & (a+\frac{\sqrt{3}}{2}, b+\frac{3}{2}, c-1) \cr  & (\sqrt{3}+\frac{a-\sqrt{3}b}{2}, \frac{\sqrt{3}a-b}{2}, 1-c) & (a+\frac{\sqrt{3}}{2}, b-\frac{3}{2}, c-1) \cr  & (\sqrt{3}-a, -b, 1+c) & (a-\sqrt{3}, -b, c+1) \cr  & (\sqrt{3}+\frac{a+\sqrt{3}b}{2}, \frac{-\sqrt{3}a-b}{2}, 1+c) & (a+\frac{\sqrt{3}}{2}, b+\frac{3}{2}, c+1) \cr  & (\sqrt{3}+\frac{a-\sqrt{3}b}{2}, \frac{\sqrt{3}a-b}{2}, 1+c) & (a+\frac{\sqrt{3}}{2}, b-\frac{3}{2}, c+1)} \textrm{.}\nonumber
\end{eqnarray}
Clearly, the equation \begin{displaymath}\omega^T \mathbf{O}(K_{6,6},p,\mathcal{C}_{3h})=0
\end{displaymath}
is satisfied for the linearly independent vectors $\omega_1^T=\left(\begin{array} {cccccc}0 & 1 & -1 & 0 & -1 & 1\end{array}\right)$ and $\omega_2^T=\left(\begin{array} {cccccc}1 & 0 & -1 & -1 & 0 & 1\end{array}\right)$. Thus, $(\mathcal{C}_{3h},\Phi)$-generic realizations of $K_{6,6}$ indeed possess a symmetry-preserving finite flex.

In general, for any dimension $d>3$, we may construct $d$-dimensional realizations of the complete bipartite graph $K_{2d,2d}$ with point group $\mathcal{C}_{dh}$ by choosing one vertex from each of the two partite sets of $K_{2d,2d}$ as representatives for the vertex orbits and placing them `generically' off the mirror plane corresponding to the reflection in $\mathcal{C}_{dh}$ and also off the rotational axis corresponding to the $d$-fold rotation in $\mathcal{C}_{dh}$. This gives rise to two vertex orbits (each of size $2d$) and $2d$ edge orbits. Since the infinitesimal rigid motions corresponding to rotations about the $d$-fold axis will always be fully-symmetric with respect to this $\mathcal{C}_{dh}$ symmetry, the orbit counts will always detect a fully-symmetric self-stress, but no fully-symmetric infinitesimal flex for these frameworks. However, these frameworks can be shown to possess a symmetry-preserving finite flex analogously as above by computing the actual ranks of the corresponding orbit matrices.\\\indent Note that for these types of realizations of $K_{2d,2d}$ with $\mathcal{C}_{dh}$ symmetry, the geometry of quadric surfaces  (see \cite{BolRot, Wbipartite}, for example) can be used to predict the existence of a fully-symmetric infinitesimal flex  and therefore also the rank properties of the corresponding orbit matrices.




\section{Fully symmetric self-stresses}
\label{sec:stress}
\subsection{The kernel of $\mathbf{O}(G,p,S)^T$}

In this section, we show that for a framework $(G,p)\in \mathscr{R}_{(G,S,\Phi)}$,
the cokernel of the orbit rigidity matrix $\mathbf{O}(G,p,S)$
is the space of all fully $(S,\Phi)$-symmetric self-stresses of $(G,p)$, restricted to the corresponding set  $\mathscr{O}_{E(G)}$ of representatives for the edge orbits  $S(e)=\{\Phi(x)(e)|\,x\in S\}$ of $G$ (Theorem \ref{thm:coker}). To prove this result, we need the following two lemmas:

\begin{lemma}
\label{lemma1}
Let $G$ be a graph, $S$ be a symmetry group in dimension $d$, $\Phi:S\to \textrm{Aut}(G)$ be a homomorphism,  and $(G,p)\in \mathscr{R}_{(G,S,\Phi)}$.
Further, let $S(e)=\{\Phi(x)(e)|\,x\in S\}$ be an edge orbit of $G$ whose representative $e=\{a,x(b)\}$ is an edge whose end-vertices lie in distinct vertex orbits of $G$.
Then there exist respective bases $\mathscr{B}_a$ and $\mathscr{B}_b$ for $U(p_a)$ and $U(p_b)$ (whose coordinate column vectors relative to the canonical basis form the $d\times c_a$ matrix $\mathbf{M}_a$ and the $d\times c_b$ matrix $\mathbf{M}_b$, respectively), a scalar $\alpha_e\in \mathbb{R}$, and two invertible $d\times d$ matrices $\mathbf{A}$ and $\mathbf{B}$ such that
\begin{equation}\label{eq:lem11}
\sum_{j:\{a,j\}\in S(e)}(p_a-p_j)^T=\frac{1}{\alpha_e}\Big(\big(p_a-x(p_b)\big)^T\mathbf{M}_{a}, 0,\ldots, 0\Big)\mathbf{A}
 \end{equation}
and
\begin{equation}\label{eq:lem12}
\sum_{j:\{b,j\}\in S(e)}(p_b-p_j)^T=\frac{1}{\alpha_e}\Big(\big(p_b-x^{-1}(p_a)\big)^T\mathbf{M}_{b}, 0,\ldots, 0\Big)\mathbf{B}\textrm{.}
 \end{equation}
\end{lemma}
\textbf{Proof.} Let $\{Id=y_0,y_1,\ldots, y_t\}$ be the stabilizer $S_{p_a}=\{x\in S|\, x(p_a)=p_a\}$ of $p_a$, $\mathbf{Y}_l$ be the matrix that represents $y_l$ with respect to the canonical basis of $\mathbb{R}^d$ for each $l$, and $\alpha_e=|S_{p_a}\cap S_{x(p_b)}|$. Then we have
\begin{displaymath}
 \sum_{j:\{a,j\}\in S(e)}(p_a-p_j)^T=\frac{1}{\alpha_e}\sum_{l=0}^{t}\Big(\mathbf{Y}_l\big(p_a- x(p_b)\big)\Big)^T\textrm{,}
\end{displaymath}
because  \begin{displaymath}\Big(y_l(\{p_a,x(p_b)\})\Big)_{l=0,\ldots,t}\end{displaymath} is the family of those bars of $(G,p)$ whose corresponding edges lie in $S(e)$ and are incident with $a$, and because $p_a-x(p_b)=y_l(p_a-x(p_b))$ if and only if $y_l$ is an element of the coset $S_{p_a}\cap S_{x(p_b)}$ of $S_{p_a}$.\\\indent
Consider the matrix representation $H:S_{p_a}\to GL(\mathbb{R}^d)$ that assigns to each $y_l\in S_{p_a}$ the $d\times d$ matrix $\mathbf{Y}_l$ which represents $y_l$ with respect to the canonical basis of $\mathbb{R}^d$. By definition, the $H$-invariant subspace $\mathbb{V}$ of $\mathbb{R}^d$ corresponding to the trivial irreducible representation of $H$  is the space $U(p_a)$.
 Thus, by the Great Orthogonality Theorem, there exists an orthogonal $d\times d$ matrix  of basis transformation $\mathbf{M}$ (i.e., $\mathbf{M}^{-1}=\mathbf{M}^T$) such that
\begin{displaymath}
\sum_{l=0}^{t}\mathbf{M}^{-1}\mathbf{Y}_l\mathbf{M}= \mathbf{M}^{-1}\Big(\sum_{l=0}^{t}\mathbf{Y}_l\Big)\mathbf{M}=\left(\begin{array}{ccc|ccc}t & & & & & \\ & \ddots & & & \mathbf{0} & \\ & & t & &  & \\ \hline & & & & & \\ & \mathbf{0} & & & \mathbf{0}  \end{array}\right)\textrm{,}
\end{displaymath}
where the first $c_a$ column vectors of $\mathbf{M}$ are the coordinate vectors of a basis for $\mathbb{V}= U(p_a)$ relative to the canonical basis of $\mathbb{R}^d$. We let $\mathbf{M}_a$ be such that
\begin{displaymath}
 \mathbf{M}=\frac{1}{t}\left(\begin{array}{c|ccc} & \vdots & &\vdots \\ \mathbf{M}_a & * & \hdots & * \\ & \vdots & &\vdots  \end{array}\right)\textrm{.}
\end{displaymath}
Then, for $\mathbf{A}=\mathbf{M}^T$, we have
\begin{eqnarray}
 \alpha_e\sum_{j:\{a,j\}\in S(e)}(p_a-p_j)^T & =& \Big(\big(\sum_{l=0}^{t}\mathbf{Y}_l\big)(p_a-x(p_b))\Big)^T \nonumber\\
& = & \Big(\mathbf{M}\left(\begin{array}{ccc|ccc}t & & & & & \\ & \ddots & & & \mathbf{0} & \\ & & t & &  & \\ \hline & & & & & \\ & \mathbf{0} & & & \mathbf{0}  \end{array}\right)\mathbf{M}^{T}(p_a-x(p_b))\Big)^T \nonumber\\
& = & \Big(\mathbf{M}\left(\begin{array}{ccc} & \mathbf{M}^{T}_a & \\ \hline \hdots & 0 &\hdots\\ & \vdots &\\\hdots & 0 &\hdots  \end{array}\right)(p_a-x(p_b))\Big)^T\nonumber\\
 & = & \Big(\big(p_a-x(p_b)\big)^T\mathbf{M}_{a}, 0,\ldots, 0\Big)\mathbf{A}\textrm{.}\nonumber
\end{eqnarray}
This proves (\ref{eq:lem11}).\\\indent
Note that if we can show that $|S_{p_a}\cap S_{x(p_b)}|=|S_{p_b}\cap S_{x^{-1}(p_a)}| $, then the proof of (\ref{eq:lem12}) proceeds completely analogously to the one of (\ref{eq:lem11}). \\\indent
Consider the map \begin{displaymath} \psi:\left\{\begin{array}{ccc}S_{p_a}\cap S_{x(p_b)} & \to  & S_{p_b}\cap S_{x^{-1}(p_a)}\\ y & \mapsto & x^{-1}yx\end{array} \right.\textrm{.}\end{displaymath}
We show that $\psi$ is well-defined. Note that $S_{x(p_b)}=xS_{p_b}x^{-1}$ and $S_{x^{-1}(p_a)}=x^{-1}S_{p_a}x$.   Thus, $y\in S_{p_a}\cap S_{x(p_b)}$ if and only if $y\in S_{p_a}$ and $y=x\hat{y}x^{-1}$ with $\hat{y}\in S_{p_b}$. We have $\psi(y)=x^{-1}yx=x^{-1}(x\hat{y}x^{-1})x=\hat{y}\in S_{p_b}$, and hence $\psi(y)\in S_{p_b}\cap S_{x^{-1}(p_a)}$. Since $\psi$ is clearly bijective, we indeed have $|S_{p_a}\cap S_{x(p_b)}|=|S_{p_b}\cap S_{x^{-1}(p_a)}| $. This gives the result.
 $\square$

\begin{lemma}
\label{lemma2}
Let $G$ be a graph, $S$ be a symmetry group in dimension $d$, $\Phi:S\to \textrm{Aut}(G)$ be a homomorphism,  and $(G,p)\in \mathscr{R}_{(G,S,\Phi)}$.
Further, let $S(e)$ be an edge orbit of $G$ whose representative $e=\{a,x(a)\}$ is an edge whose end-vertices lie in the same vertex orbit of $G$.
 Then there exists a basis $\mathscr{B}_a$ for $U(p_a)$ (whose coordinate column vectors relative to the canonical basis form the $d\times c_a$ matrix $\mathbf{M}_a$), a scalar $\alpha_e\in \mathbb{R}$, and an invertible $d\times d$ matrix $\mathbf{A}$ such that
\begin{equation}\label{eq:lem22}
\sum_{j:\{a,j\}\in S(e)}(p_a-p_j)^T=\frac{1}{\alpha_e}\Big(\big(2p_a-x(p_a)-x^{-1}(p_a)\big)^T\mathbf{M}_{a}, 0,\ldots, 0\Big)\mathbf{A}\textrm{.}
\end{equation}
\end{lemma}
\textbf{Proof.}  Let $\{Id=y_0,y_1,\ldots, y_t\}$ be the stabilizer $S_{p_a}=\{x\in S|\, x(p_a)=p_a\}$ of $p_a$, and let $\mathscr{F}_1$ and $\mathscr{F}_2$ be the families of bars of $(G,p)$ defined by
\begin{eqnarray}\mathscr{F}_1 &= &\Big(y_l(\{p_a,x(p_a)\})\Big)_{l=0,\ldots,t}\nonumber\\
\mathscr{F}_2 &=& \Big(y_l(\{p_a,x^{-1}(p_a)\})\Big)_{l=0,\ldots,t}\nonumber\end{eqnarray}
The range of the families $\mathscr{F}_1$ and $\mathscr{F}_2$ are the bars that correspond to the summands in the left hand side of equation (\ref{eq:lem22}). Note that we either have $\mathscr{F}_1=\mathscr{F}_2$ or $\mathscr{F}_1\cap\mathscr{F}_2=\emptyset$ ($\mathscr{F}_1=\mathscr{F}_2$ if and only if there exists $y_l\in S_{p_a}$ such that $y_l(x(p_a))=x^{-1}(p_a)$). Moreover, we have $|S_{p_a}\cap S_{x(p_a)}|=|S_{p_a}\cap S_{x^{-1}(p_a)}|$ since, by a similar argument as in the proof of Lemma \ref{lemma1}, the map \begin{displaymath} \psi:\left\{\begin{array}{ccc}S_{p_a}\cap S_{x(p_a)} & \to  & S_{p_a}\cap S_{x^{-1}(p_a)}\\ y & \mapsto & x^{-1}yx\end{array} \right.\end{displaymath}
is well-defined and bijective.
\\\indent
Suppose first that $\mathscr{F}_1\cap\mathscr{F}_2=\emptyset$. Then we have
\begin{eqnarray}\label{eqlem1}
 & &\sum_{j:\{a,j\}\in S(e)}(p_a-p_j)^T \nonumber\\ &= & \frac{1}{\alpha_e}\Big(\sum_{l=0}^{t}\Big(\mathbf{Y}_l\big(p_a- x(p_a)\big)\Big)^T + \sum_{l=0}^{t}\Big(\mathbf{Y}_l\big(p_a- x^{-1}(p_a)\big)\Big)^T\Big)\nonumber\\
 &=& \frac{1}{\alpha_e}\Big(\big(\sum_{l=0}^{t}\mathbf{Y}_l\big)\big(2p_a-x(p_a)-x^{-1}(p_a)\big)\Big)^T\textrm{,}
\end{eqnarray}
where  $\mathbf{Y}_l$ is the matrix that represents $y_l$ with respect to the canonical basis of $\mathbb{R}^d$ for each $l$, and $\alpha_e=|S_{p_a}\cap S_{x(p_a)}|=|S_{p_a}\cap S_{x^{-1}(p_a)}|$.
\\\indent Suppose next that $\mathscr{F}_1=\mathscr{F}_2$. Then
\begin{displaymath}
\sum_{l=0}^{t}\Big(\mathbf{Y}_l\big(p_a- x(p_a)\big)\Big)^T=\sum_{l=0}^{t}\Big(\mathbf{Y}_l\big(p_a- x^{-1}(p_a)\big)\Big)^T\textrm{,}
\end{displaymath}
and hence
\begin{equation}\label{eqlem2}
 \sum_{j:\{a,j\}\in S(e)}(p_a-p_j)^T=\frac{1}{\alpha_e}\Big(\big(\sum_{l=0}^{t}\mathbf{Y}_l\big)\big(2p_a-x(p_a)-x^{-1}(p_a)\big)\Big)^T\textrm{,}
\end{equation}
where $\alpha_e=2|S_{p_a}\cap S_{x(p_a)}|$.\\\indent
Now, by the same argument as in the proof of Lemma \ref{lemma1}, it follows from equations (\ref{eqlem1}) and (\ref{eqlem2})  that for the scalars $\alpha_e$ defined above and the matrices  $\mathbf{M}_a$ and $\mathbf{A}$  defined in Lemma \ref{lemma1}, equation (\ref{eq:lem22}) holds. $\square$

\begin{theorem}
\label{thm:coker}
Let $G$ be a graph with $V(G)=\{1,\ldots,n\}$, $S$ be a symmetry group in dimension $d$, $\Phi:S\to \textrm{Aut}(G)$ be a homomorphism, $\mathscr{O}_{V(G)}=\{1,\ldots, k\}$ and $\mathscr{O}_{E(G)}=\{e_1,\ldots, e_r\}$  be sets of representatives for the vertex orbits $S(i)=\{\Phi(x)(i)|\,x\in S\}$ and edge orbits  $S(e)=\{\Phi(x)(e)|\,x\in S\}$ of $G$, respectively, and $(G,p)\in \mathscr{R}_{(G,S,\Phi)}$. If the scalars  $\alpha_{e_i}$, $i=1,\ldots, r$, and the bases for the spaces $U(p_i)$, $i=1,\ldots, k$, are defined as in Lemmas \ref{lemma1} and \ref{lemma2}, then $\tilde{\omega}\in \mathbb{R}^{r}$ is an element of the kernel of $\mathbf{O}(G,p,S)^T$ if and only if  \begin{displaymath}\overline{\omega}=\left(\begin{array}{c}\alpha_{e_1}(\tilde{\omega})_1\\\vdots\\ \alpha_{e_r}(\tilde{\omega})_r \end{array}\right)\end{displaymath}  is the restriction $\omega|_{\mathscr{O}_{E(G)}}$ of a fully $(S,\Phi)$-symmetric self-stress $\omega$ of  $(G,p)$ to $\mathscr{O}_{E(G)}$.
\end{theorem}
\textbf{Proof.} We let $\mathbf{O}_{i,a}$ be the $c_a$-dimensional row vector which consists of those components of the $i$th row of $\mathbf{O}(G,p,S)$ that correspond to the vertex $a\in V(G)$. We further let $\overline{\mathbf{O}}_{i,a}$ be the $d$-dimensional row vector $(\mathbf{O}_{i,a}, 0,\ldots,0)$.\\\indent
Suppose first that $\overline{\omega}$ is the restriction $\omega|_{\mathscr{O}_{E(G)}}$ of a fully $(S,\Phi)$-symmetric self-stress $\omega$ of  $(G,p)$. Then for every vertex $a=1,\ldots, k$, we have
\begin{displaymath}
 \sum_{i=1}^{r}\sum_{j:\{a,j\}\in S(e_i)}(\overline{\omega})_i(p_a-p_j)^T=0^T\textrm{.}
\end{displaymath}
By (\ref{eq:lem11}), (\ref{eq:lem12}),  and (\ref{eq:lem22}), for every vertex $a=1,\ldots, k$, we have
\begin{eqnarray}
 \sum_{i=1}^{r}\sum_{j:\{a,j\}\in S(e_i)}(\overline{\omega})_i(p_a-p_j)^T &=& \sum_{i=1}^{r}(\overline{\omega})_i \sum_{j:\{a,j\}\in S(e_i)} (p_a-p_j)^T \nonumber\\
& = &  \sum_{i=1}^{r}(\tilde{\omega})_i \big(\overline{\mathbf{O}}_{i,a}\mathbf{A}\big)\nonumber\\
& = & \big( \sum_{i=1}^{r}(\tilde{\omega})_i \overline{\mathbf{O}}_{i,a}\big)\mathbf{A}\textrm{,}\nonumber
\end{eqnarray}
where $\mathbf{A}$ is defined as in Lemmas \ref{lemma1} and  \ref{lemma2}. (In particular, if $p_a$ is not fixed by any non-trivial symmetry operation in $S$, then $\mathbf{A}$ is the $d \times d$ identity matrix and $\overline{\mathbf{O}}_{i,a}=\mathbf{O}_{i,a}$.) Since $\mathbf{A}$ is invertible, it follows that
\begin{displaymath}
 \sum_{i=1}^{r}(\tilde{\omega})_i \overline{\mathbf{O}}_{i,a}=0^T\textrm{,}
\end{displaymath}
and hence
\begin{displaymath}
 \sum_{i=1}^{r}(\tilde{\omega})_i \mathbf{O}_{i,a}=0^T\textrm{.}
\end{displaymath}

Conversely, if $\tilde{\omega}$ is an element of the kernel of $\mathbf{O}(G,p,S)^T$, then for every vertex $a=1,\ldots, k$, we have
\begin{displaymath}
 \sum_{i=1}^{r}(\tilde{\omega})_i \mathbf{O}_{i,a}=0^T\textrm{.}
\end{displaymath}
and hence, by the same argument as above,
\begin{displaymath}
 \sum_{i=1}^{r}\sum_{j:\{a,j\}\in S(e_i)}(\overline{\omega})_i(p_a-p_j)^T=0^T\textrm{.}
\end{displaymath}
Moreover, for every $x\in S$, we have
\begin{eqnarray}
 \sum_{i=1}^{r}\sum_{j:\{a,j\}\in S(e_i)}(\overline{\omega})_i\big(\mathbf{X}(p_a-p_j)\big)^T & = & \big(\sum_{i=1}^{r}\sum_{j:\{a,j\}\in S(e_i)}(\overline{\omega})_i(p_a-p_j)^T\big)\mathbf{X}^T\nonumber\\
& = & 0^T\mathbf{X}^T\nonumber\\
& = & 0^T\nonumber\textrm{,}
\end{eqnarray}
where $\mathbf{X}$ is the matrix that represents $x$ with respect to the canonical basis of $\mathbb{R}^d$. This completes the proof. $\square$

\subsection {Fully symmetric tensegrities}
\label{sec:SymmetricTensegrities}
It is natural to investigate how stressed symmetric frameworks can convert to tensegrity frameworks, with cables (members that can get shorter but not longer), struts (members that can get longer but not shorter) as well as bars (whose length is fixed) \cite{RW1}.  A number of the classical
 tensegrity frameworks are based on symmetric frameworks, and the Robert Connelly's web site \cite{conweb}  permits an interactive exploration of a range of examples of symmetric tensegrity frameworks.

We give a few basic definitions and translate some standard results to the symmetric setting.

A {\em tensegrity graph} $\hat G$ has a partition of the edges of $G$ into three disjoint parts $E(G)=E_+(G)\cup E_-(G)\cup E_0(G)$.  $E_+(G)$ are the edges that are {\em cables}, $E_-(G)$ are the {\em struts} and $E_0(G)$ are the {\em bars}.  For a tensegrity framework $(\hat G, p)$, a {\em proper self-stress} is a self-stress on the underlying framework $(G,p)$ with the added condition that   $\omega_{ij}\geq 0, \{i,j\}\in E_+$, $\omega_{ij}\leq 0, \{i,j\}\in E_-$ \cite{RW1}.

Given a symmetric framework $(G,p)\in \mathscr{R}_{(G,S,\Phi)}$, it is possible to use a fully $(S,\Phi)$-symmetric self-stress on the bar and joint framework $(G,p)$ to investigate both the infinitesimal rigidity of $(G,p)$,  and  the infinitesimal rigidity of an associated \emph{fully symmetric} tensegrity framework $(\hat{G},p)$ (i.e., the edges of an edge orbit are either all cables, or all struts, or all bars),  with all members with $\omega_{ij}>0$ as cables and all members with $\omega_{ij}<0$ as struts.

The standard result for the infinitesimal rigidity of such frameworks is

\begin{theorem} [Roth, Whiteley \cite{RW1}] A tensegrity framework $(\hat{G},p)$ is infinitesimally rigid if and only if the underlying bar framework $({G},p)$ is infinitesimally rigid as a bar and joint framework and $({G},p)$ has a self-stress which has  $\omega_{ij}>0$ on cables and  $\omega_{ij}<0$ on struts.
\end{theorem}

Translated in terms of the orbit matrix for a symmetric framework, this says:

\begin{cor}  A fully symmetric tensegrity framework $(\hat{G},p)$ is infinitesimally rigid if and only if the underlying bar framework $({G},p)\in \mathscr{R}_{(G,S,\Phi)}$ is infinitesimally rigid as a bar and joint framework and the orbit matrix $\mathbf{O}({G},p,S)$ has a self-stress which has  $\omega_{ij}>0$ on cables and  $\omega_{ij}<0$ on struts.
\end{cor}

Often, tensegrity frameworks are built which are rigid, but not infinitesimally rigid \cite{conrigener,conweb}.  Clearly, the underlying framework  $({G},p)$ is not generic (where rigidity is equivalent to infinitesimal rigidity), so $({G},p)$ has some self-stress.   The results of Connelly \cite{conrigener} tell us that $(\hat{G},p)$ has a non-zero proper self-stress.

\begin{theorem} [Connelly \cite{conrigener} Theorem 3] Let $(\hat{G},p)$ be a rigid tensegrity framework with a cable or strut.  Then there is a proper self-stress in the tensegrity framework (with  $\omega_{ij}>0$ on cables and  $\omega_{ij}<0$ on struts).
\end{theorem}\label{thm:Con}

Given a fully symmetric rigid tensegrity framework $(\hat{G},p)$, we can show that the guaranteed self-stress can be chosen to be fully symmetric.

\begin{cor}  Let $(\hat{G},p)$ be a fully symmetric rigid tensegrity framework with a cable or strut whose underlying bar framework lies in $\mathscr{R}_{(G,S,\Phi)}$.  Then there is a fully $(S,\Phi)$-symmetric non-zero proper self-stress in the tensegrity framework (with  $\omega_{ij}>0$ on cables and  $\omega_{ij}<0$ on struts).
\end{cor}
\noindent \textbf{Proof.}  By Theorem~\ref{thm:Con}, there is a non-zero proper self-stress in $(\hat{G},p)$.  We want to symmetrize this self-stress.    For each element of the group $x\in S$, and each edge $\{i,j\}$, we have the coefficient $\omega_{x(i,j)}$ of the corresponding element of the orbit.  If we add over all elements of the group, this is a finite sum, and we have a combined coefficient $ \omega_{S(i,j)} $.  It is a direct computation to confirm that these coefficients are a self-stress (the sum of self-stresses is a self-stress) and that they form a fully symmetric self-stress.  Since the original stress was proper on a fully symmetric tensegrity framework, all the $\omega_{x(i,j)}$ for a given edge have the same sign, so there is no cancelation.  We conclude that this is the required non-zero proper fully symmetric self-stress.
$\square$

The following example illustrates how these pieces fit together in the layers of symmetry-preserving finite flexes in symmetry generic configurations, non-symmetric finite flexes for symmetry generic configurations for a larger group, and fully symmetric stresses giving rigidity for an even larger group.

\begin{examp} Consider the graph and frameworks illustrated in Figure~\ref{fig:cube}.

\begin{figure}
    \begin{center}
  \subfigure[] {\includegraphics [width=.35\textwidth]{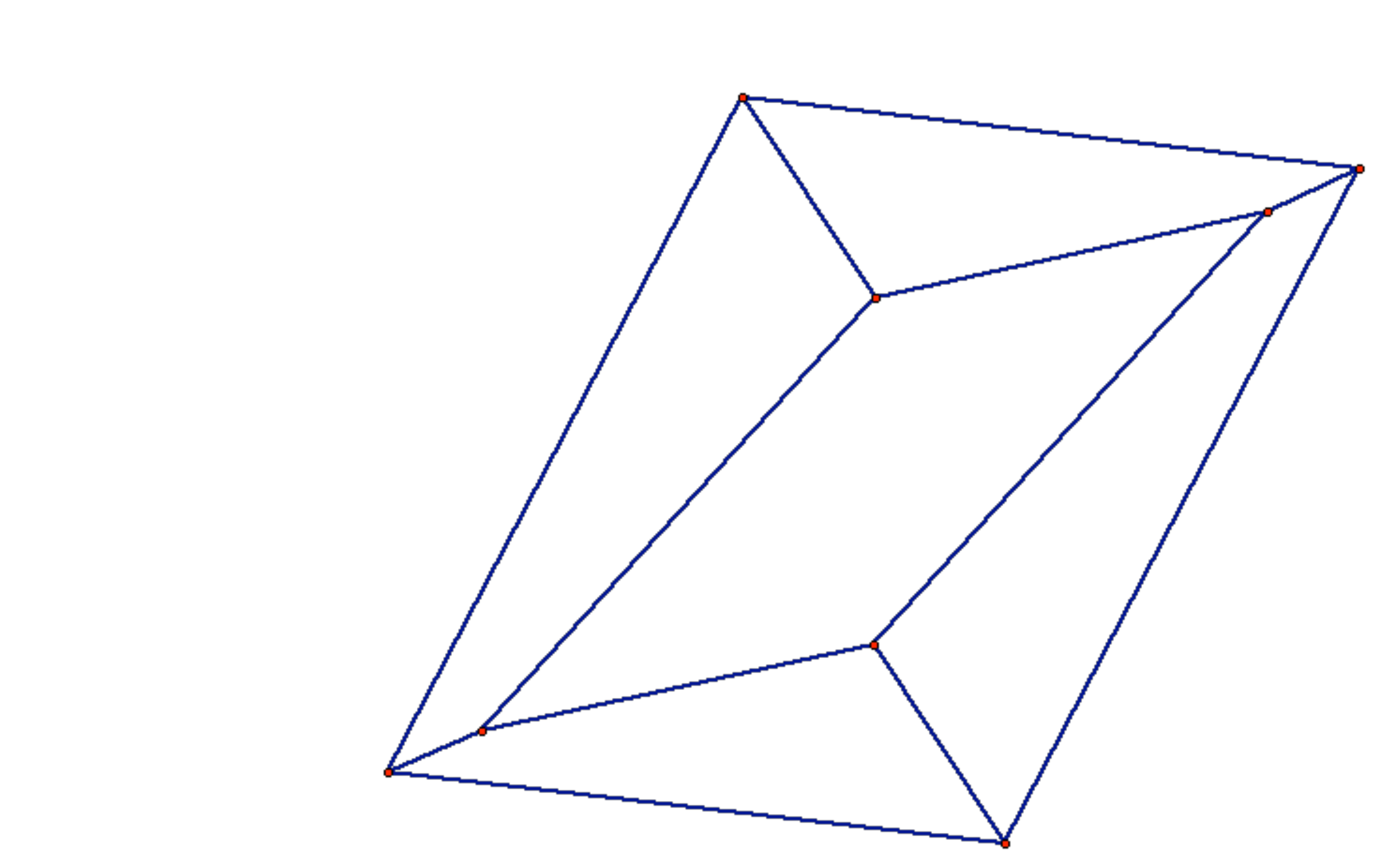}} \quad
  \subfigure[] {\includegraphics [width=.28\textwidth]{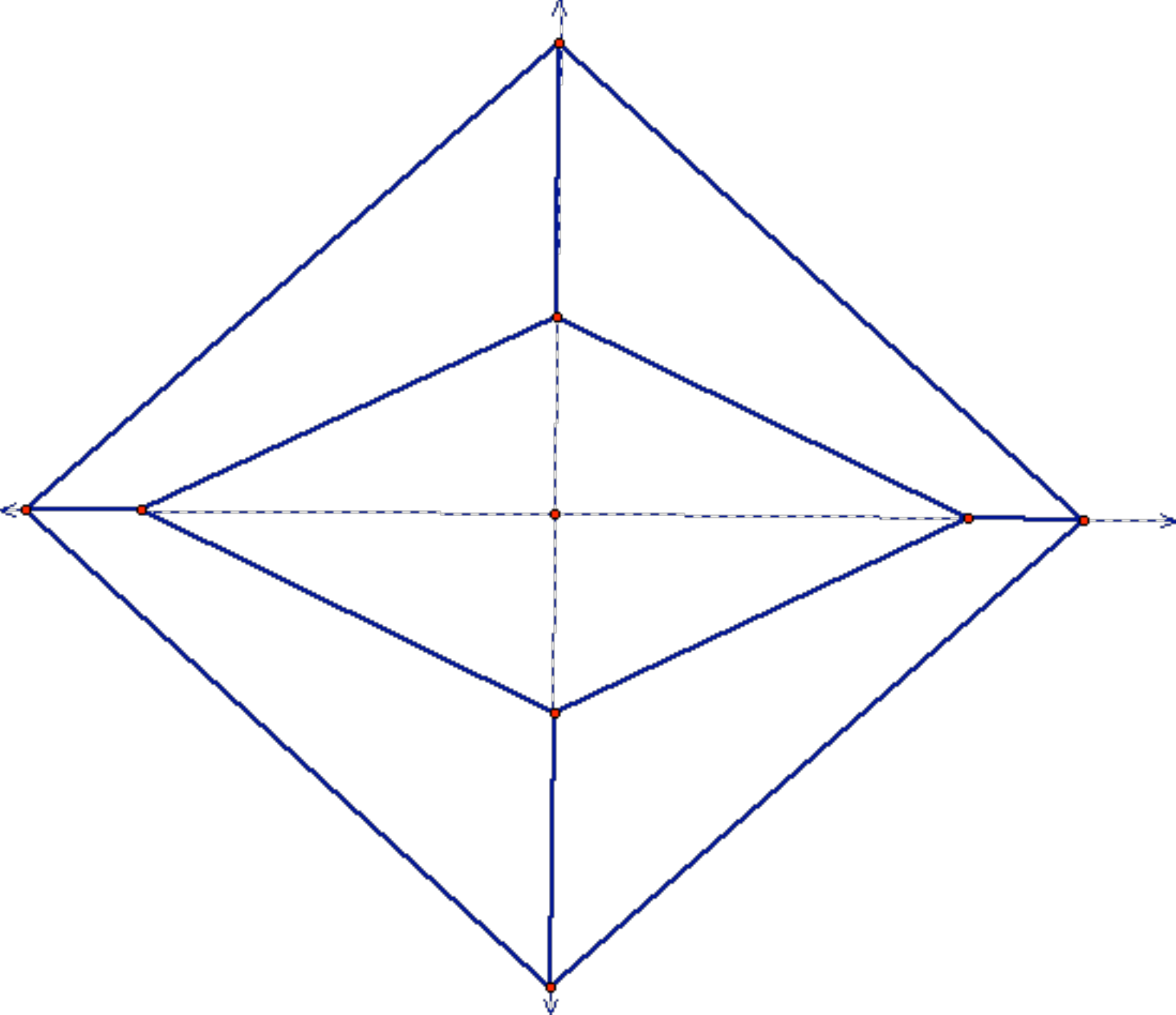}}\quad \quad
  \subfigure[] {\includegraphics [width=.25\textwidth]{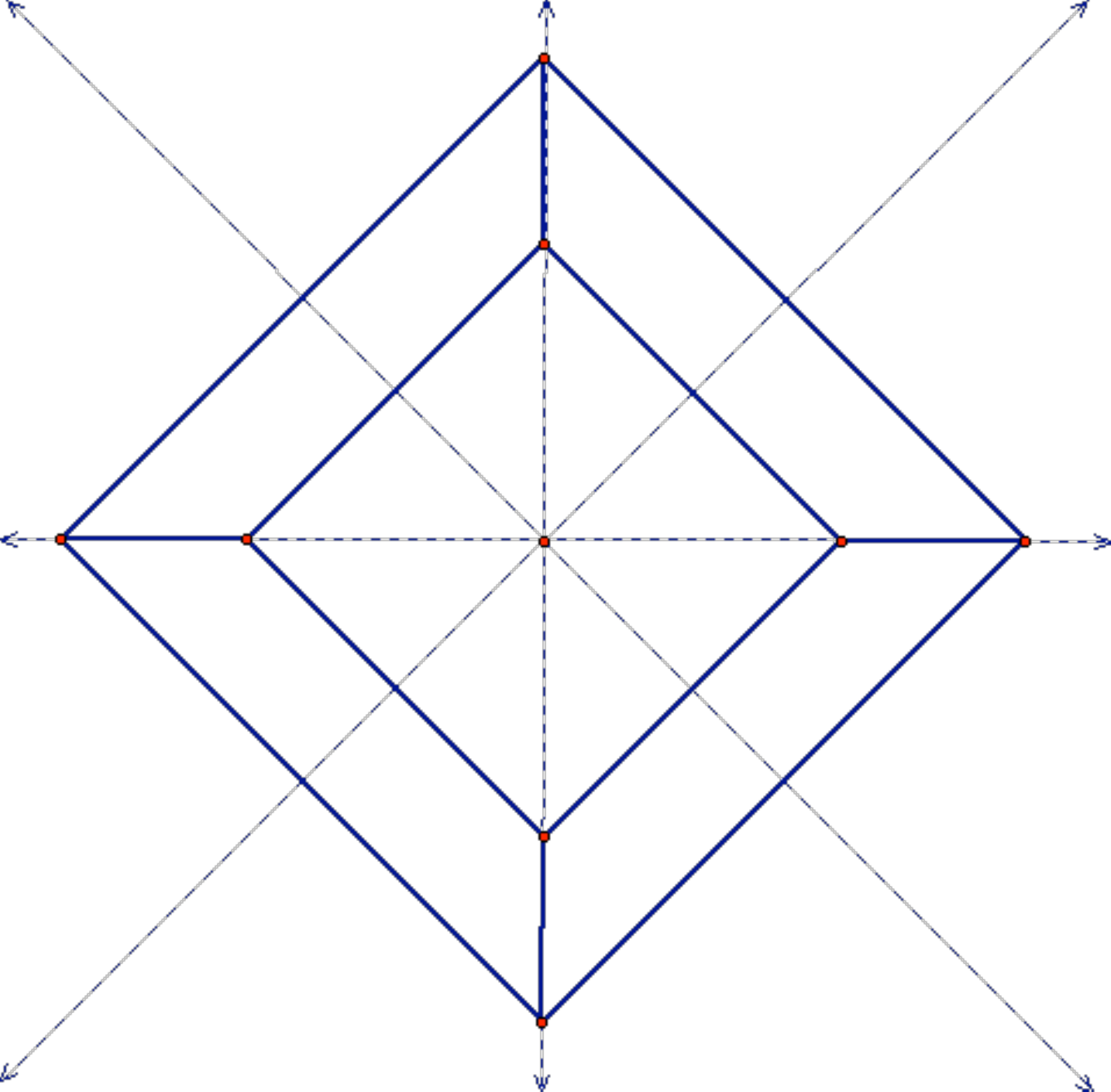}}
      \end{center}
    \caption{We give three plane configurations for the edge graph of a cube. In (a) there is $C_{2}$ symmetry, in (b) there is $C_{2v}$ symmetry, and in (c) there is $C_{4v}$ symmetry.  (a) has a symmetry-preserving finite flex, (b) has a finite flex which breaks the mirror symmetry, and (c) has a fully symmetric self-stress which makes it rigid.}
    \label{fig:cube}
    \end{figure}

\begin{enumerate}
\item Figure ~\ref{fig:cube}(a) shows the graph $G$ realized at a generic configuration with $C_{2}$ symmetry.  The counts for the rank of the orbit matrix are:  $r=6$, $c=8$ and $m=1$.  This guarantees a symmetry-preserving finite flex.
	While it is not immediate, the standard result for such planar graphs \cite{CWMax} shows that this framework only has a self-stress if it is the projection of a plane faced polyhedron - which this is not (there is no consistent line of intersection of the outside quadrilateral and the inside face).   The symmetry-preserving finite flex is also the flex guaranteed by the basic generic counts:   $|E(G)|=12$, $2|V(G)|=16$ and $|E(G)|=12< 16 -3= 2|V(G)|-3$.

\item Figure~\ref{fig:cube}(b)  shows a symmetry generic configuration for $C_{2v}$.  The new counts are $r=4$, $c=4$ and $m=0$.  The corresponding orbit matrix counts to be independent - and in fact the framework still has no self-stress.  The framework is still not the projection of a plane faced polyhedron.  There is a finite flex, but it is not symmetry-preserving for this $C_{2v}$ symmetry (only for $C_{2}$).

\item Figure~\ref{fig:cube}(c) shows a symmetry generic configuration for $C_{4v}$.  The revised counts are:  $r=3$, $c=2$ and $m=0$. We are guaranteed a fully symmetric self-stress.  (One can also see this as the projection of a plane faced cube-line polyhedron).   It is now possible that this is rigid (and remains rigid with cables and struts following the signs of the self-stress).   With cables on the interior of the framework, this is a spider web, and the approach of \cite{conrigener} just works to confirm that these are rigid (though not infinitesimally rigid).
\end{enumerate}

\end{examp}

As the example illustrates, and the many structures on \cite{conweb} confirm, a fully symmetric self-stress can be the way of forming a rigid tensegrity framework which is too undercounted to be infinitesimally rigid.

We conjecture that a further analog of Connelly's Theorem also holds, and that the basic proof can be symmetry adapted:

\begin{conj}  Let $(\hat{G},p)$ be a fully symmetric  tensegrity framework with a cable or strut which has no symmetry-preserving finite flex.  Then there is a fully symmetric non-zero proper self-stress in the tensegrity framework (with  $\omega_{ij}>0$ on cables and  $\omega_{ij}<0$ on struts).
\end{conj}

\section{Further work}
\label{sec:further}
As mentioned in the introduction, the analysis of the orbit matrix opens up a number of questions which are analogs of the previous work for the standard rigidity matrix.   The following samples are not exhaustive, and we find new possibilities keep opening up for us as we continue to work with the tools and reflect on the possibilities.

\subsection{Necessary and sufficient conditions for a full rank orbit matrix}
An important question for  the standard rigidity matrix has been deriving necessary and sufficient conditions on the graph for the rigidity matrix to be of full rank (generic rigidity), or independent, or to have a self-stress.  The most famous example is Laman's Theorem characterizing generic rigidity in the plane \cite{Laman}.  Within the context of symmetric frameworks, there are generalizations for key plane groups ($C_3$, $C_s$, and $C_2$) presented in \cite{BS4, BS3}.  With these combinatorial calculations come fast algorithms for verifying the generic rigidity.

It is natural to seek necessary and sufficient conditions for the orbit matrix of $(G,p) \in \mathscr{R}_{(G,S,\Phi)}$ to be of full rank (i.e., for $(G,p)$ to have only trivial fully $(S,\Phi)$-symmetric infinitesimal motions) for a symmetry generic $p$, or to be independent (i.e., for $(G,p)$ to have no fully $(S,\Phi)$-symmetric self-stresses).  Of course, given a symmetric framework $(G,p) \in \mathscr{R}_{(G,S,\Phi)}$ which is independent and infinitesimally rigid with the usual rigidity matrix, its orbit matrix will also be independent and of maximal rank.  However, we have seen that there are frameworks which are dependent but the lack of a fully $(S,\Phi)$-symmetric self stress means that the orbit matrix is independent, as well as frameworks which have infinitesimal flexes but the lack of a fully $(S,\Phi)$-symmetric infinitesimal flex means that the orbit matrix is of full rank.  So we are seeking new results and will need new techniques.

The Fully Symmetric Maxwell's Rule ($r=c-m$) gives the standard necessary counts on $G$, $S$, and $\Phi$ for independence and full rank of an orbit matrix with $c$ columns, $r$ rows, and a space of trivial fully $(S,\Phi)$-symmetric infinitesimal motions (kernel of the orbit matrix for the complete graph) of dimension $m$.  As usual, there are some added necessary conditions for independence of the rows which come from subgraphs $G'$ of the graph $G$:
\begin{enumerate}
\item If the rows of the orbit matrix are independent, then for each fully symmetric subgraph $G'$ (generating $r'$  rows and $c'$ columns, as well as $m'$ trivial infinitesimal motions for these columns), we have $r' \leq c' - m'$;
\item If $H$ is a subgraph of $G$ such that $H$ and $x(H)$ are disjoint for each $x\in S$,   then $|E(H)| \leq d|V(H)| - {d+1 \choose 2}$, where the framework is in dimension $d$, with $|V(H)| \geq d$.
\end{enumerate}
Notice that we do not add special conditions for `small' subgraphs in part 1 above.  The reference to $m'$ actually codes for all those special cases.

How could we generate sufficient conditions?  One traditional way for the
standard
 rigidity matrix has been to start with minimal examples, and use inductive techniques which preserve the independence and full rank of the rigidity matrix.   These techniques include versions of vertex addition, edge splitting, and vertex splitting.   This has been extended to fully symmetric inductive techniques, still with the standard rigidity matrix, in \cite{BS4, BS3}.   Transferred to the orbit matrix, such fully symmetric techniques will still preserve the independence and the full rank of the orbit matrix.   However, there are many more inductive techniques which preserve the full rank of the orbit matrix - but would not preserve the full rank of the original rigidity matrix, since they would leave infinitesimal flexes which are not fully symmetric.   For example, simply adding a vertex along the axis of a 2-fold rotation in 3-space (which adds one column) will only require one added edge orbit - which could be one edge (along the axis) or two edges (the orbit of a single edge) and this would definitely not generate an infinitesimally rigid framework in 3-space!

It is unclear whether there are symmetry groups for which the full characterization is accessible.  When we find such a characterization, we will have a fully symmetrized version of the pebble game, for the orbit multi-graph.

\subsection{Geometric conditions for lower rank in the orbit matrix}
For standard rigidity, there has been an algebraic geometric exploration of when a specific configuration $p$ makes a generically rigid graph $G$ into an infinitesimally flexible framework  $(G,p)$.  The conditions are expressed in terms of a polynomial pure condition in the coordinates of $p$ which is $=0$ if and only if $(G,p)$ is infinitesimally flexible.  There will be a comparable theory for when configurations lower the symmetry generic rank of the orbit matrix.   We illustrate the layers of this for a specific plane example with $C_3$ symmetry.

\begin{examp}
Consider the framework illustrated in Figure~\ref{fig:Tri}(a).   The graph is generically rigid, and the pure condition for a lower rank of the rigidity matrix can be simplified to:  make any of the four triangles collinear or make the induced points  in Figure~\ref{fig:Tri}(c) collinear.
\begin{figure}
    \begin{center}
  \subfigure[] {\includegraphics [width=.20\textwidth]{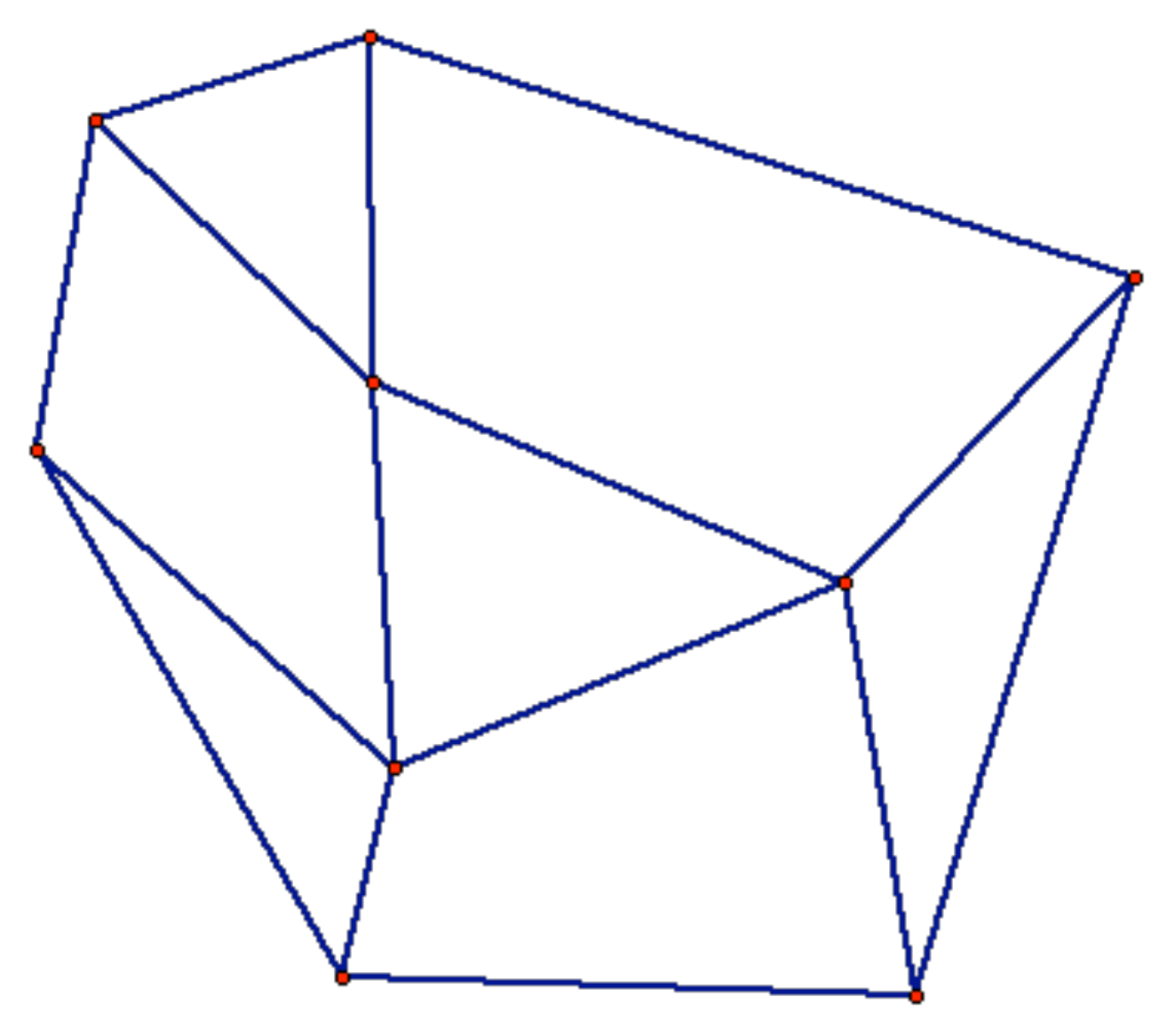}} \quad
  \subfigure[] {\includegraphics [width=.20\textwidth]{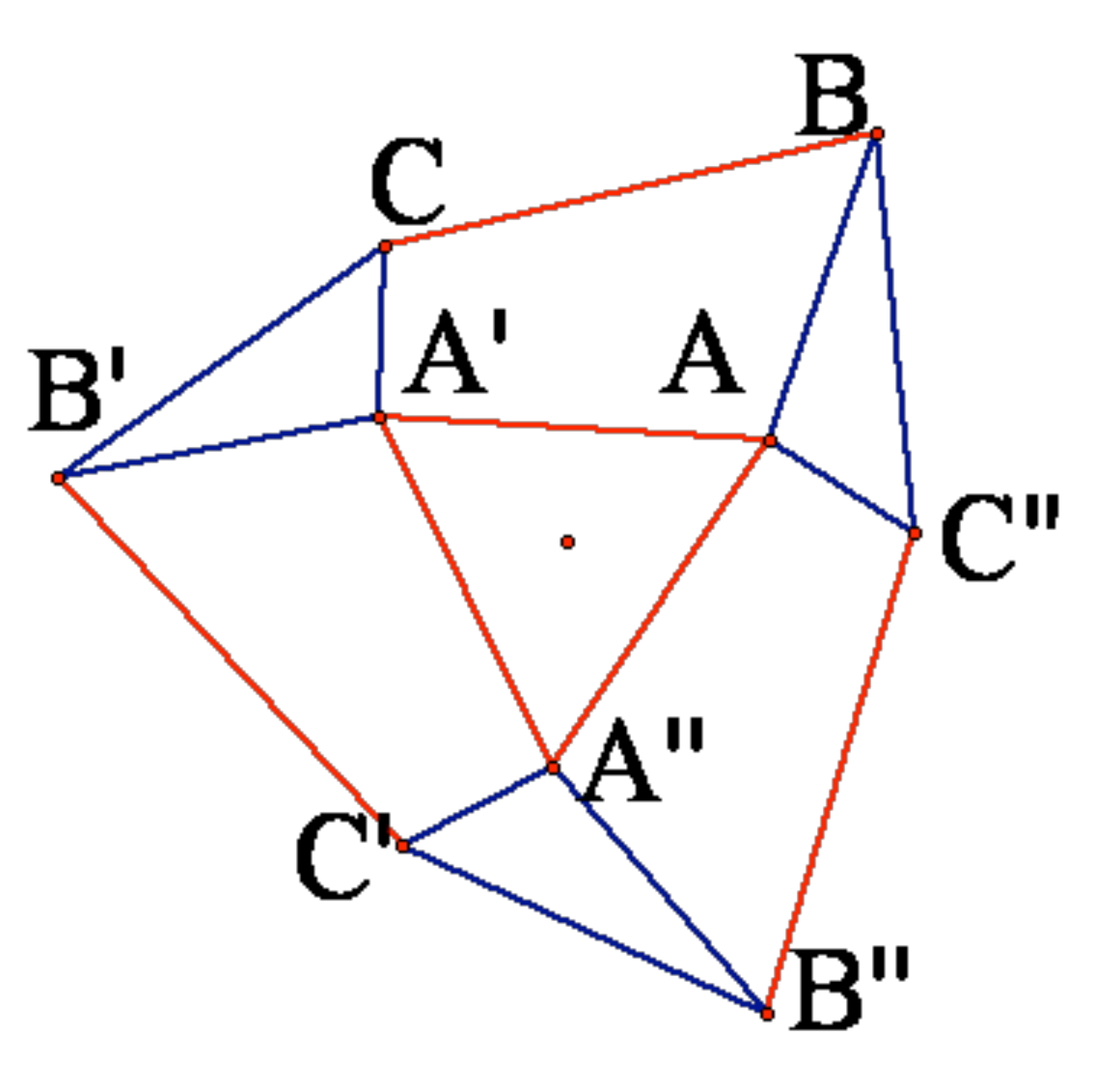}}
  \subfigure[] {\includegraphics [width=.35\textwidth]{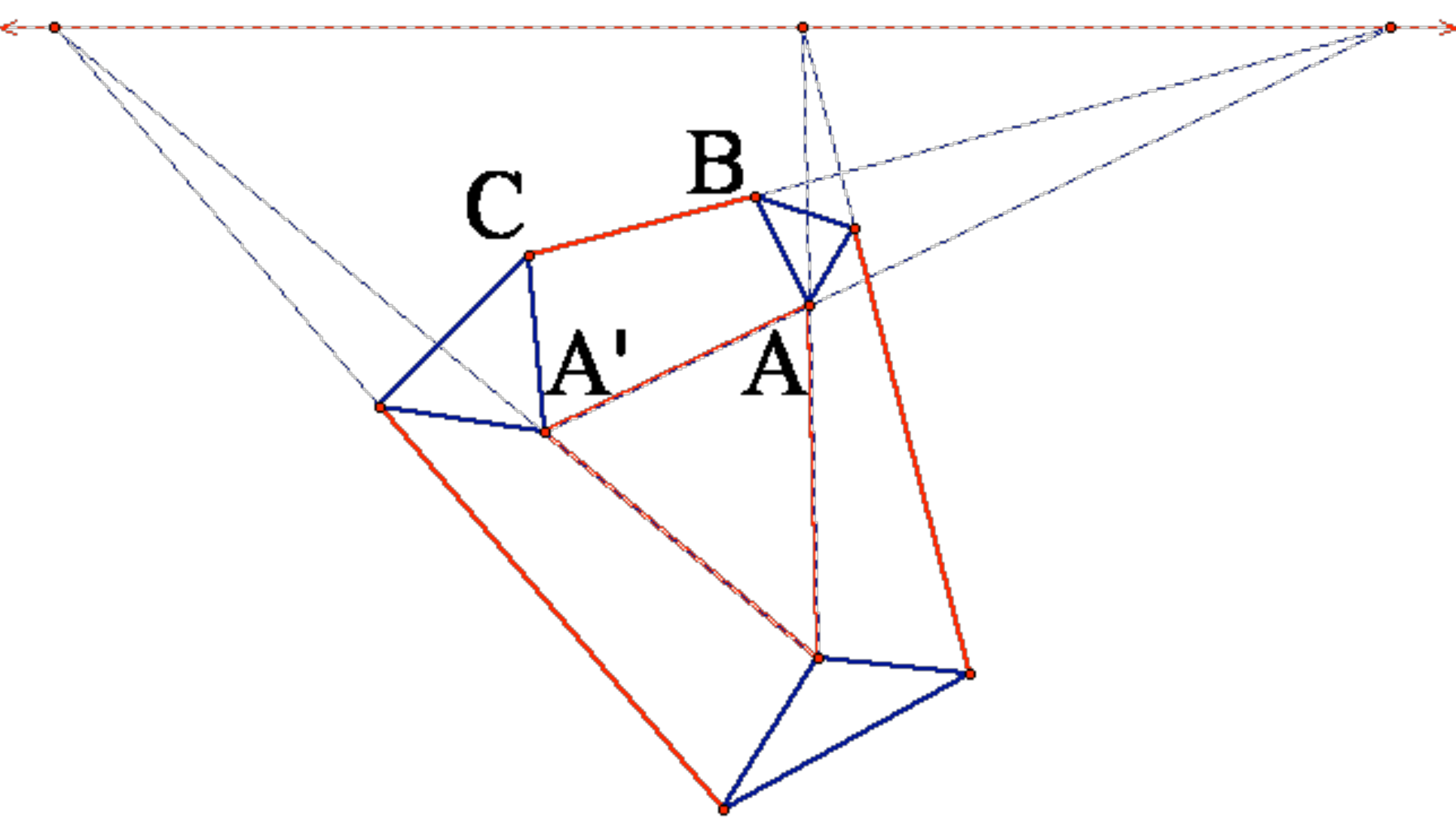}}
  \\
    \subfigure[] {\includegraphics [width=.25\textwidth]{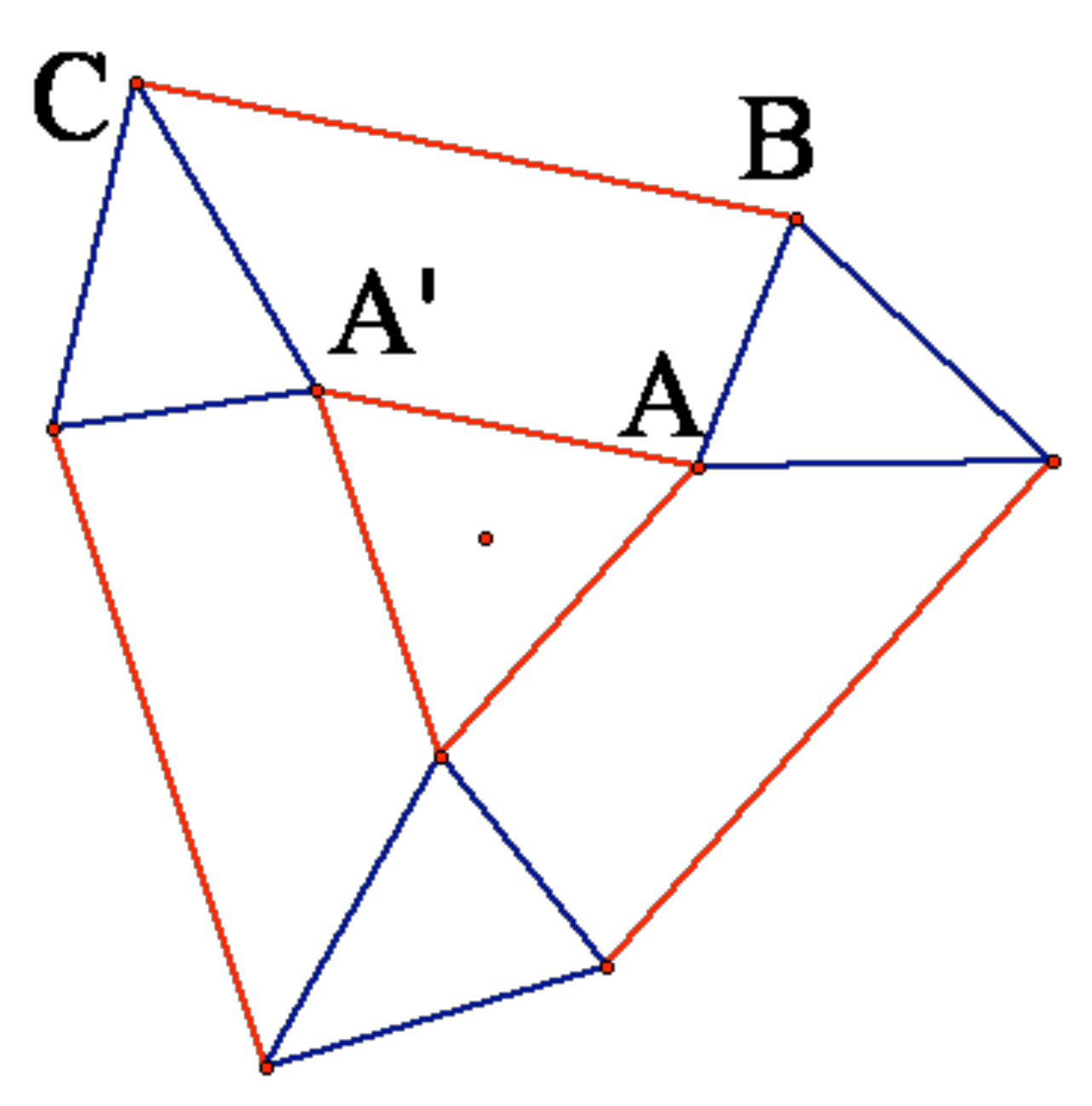}} \quad\quad
  \subfigure[] {\includegraphics [width=.28\textwidth]{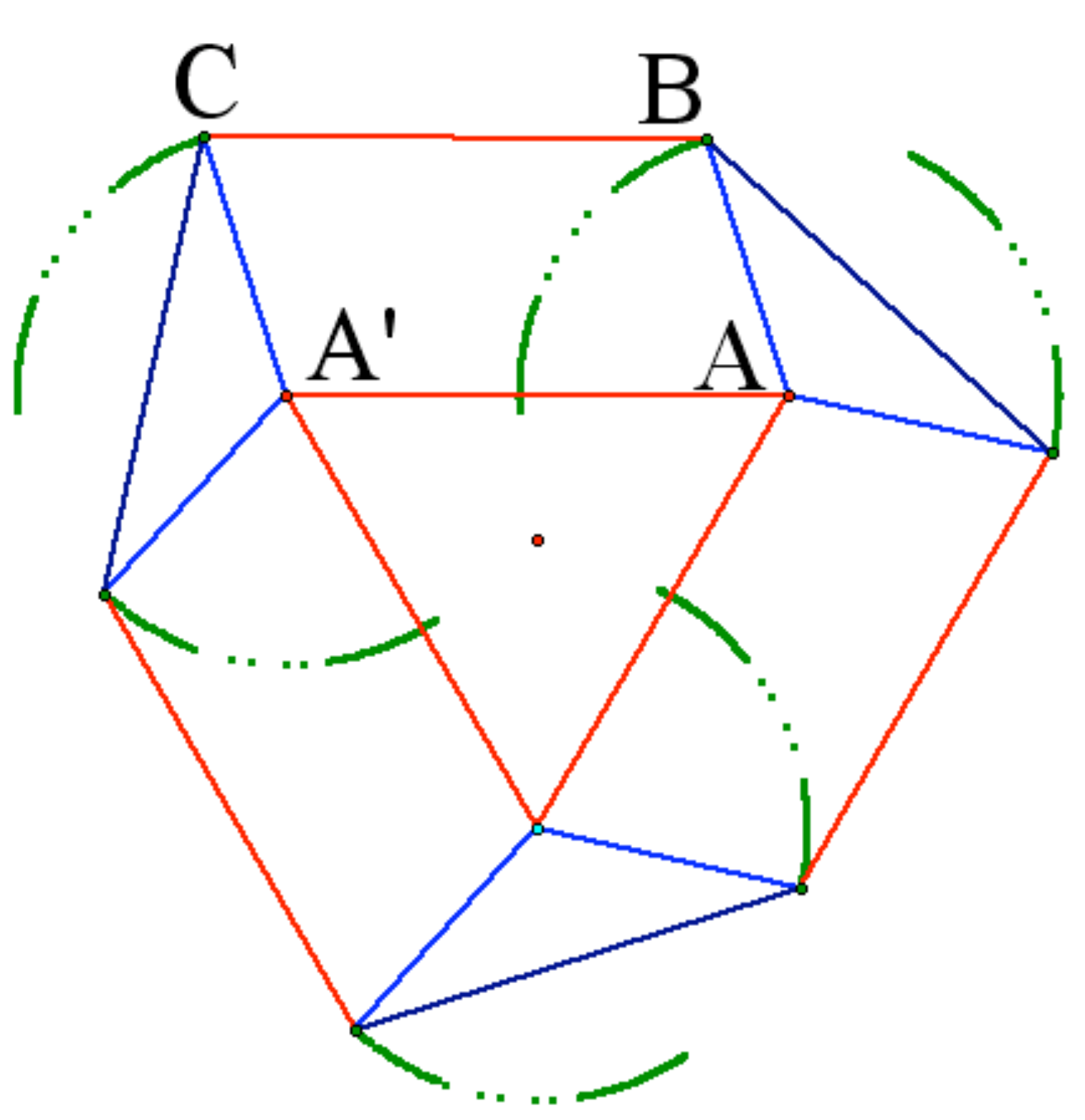}}
      \end{center}
    \caption{The graph in (a) is generically rigid, and is also symmetry generically rigid for 3-fold rotation (b).  The geometric condition for a non-trivial infinitesimal motion is three collinear induced points (c), which can also be achieved with 3-fold symmetry (d).  Only the configuration (e), with parallelograms $ABCA'$ have symmetry-preserving finite flexes.}
    \label{fig:Tri}
    \end{figure}

Symmetry generic realizations with $C_3$ symmetry are still infinitesimally rigid (Figure~\ref{fig:Tri}(b)).  Assuming  $C_3$ symmetry, the condition for an infinitesimal flex is that the three collinear points lie at infinity - or equivalently that the pairs of lines $AA', BC$ are parallel (Figure~\ref{fig:Tri}(d)).

This is not enough for a fully symmetric infinitesimal flex - or equivalently for a drop in the rank of the orbit matrix.  A direct geometric analysis verifies that the geometric condition for a fully symmetric infinitesimal flex (i.e. a drop in the rank of the orbit matrix) is that the three congruent faces $A,B,C,A'$ are parallelograms (Figure~\ref{fig:Tri}(e)).  From the geometric theory of such structures of parallelograms and triangles, it is known that this infinitesimal flex is a finite,  symmetry-preserving flex.  Thus, we can express the condition on a configuration lowering the rank of the orbit matrix in terms of a polynomial in the representative vertices $A,B,C$ and image of $A$ under $C_3$.
\end{examp}

This example suggests  that there is some interesting algebraic geometry to explore here.  In previous work \cite{WW} the polynomial conditions were extracted by using 'tie-downs' (equivalently striking out some columns) to square up the matrix.  This approach is still relevant, but rules for which tie-downs or pinned vertices remove all fully symmetric trivial motions are more complex.


\subsection{Transfer to other metrics}
\label{sec:metrictransf}
The paper \cite{SalWW} presents results about the transfer of first-order rigidity properties (essentially all properties of the rigidity matrix) among frameworks which realize a given graph, on the same projective configuration, in the metric spaces $\mathbb{E}^d$, $\mathbb{S}^d$ and $\mathbb{H}^d$.  What about a transfer of the orbit matrix for a symmetry group in $\mathbb{E}^d$ to the other metrics with the same symmetry groups?

In $\mathbb{E}^d$, all groups of isometries for a framework are point groups (there is a fixed point).  These other spaces also share these same point groups - a connection that can be seen by coning up a dimension and then slicing the cone along a corresponding unit sphere.  $\mathbb{S}^d$ and $\mathbb{H}^d$ have additional groups of isometries which do not fix a point and these can vary from space to space.

For simplicity, consider a point group in $\mathbb{E}^d$ and a sphere $\mathbb{S}^d$ tangent to the Euclidean space at the central point of the group.  It is not hard to give a correspondence to a point group in the spherical space as well as a correspondence between symmetry generic frameworks in the two spaces.  This correspondence will conserve fully symmetric infinitesimal flexes, fully symmetric trivial infinitesimal motions, and fully symmetric self-stresses.  In short, the orbit matrices of the two configurations in the two metrics will have a simple invertible correspondence generated by multiplication on the right and left by appropriate invertible matrices \cite{SchuWW}.

Underlying this transfer is the operation of symmetric coning - with a new vertex in the next dimension, which is on the normal to the lower dimension and extends the axes and mirrors in the lower space in a way that conserves the group, and preserves symmetry, including finite flexes.   A particular byproduct of this is the observation that repeated coning of the flexible octahedron or the flexible cross-polytope  will generate flexible polytopes in every dimension \cite{SchuWW}.

A similar process transfers orbit matrices  and the  predictions of finite flexes among
 $\mathbb{E}^d$, $\mathbb{S}^d$, and $\mathbb{H}^d$. This transfer gives a simple derivation of prior results on the flexibility of classes of Bricard
 octahedra and cross-polytopes in the spherical and hyperbolic metrics \cite{alex}.  It is unusual for flexibility to transfer - so symmetry is a special situation.  This transfer extends to other spaces with the same underlying projective geometry, such as the Minkowskian metric, provided that the point group is also realized as isometries in this metric.  The full exploration of this transfer is the topic of continuing exploration, and further details and results will be presented in \cite{SchuWW}.

 These other spaces such as $\mathbb{S}^d$  have additional symmetry groups which are not point groups (do not fix any point, or pair of antipodal points)
and hence do not correspond to the symmetries in  $\mathbb{E}^d$.
There will be orbit matrices for these groups as well, and hence we can study these cases using a direct extension of the methods presented in this paper.   These connections will be further explored in \cite{SchuWW}.

\subsection{Extensions to body-bar frameworks}
One now standard extension of bar and joint frameworks are the body-bar frameworks \cite{WWbb,CJW}.  These are a special class of frameworks, which in dimensions 3 and higher have a complete characterization for the multi-graphs which are generically isostatic (rigid, independent).  The basic analysis of symmetry adapted rigidity matrices for these structures has been presented in \cite{gsw}.

It is clear that there are corresponding orbit matrices for body-bar frameworks, since they have bar and joint models, and the desired orbit matrix can, in principle, be extracted from that.  The counting of columns and rows can also be adapted - though it would be helpful to give this in full detail.

A further extension studies body-hinge frameworks, with an emphasis on molecular models, where bodies (atoms) are connected by bonds (sets of 5 bars).  The molecular models also have bar and joint models, so in principle there are corresponding orbit matrices, and counts to predict finite flexes.  The classical `boat and chair' configurations of cyclohexane in chemistry  (a ring of six carbons) is an example where 3-fold symmetry (the chair) keeps the generic first-order rigidity and independence, and the 2-fold symmetry (the boat) is a model of the flexible octahedron.

Theorem~\ref{3dhalfturn} showing the flexibility of generically isostatic graphs in 3-space realized with 2-fold symmetry, extends from this example to general molecules in 3-space with 2-fold symmetry and no atoms or bonds intersecting the axis.  This is a common occurrence among dimers of proteins, so it has potential applications to the study of proteins \cite{W5}.

\subsection{Orbit matrices for other geometric constraint systems}
Owen and Power have investigated other examples of geometric constraints in CAD under symmetry \cite{owen}.   In general, constraint systems with matrix representations are open to analysis using group representations and symmetric block decompositions of their matrices.  However, there are some surprises which confirm that the analysis of corresponding orbit matrices may not be a simple translation of the results given here.

\begin{figure}
    \begin{center}
  \subfigure[] {\includegraphics [width=.30\textwidth]{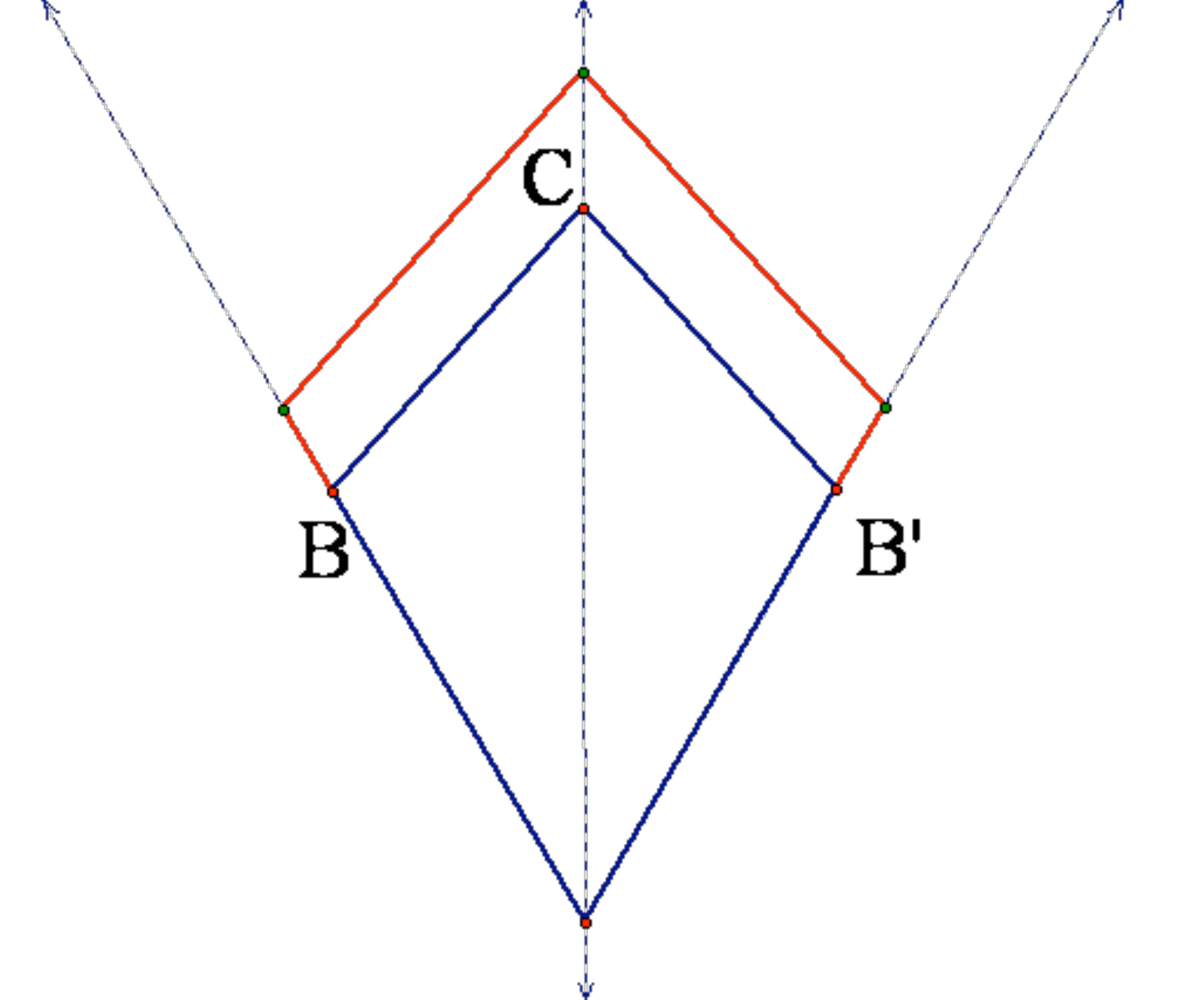}}
  \subfigure[] {\includegraphics [width=.40\textwidth]{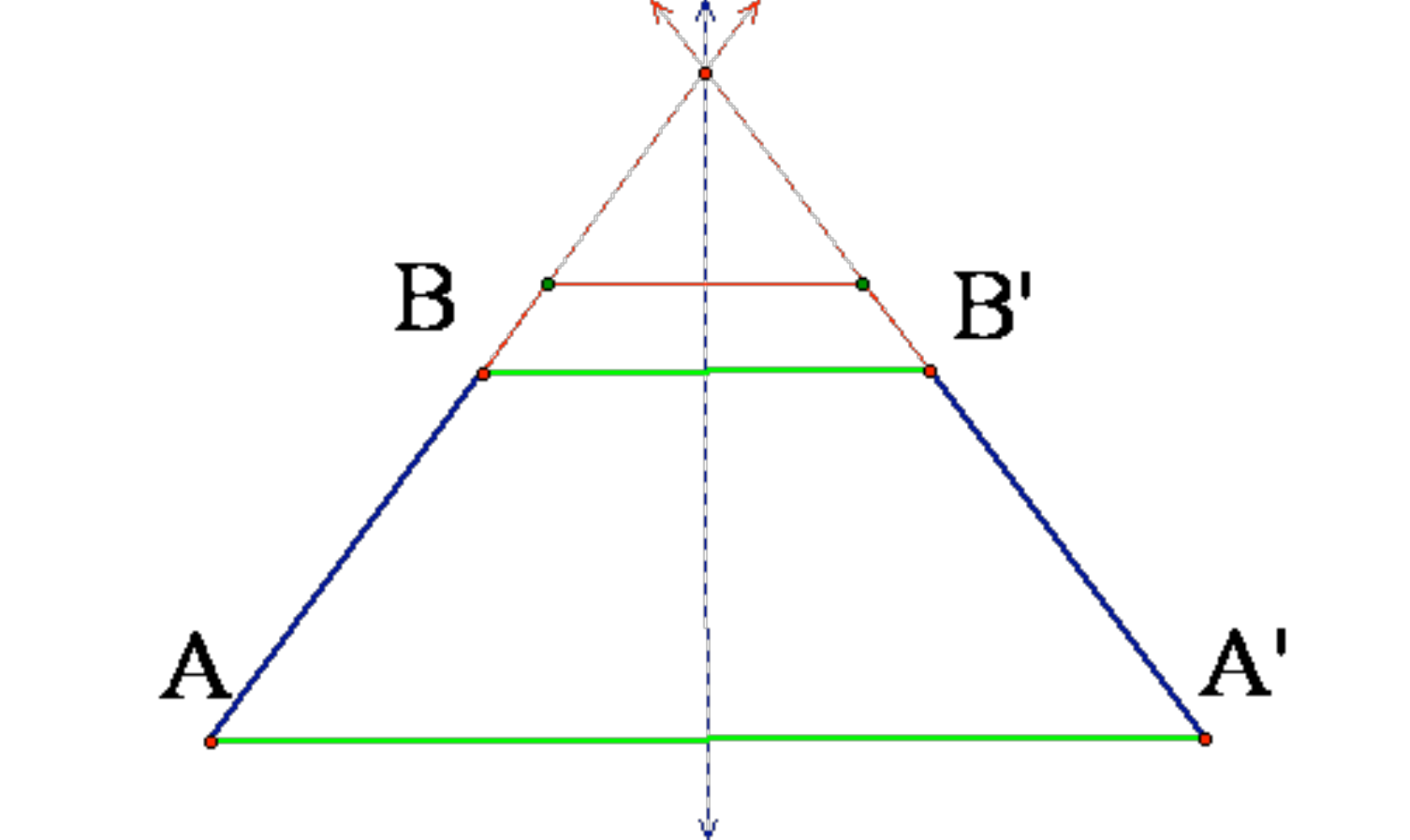}}
       \end{center}
    \caption{For symmetric frameworks, the space of fully  symmetric trivial parallel drawings may be larger than the space of fully  symmetric trivial infinitesimal motions (a), and frameworks without fully symmetric infinitesimal flexes may have fully symmetric parallel drawings (b).  }
    \label{fig:Parr}
    \end{figure}

It is well known that in the plane, infinitesimal motions correspond to parallel drawings of the same geometric graph and  configuration.  The correspondence involves turning all the velocities by $90^\circ$, which takes a trivial rotation to a trivial dilation.  For symmetry, this turn takes an infinitesimal motion which is fully symmetric for a rotation to a parallel drawing which is fully symmetric for the same rotation.  However, this operation takes an infinitesimal motion which is fully symmetric for a mirror to a parallel drawing which is anti-symmetric for the same mirror (and vice versa).  Clearly, there are changes in the development of the orbit matrix, even for this special example.  There are also additional fully symmetric trivial motions (dilations about the center of the point group are trivial, for the mirror).  More surprisingly, some edge orbits seem to disappear in the obit matrix (edge $AA'$ in Figure~\ref{fig:Parr}(b)).   Figure ~\ref{fig:Parr} illustrates two examples.
This may be enough to confirm that the extensions to other constraint systems are non-trivial, and worth carrying out!

\providecommand{\bysame}{\leavevmode\hbox to3em{\hrulefill}\thinspace}
\providecommand{\MR}{\relax\ifhmode\unskip\space\fi MR }
\providecommand{\MRhref}[2]{%
  \href{http://www.ams.org/mathscinet-getitem?mr=#1}{#2}
}
\providecommand{\href}[2]{#2}

\end{document}